\documentclass[11pt,a4]{article}

\usepackage[numbers]{natbib}

\usepackage{subfigure}
\RequirePackage{fix-cm}
 
\usepackage[a4paper]{geometry}

\usepackage{latexsym}
\usepackage{amsmath}
\usepackage{amssymb}
\usepackage{theorem}
\usepackage{xspace}
\usepackage{array}
\usepackage{subfigure}
\usepackage{epsfig}
\usepackage{hhline}
\usepackage{bbm}
\usepackage{color}
\usepackage{graphicx}
\usepackage{multirow}
\usepackage{dcolumn}
\usepackage{caption}

\usepackage[table]{xcolor}
\usepackage{colortbl}
\definecolor{Gray}{gray}{0.9}
\definecolor{Gray2}{gray}{0.85}

\def\cro#1{\left[#1\right]}

\def\eeX{\mathbb{X}}

\newcommand{\eX}{{\bf X}}

\def\eu{\mathbf{{u}}}

\newcommand{\N}{\mathbb{N}}
\newcommand{\R}{\mathbb{R}}
\newcommand{\PP}{\mathbb{P}}

\def\argmin{\mathop{\mathrm{arg\,min}}}

\def\ex{\mathbf{x}}

\def\ef{\textrm{\mathversion{bold}$\mathbf{\phi}$\mathversion{normal}}}  
\def\ee1{\textrm{\mathversion{bold}$\mathbf{\varepsilon}$\mathversion{normal}}}  
\def\eth{\textrm{\mathversion{bold}$\mathbf{\theta}$\mathversion{normal}}}

\def\eb{\textrm{\mathversion{bold}$\mathbf{\beta}$\mathversion{normal}}}

\def\eE{I\!\!E}

\def\e1{1\!\!1}

\theoremstyle{plain}

\newtheorem{theorem}{Theorem}[section]
\newtheorem{corollary}{Corollary}[section]
\newtheorem{lemma}{Lemma}[section]
\newtheorem{remark}{Remark}[section]

\definecolor{gray2}{rgb}{0.4,0.4,0.4}

\usepackage[title]{appendix}

\renewcommand{\eqref}[1]{(\ref{#1})}

\begin{document}

 \title{Change-point detection in a linear model by adaptive fused  quantile  method}

 \author{Gabriela CIUPERCA$^1$ and Mat\'u\v{s} MACIAK$^{2}$}
 
 \maketitle
 \setcounter{footnote}{1}
\footnotetext{\noindent \small{Universit\'e de Lyon, Université Lyon 1, CNRS, UMR 5208, Institut Camille Jordan, Bat.  Braconnier, 43, blvd du 11 novembre 1918, F - 69622 Villeurbanne Cedex, France\\
	\textit{Email address}: Gabriela.Ciuperca@univ-lyon1.fr\\
	\indent $^2$Charles University, Faculty of Mathematics and Physics, Department of Probability and Mathematical Statistics, Sokolovsk\'a 83, Prague, 186 75, Czech Republic\\
	\textit{Email address}: Matus.Maciak@mff.cuni.cz}}


 
\begin{abstract}
	A novel approach to quantile estimation in multivariate linear regression models with change-points is proposed: the change-point detection and the model estimation are both performed automatically, by adopting either the quantile fused penalty or the adaptive version of the quantile fused penalty. These two methods combine the idea of the check function used for the quantile estimation and the $L_1$ penalization principle known from the signal processing and, unlike some standard approaches, the presented methods go beyond typical assumptions usually required for the model errors, such as sub-Gaussian or Normal distribution. They can effectively handle heavy-tailed random error distributions, and, in general, they offer a more complex view on the data as one can obtain any conditional quantile of the target distribution,
not just the conditional mean. The consistency of detection is proved and proper convergence rates for the parameter estimates are derived. The empirical performance is investigated via an extensive comparative simulation study and practical utilization is demonstrated using a real data example.
\end{abstract}


\section{Introduction}
In this paper we consider an automatic detection of change-points in a multivariate linear model using the fused penalty  technique. The proposed method  covers a large spectrum of scenarios for various model errors and, above all, it allows to simultaneously detect the number of change-points and their locations in the underlying model. This avoids using  multistage procedures where firstly one needs to detect the number of change-points on a basis of some criterion, find their locations, and to estimate the corresponding model parameters (see, for instance, \cite{Ciuperca-16}).

The change-point detection assessed by penalizing the sum of squares  with the fused group LASSO penalty was firstly considered by \cite{Zhang-Xiang-15} for detecting one change-point with the Gaussian error terms. The detection of several change-points was later considered in \cite{Zhang.Geng.15}.
In general, the least squares regression with the LASSO type penalties used for detecting multiple change-points is proposed in \cite{Jin-Wu-Shi-16} or \cite{Hyun-GSell-Tibshirani-18}. Independent and identically distributed ($iid$) centered errors with some bounded variance and the LASSO-type penalty are discussed in \cite{Leonardi-Buhlmann.16} and the idea is further elaborated for the fused penalty and strong mixing centered errors under some specific moment conditions in \cite{Qian.Su.16}.

 For a particular case of a piecewise constant model with the Gaussian or sub-Gaussian errors, there is \cite{Rinaldo.09} who penalizes the $L_2$-norm with the fused LASSO penalty. Under some more general assumptions of $iid$ zero-mean and bounded variance errors, \cite{Harchaoui.Levy.10} propose a~method to automatically detect the number of change-points and their locations by the fused penalty  method. The results have been further deepened in \cite{lin2016}. On the other hand, if the model error terms do not satisfy some standard conditions, the penalized  least squares methods are no longer applicable. An alternative approach for this case is considered, for instance, in \cite{Ciuperca-Maciak-17} where the authors proposed an idea of the quantile LASSO instead. The quantile linear regression with the LASSO method for detecting a change is also investigated in \cite{Lee-Liao-Seo-Shin-18} but the authors only focus on situations where one change-point occurs in the model. 

However, in contrast to  \cite{Ciuperca-Maciak-17}, where the automatic detection of change-points in a piece-wise constant model is studied for the fused quantile penalty, in the present work we consider a multiple regression model with multiple change-points and the adaptive fused penalty is used instead to recover these change-points. The adaptive penalty improves the performance of the change-point detection and decreases the shrinkage and, therefore, it can applied for an automatic detection and simultaneous estimation in general linear models with a fixed number of explanatory variables.

The rest of the paper is organized as follows: in the next section we introduce some notation and the underlying model is defined together with some necessary assumptions.
In Section \ref{section3} we introduce the change-point detection approach based on the quantile fused penalty and some theoretical results are derived. An adaptive fused quantile method is discussed in Section \ref{section4}, the consistency of the detection is proved and the proper convergence rate for the parameter estimates are given. Both methods are investigated in terms of an extensive comparative simulation study in Section \ref{section5} and a real data example is presented in Section \ref{regression_example}. All technical details and  theoretical proofs are given in the appendix section.

\section{Model and assumptions}
Let us start by introducing some necessary notation. We use $C$ to denote a  positive generic constant which does not depend on $n \in \mathbb{N}$. Moreover, for any set of elements $E$, we also denote its complement by  $\overline E$. For any vector, we use $\|.\|$ to denote the Euclidean norm and $\|.\|_\infty$ to denote the maximum norm. Similarly, for some matrix, we use $\|.\|$ to denote its  spectral norm. Moreover, for a positive definite matrix, we  use $\mu_{min}(.)$ (and $\mu_{max}(.)$ respectively)  to denote its largest (or smallest respectively) eigenvalue. Finally, for some positive sequences $(s_n)_n$ , $(r_n)_n$ we denote by $s_n \gg r_n$ the fact that $\lim_{n \rightarrow \infty} s_n/r_n=\infty$ and for any real number $x \in \mathbb{R}$, we use $[x]$ to denote its integer part.

Consider now a linear model for which the parameters can change along observations $i \in \{1, \dots, n\}$, such that
\begin{equation}
\label{eq1}
Y_i =\ex_i^\top \eb_i+\varepsilon_i, \qquad i=1, \cdots , n,
\end{equation}
where  $\eb_i \in \R^p$, with $p$ not depending on $n \in \mathbb{N}$, and $\ex_i=(c^{(0)}, x_{2i}, \cdots , x_{pi})^\top$, for $c^{(0)}$ being some nonzero constant, is the vector of the subject's specific explanatory variables for some observation (index) $i \in \{1, \dots , n\}$. In other words, the model defined by \eqref{eq1} is assumed to contain the intercept term by default. This assumption is, however, standard for quantile models in high-dimension in general (see, for instance, \cite{Zheng-Gallagher-Kulasekera-13}, \cite{Zheng-Peng-He-15}). In addition, the model in \eqref{eq1} is assumed to have $K^* \in \mathbb{N}$ changes located at $t^*_1 < \cdots < t^*_{K^*} \in \{1, \dots, n\}$, such that 
\begin{equation}
\label{eq2}
\eb_i=\eb_{t_k}, \qquad \forall i=t^*_k, t_k^*+1, \cdots , t^*_{k+1}-1, \qquad k=0,1, \cdots , K^*, 
\end{equation}
with $t^*_0=1$, $t^*_{K^*+1}=n$, and $\eb_n=\eb_{t^*_{K^*+1}}$. For simplicity we can define an overall $np$-dimensional vector of parameters $\eb^n=(\eb_1^\top, \cdots, \eb_n^\top)^\top \in \mathbb{R}^{np}$. In general, the number of change-points $K^* \in \mathbb{N}$ and the locations $t^*_1, \cdots , t^*_{K^*}$ where change-points occur, are all left unknown.  Alternatively, for each model phase $k=0, \cdots, K^*$, we have the corresponding vector parameters $\ef_{1}, \dots, \ef_{K^* + 1} \in \mathbb{R}^p$, where
\[
\ef_{k+1} \equiv  \eb_{i}, \qquad \textrm{for} \quad  i=t^*_k, t^*_k+1, \cdots , t^*_{k+1}-1, \quad \textrm{and} \quad k=0, \cdots, K^*,
\]
and, analogously, $\ef_{K^*+1} = \eb_n$ for the last phase. The true values of parameters $\ef_k$, for $k=1, \cdots K^*+1$, are also unknown and they are denoted by  $\ef^*_k$. The corresponding true values of the vector parameters $\eb_i$ are denoted by $\eb^*_i$, for $i=1, \cdots , n$. 
It is assumed that  the number of true change-points $K^* \equiv Card\{ i \in \{2, \cdots , n \}; \; \eb^*_i \neq \eb^*_{i-1}  \}$ is bounded, but it is unknown. 

For the piece-wise constant model considered in \cite{Ciuperca-Maciak-17}, there is $p=1$ and, hence, the dimension of $\eb$ is same as the overall number of observations. However, in a general case where $p \in \mathbb{N}$, the number of parameters can heavily exceed the number of observations.
Let us define the set  ${\cal A}^*$ which  contains all indexes (locations) of the true change-points
\[
{\cal A}^*=\{t^*_1, \cdots , t^*_{K^*}  \}.
\]
Obviously, it holds that $|{\cal A}^*|=K^*$. Let us also define the empirical quantile process 
\begin{equation}
\label{eqQ}
G_n(\eb^n)=\sum^n_{i=1} \rho_\tau(Y_i- \ex_i^\top \eb_i),
\end{equation} 
which is associated with the model in \eqref{eq1}. Function $\rho_\tau :\R \rightarrow \R_+$ is used to denote the standard check function $\rho_\tau(u)=u(\tau -\e1_{\{u< 0\}})$, for some fixed quantile level $\tau \in (0,1)$, and any $u \in \R$. 

The model above can be also equivalently expressed by using consecutive differences between the unknown parameters $\eb_1, \dots, \eb_n$. Let us define parameters $\eth_1, \dots, \eth_n \in \mathbb{R}^p$, where $\eth_1=\eb_1$ and $\eth_j=\eb_j- \eb_{j-1}$, for any $j = 2, \cdots , n$. These vector parameters can be again joined into just one overall vector   $\eth^n=(\eth_1^\top, \cdots, \eth_n^\top)^\top \in \mathbb{R}^{np}$,  with the true values  $\eth^*=(\eth^{*\top}_1, \cdots , \eth^{*\top}_n)^\top$ being associated with $\eb_1^*, \cdots , \eb_n^*$ and their consecutive differences in particular. Hence, by  using the model matrices 
\[
\eX = \left[  
\begin{array}{cccc}
\ex_1 & \textbf{0} & \cdots & \textbf{0 } \\
\textbf{0} & \ex_2 &  \cdots & \textbf{0 } \\
\vdots& \vdots & \ddots & \vdots \\
\textbf{0 } & \textbf{0 } & \ldots & \ex_n
\end{array}
\right], 
\qquad 
\textbf{A}=\left[  
\begin{array}{cccc}
\textbf{I}_p & \textbf{0} & \cdots & \textbf{0 } \\
\textbf{I}_p & \textbf{I}_p &  \cdots & \textbf{0 } \\
\vdots& \vdots & \ddots & \vdots \\
\textbf{I}_p & \textbf{I}_p & \ldots & \textbf{I}_p
\end{array}
\right],
\qquad \textrm{and} \quad 
\eeX=\eX \textbf{A},
\] 
of the  dimensions  $n \times (np)$, $(np) \times (np)$,  and $n \times (np)$, we can  rewrite the model in \eqref{eq1} in an equivalent form 
\begin{equation}
\label{eq5}
\textbf{Y}^n=\eX \eb^n+\ee1^n=\eeX \eth^n+\ee1^n,
\end{equation}  
or, alternatively, also in a cumulative form
\begin{equation}
\label{eq5bis}
Y_i=\ex_i^\top \sum^i_{s=1} \eth_s+\varepsilon_i, \qquad i=1, \cdots , n.
\end{equation}

\noindent\textbf{Assumptions:}
\begin{itemize}
	\item[\textbf{(A1)}] There exists (eventually after a change of scale) a constant $0< c^{(1)}< 1$ not depending on $n$ such that  $\max_{1 \leqslant i \leqslant n} \| \ex_i\|  \leq c^{(1)}$.
	
	\item[\textbf{(A2)}] The random error  terms $(\varepsilon_i)$ are $iid$ with the distribution function $F(x) >0$, for all $x \in \R$, such that $\PP[\varepsilon <0]=\tau$.  Moreover, the corresponding density function $f(x)$ is bounded;
	
	\item[\textbf{(A3)}] There exist two constants $0 < m_0 \leq M_0 < \infty$, such that  
	\begin{align*}
	m_0 &  \leq \inf_{1 \leq n_1 < n_2 \leq n}  \mu_{\min} \Big( (n_2-n_1)^{-1} \sum^{n_2-1}_{i=n_1} \ex_{i} \ex_{i}^\top\Big)\\
	&  \leq \sup_{1 \leq n_1 < n_2 \leq n}  \mu_{\max} \Big( (n_2-n_1)^{-1} \sum^{n_2-1}_{i=n_1} \ex_{i} \ex_{i}^\top\Big) \leq M_0,
	\end{align*} 
	for any $n_1, n_2 \in \mathbb{N}$, such that  $1 \leqslant n_1 < n_2 \leqslant n$. Moreover,  the minimal distance between two consecutive change-points is $I^*_{min}\equiv \min_{1 \leqslant k \leqslant K^*}(t^*_{k-1} - t^*_k)$ and, analogously $I^*_{max}\equiv \max_{1 \leqslant k \leqslant K^*}(t^*_{k-1} - t^*_k)$ for the maximum distance.
	
	\item[\textbf{(A4)}] Let $I^*_{min} \geq n \delta_n$, for some decreasing sequence $(\delta_n)$, such that $\delta_n \rightarrow 0$ and also $(\log n)^{-1}n \delta_n \rightarrow \infty$, for $n \rightarrow \infty$. \\

	\item[\textbf{(A5)}] There exist two bounded constant $0 < c^{(a)}, c^{(b)} < \infty$, not depending on $n$, such that 
	\begin{itemize}
		\item[(a)] $\max_{1 \leqslant j \leqslant K^*}\| \ef^*_{j+1} - \ef^*_j \| <c^{(a)}$;
		\item[(b)] $\min_{1 \leqslant j \leqslant K^*}\| \ef^*_{j+1} - \ef^*_j \| > c^{(b)}$.
	\end{itemize}
	
	\item[\textbf{(A6)}] The overall number of change-points in the model,  $K^* \in \mathbb{N}$, is bounded and does not depend on $n \in \mathbb{N}$.	
\end{itemize}

Since $c^{(0)}\neq 0$, Assumption (A1) implies that $0< |c^{(0)}|<1$ and it is required to control the quantity $x_{ij}\e1_{\{\varepsilon \leq u\}}$ for $u \in \R$ (see, for instance, \cite{Leonardi-Buhlmann.16}).  Assumption (A2) is standard for the quantile regression models
and the independence of the error terms is commonly considered also in various change-points models in \cite{Lee-Liao-Seo-Shin-18}, \cite{Jin-Wu-Shi-16}, \cite{Zhang.Geng.15}, \cite{Ciuperca-16}, and others. Assumption (A3) imposes restrictions on the eigenvalues of design matrix such that  the matrix is well defined (see also \cite{Qian.Su.16} or  \cite{Zhang.Geng.15}).  Assumption (A4) is common for ensuring a proper change-point detection by the LASSO methods (see \cite{Harchaoui.Levy.10} or  \cite{Ciuperca-Maciak-17}) and it postulates that the true change-points are at a mutual distance which is big enough where the sequence $(\delta_n)$ controls the convergence rate of the change-point estimator when the number of changes is correctly estimated. Assumption (A5) is necessary  to distinguish the existing change-points in the model and Assumption (A6), which is also considered in \cite{Rinaldo.09} or \cite{Zhang-Xiang-15} for Gaussian errors, is needed for the detection  of multiple change-points in the linear regression model using the least squares and the fused group LASSO penalty.

\section{Change-point detection by the quantile fused method}
\label{section3}
In this section we firstly propose the quantile fused estimation approach and  we study the properties of the obtained estimates:  the estimates of the change-point locations and the estimates of the corresponding regression parameters between two consecutive change-points. These estimators will used later, in the next section, to define 
the weights for the adaptive fused quantile approach which can provide better asymptotic and finite sample results. Considering the model in \eqref{eq1}, the unknown parameter vector $\eb^n \in \R^{np}$ is estimated by minimizing the objective function
\begin{equation}
\label{eq3}
R(\eb^n)=\sum^n_{i=1} \rho_\tau(Y_i- \ex_i^\top \eb_i)+n\lambda_n \sum^n_{i=2} \|\eb_i-\eb_{i-1}\|,
\end{equation}
with the fused group LASSO type penalty. The value of the tuning parameter $\lambda_n  > 0$  controls the number of changes appearing in the final model: for $\lambda_n \to 0$ there will be a change-point detected at each available observation while the scenario with $\lambda_n \to \infty$ will result in a simple ordinary linear regression fit with no change-points at all.

In addition, let us assume that the tuning parameter $\lambda_n$ converges to zero with a slower rate than the sequence $(\delta_n)$. Let the following holds:
\begin{itemize}
	\item[\textbf{(A7)}] Let $(\delta_n)$ and $(\lambda_n)$ be two positive sequences satisfying the following: $n \lambda_n \rightarrow \infty$ and  $ \lambda_n / \delta_{n} \rightarrow 0$, as $n \to \infty$.
\end{itemize}
One possible option how to choose sequences $(\lambda_n)$ and $(\delta_n)$ such that they will satisfy Assumptions (A4) and (A7) is, for instance, to take $\lambda_n=n^{-1} (\log n)^{5/2} $ and $\delta_n=n^{-1} (\log n)^3 $. It is straightforward to see that instead of minimizing the objective function in \eqref{eq3} with respect to  $\eb^n \in \mathbb{R}^{np}$, one can equivalently deal with the objective function
\begin{equation}
\label{eq4}
\widetilde R(\eth^n) \equiv \sum^n_{i=1} \rho_\tau\bigg(Y_i - \ex_i^\top \sum^i_{k=1} \eth_k\bigg) +n \lambda_n \sum^n_{i=2} \| \eth_i\|,
\end{equation}
where the minimization now takes place with respect to $\eth^n \in \mathbb{R}^{np}$. The formulation in \eqref{eq4} can be also seen in terms of the  quantile LASSO problem with grouped variables: the number of groups is (up to one) equal to the number of observations $n \in \mathbb{N}$ and the number of parameters in each group is $p$. Thus, we can define the quantile fused estimators for the unknown parameters $\eb^n$, and  $\eth^n$ respectively, such that 
\[
\overset{\vee}{ \eb^n } \equiv  \argmin_{\eb^n \in \R^{np}} R(\eb^n), \qquad \overset{\vee}{ \eth^n} \equiv \argmin_{\eth^n \in \R^{np}}  \widetilde R(\eth^n),
\]
where $\overset{\vee}{ \eb^n} =(\overset{\vee}{\eb_1^\top}, \cdots , \overset{\vee}{\eb_n^\top})^\top$ and $\overset{\vee}{ \eth^n}=(\overset{\vee}{ \eth_1^\top}, \cdots , \overset{\vee}{ \eth_n^\top})^\top$. 
The estimators of the change-point  locations are the  observation indexes $i \in \{2, \cdots , n \}$, where $\overset{\vee}{ \eb}_i \neq \overset{\vee}{ \eb}_{i-1} $, that is $ \overset{\vee}{\eth}_i \neq \textbf{0}$. Let us define the set of estimated change-point locations as
\begin{equation}
\label{hAn}
\overset{\vee}{ {\cal A}}_n \equiv \{ i \in \{2, \cdots , n \}; \; \overset{\vee}{ \eb}_i \neq \overset{\vee}{ \eb}_{i-1}  \} = \{ \overset{\vee}{  t}_1 < \cdots < \overset{\vee}{  t}_{|\overset{\vee}{ \cal A}_n|} \},
\end{equation}
where $|\overset{\vee}{ \cal A}_n|$ denotes the cardinality of $\overset{\vee}{ {\cal A}}_n$. Now, for any $k=0, \cdots , |\overset{\vee}{ {\cal A}}_n|$, thus for any $(k+1)$-th model phase (i.e., observations starting with $\overset{\vee}{  t}_k$ until $\big(\overset{\vee}{  t}_{k+1}-1 \big)$, for $\overset{\vee}{t}_0=1$ and $\overset{\vee}{  t}_{|\overset{\vee}{ {\cal A}}_n|+1}=n$), we have the corresponding quantile fused  parameter estimator  $\overset{\vee}{ \ef}_{k+1} \in \mathbb{R}^p$. In this section we  study the consistency properties of the change-point location estimators in \eqref{hAn} and the corresponding convergence rate of the given regression parameter estimators. We also show that the proposed method overfits the true model (with probability converging to one) with respect to the number of change-points being detected. However, the proofs of the results from this section are omitted because they follow in a straightforward way from the proofs of the next section where the weights are all set to one. On the other hand, for some identifiability purposes, we impose the following  assumption on the distribution function of the error terms:
\begin{itemize}
	\item[\textbf{(A8)}] For any $k=1, \cdots , K^*$, the limit  $\lim_{n \rightarrow \infty} (n \delta_n)^{-1} \sum^{t^*_k -1}_{i=t^*_k-[n \delta_n]} F\big(\ex_i^\top(\ef^*_{k+1} -\ef^*_k)\big)\equiv L_k$ exists and $\tau \neq L_k$.
\end{itemize}
In case of $p=1$, Assumption (A8) becomes $F\big(\ef^*_{k+1} -\ef^*_k\big) \neq F(0)$, which implies, taking into account that $F(x)>0$ from Assumption (A2), that $\ef^*_{k+1} \neq \ef^*_k$, which is a classical condition required for the model between two successive change-points.

The following theorem shows that if the number of estimated  change-points in the model is equal to $K^*$, then for each true change-point location $t_k^* \in \{1, \dots, n\}$ we have a corresponding estimator $\overset{\vee}{ t}_k \in \{1, \dots, n\}$, such that their mutual distance is less then $n \delta_n$ with probability converging to one. 

\begin{theorem}
	\label{Theorem 3.1}
	Under Assumptions (A1)-(A8), if $|\overset{\vee}{ \cal A}_n|=K^*$,  then it holds that
	\[
	\lim_{n  \rightarrow \infty } \PP \left[\max_{1 \leq k \leq K^*} |\overset{\vee}{ t}_k - t^*_k | \leq n \delta_n  \right]  =  1.
	\]
\end{theorem}  

In addition, the next theorem provides the corresponding convergence rates  for the regression coefficient estimates, however, under the situation where  $|\overset{\vee}{ \cal A}_n |=K^*$.  The  convergence rates 
depend on the value of the regularization parameter $\lambda_n > 0$ and the minimal distance between two consecutive   change-points in the true model. Therefore, for convenience, we define the sequence 
\[b_n=n \lambda_n (I^*_{min})^{-1}+ \big(I^*_{min}\big)^{-1/2},\]
for any $n \in \mathbb{N}$. The convergence rate for the parameters estimates of the regression coefficients in the true model is given by the following theorem.

\begin{theorem}
	\label{Theorem 3.1ii}
	Under the same assumptions as in Theorem  \ref{Theorem 3.1}, if, in addition, it holds that $(n \delta_n)^{-1}I^*_{min} \rightarrow \infty$, as $n \rightarrow \infty$, then 
	\[
	\big\|\overset{\vee}{ \ef}_k-\ef^*_k \big\|=O_{\PP} (b_n),
	\]
	for any $k = 1, \cdots , K^*+1$.
\end{theorem}

Let us note, that the condition $(n \delta_n)^{-1}I^*_{min} \rightarrow \infty$, for $n \rightarrow \infty$, required in Theorem \ref{Theorem 3.1ii}, is necessary to separate the estimators from two consecutive change-points.

Next, we deal with the overestimation case: if the number of estimated change-points is strictly greater than $K^*$, we suppose that it is still inferior to an arbitrary number $K_{max}$, but bounded.
Thus, we only consider cases with a bounded set $|\overset{\vee}{ {\cal A}}_n|$.
Theorem \ref{Theorem 3.1} is, therefore, a special case of Theorem \ref{Theorem 3.2}. Let us consider the distance ${\cal E}$ between two  sets $A$ and $B$ defined as
\[
{\cal E}(A|| B) \equiv \sup_{b \in B} \inf_{a \in A} |a-b|,
\]
which is the analogy of the set distance used in \cite{Harchaoui.Levy.10}.

\begin{theorem}
	\label{Theorem 3.2}
	Under Assumptions (A1)-(A8), for $K^*=  |{\cal A}^*| \leq |\overset{\vee}{ {\cal A}}_n| \leq K_{max} < \infty $,   we have that
	$$
	\lim_{n \rightarrow \infty}  \PP \big[{\cal E}(\overset{\vee}{ {\cal A}}_n|| {\cal A}^*) \leq n \delta_n \big] = 1.
	$$
\end{theorem}

To summarize the results above, the proposed fused penalty method can be effectively used to detect all existing change-points in the model if there are at least  as many change-point detected as the number of true change-points $K^* \in \mathbb{N}$. By the following theorem we show that the last scenario where the method underestimates the number of change-points, only occurs with probability converging to zero as $n$ tends to infinity.

\begin{theorem}
	\label{Theorem 3.3}
	Under Assumptions (A1)-(A8), we have that
	$$\PP[|\overset{\vee}{ \cal A}_n | < K^*] {\underset{n \rightarrow \infty}{\longrightarrow}} 0.$$
\end{theorem}

By Theorems  \ref{Theorem 3.2} and \ref{Theorem 3.3}, we deduce that for each $k=1, \cdots , K^*$, the true change-point location $t^*_k$ has (with probability tending to one) at least one estimator $\overset{\vee}{ t}_j$, $j \in \{1, \cdots ,  |\overset{\vee}{ \cal A}_n |\}$ at a distance less than $[n \delta_n]$.   On the other hand, considering the convergence rate of the regression parameter estimates between two consecutive change-points obtained in the proof of Theorem \ref{Theorem 3.2}, we can only consider the elements  of $\overset{\vee}{ \cal A}_n $, for which the mutual distance converges to infinity as $n \rightarrow \infty$. Hence, instead of $\overset{\vee}{ \cal A}_n $ we can consider a smaller set 
\[
\overset{\smile}{ {\cal A}_n}=\left\{\overset{\vee}{t}_j \in \overset{\vee}{ \cal A}_n,~j \in \{1, \cdots , |\overset{\vee}{ \cal A}_n| \};~ \overset{\vee}{t}_j - \overset{\vee}{ t}_{j-1}  \overset{\PP} {\underset{n \rightarrow \infty}{\longrightarrow}} \infty \right\}.\]

We denote the elements of $\overset{\smile}{ {\cal A}_n}$ as $\left\{\breve{ t}_1, \cdots ,\breve{ t}_{|\overset{\smile}{ {\cal A}_n} |} \right\}$ and the corresponding estimator of the regression parameter for some segment between two consecutive change-point estimates $\breve{ t}_{j-1}$ and $\breve{ t}_j$, is denoted as $\overset{\smile}{\ef}_{\breve{ t}_j}$.
Thus, we have that 
\[
\overset{\smile}{\eb}_i=\overset{\smile}{\ef}_{\breve{ t}_j}, \qquad \textrm{for} \quad i=\breve{ t}_{j-1}, \cdots , \breve{ t}_j-1, \quad \textrm{and} \quad j=1, \cdots |\overset{\smile}{ {\cal A}_n}|.
\] 
From the proof of Theorem \ref{Theorem 3.3}, we also deduce that $\lim_{n  \rightarrow \infty} \PP\big[|\overset{\smile}{ {\cal A}_n}| <K^* \big]=0$.
Let us now denote the consecutive differences as  $\overset{\smile}{\eth}_i=\overset{\smile}{\eb}_i - \overset{\smile}{\eb}_{i-1}$, for $i=2, \cdots , n$, where $\overset{\smile}{\eth}_1=\overset{\smile}{\eb}_1$. For instance, if there are several consecutive change-point estimates $\overset{\vee}{t}_j$, which are, asymptotically,  all within a bounded distance from each other then we unify them into just one set $\overset{\vee}{T}$ and we only consider  the change-point estimate which is the smallest one among them (which will be the element of  $\overset{\smile}{ {\cal A}_n}$). Consequently, for the estimators of the vector parameters $\eb_i$, we take into account the   quantile fused estimator obtained between the last element of $\overset{\vee}{T}$ and the first of the consecutive set of analogous indexes.\\

\noindent\underline{\textit{Example}}\\
Suppose we have the following situation: $\overset{\vee}{ t}_{l} < \overset{\vee}{ t}_{l+1}< \overset{\vee}{ t}_{l+2}< \overset{\vee}{ t}_{l+3}< \overset{\vee}{ t}_{l+4}$, such that $ \overset{\vee}{ t}_{l+1}-  \overset{\vee}{ t}_{l} \rightarrow \infty$, $ \overset{\vee}{ t}_{l+2}-  \overset{\vee}{ t}_{l+1} $ and $\overset{\vee}{ t}_{l+3}-  \overset{\vee}{ t}_{l+2} $ are bounded, and  $ \overset{\vee}{ t}_{l+4}-  \overset{\vee}{ t}_{l+3} \rightarrow \infty$, all in probability, as $n$ tends to infinity. Thus, only $\overset{\vee}{ t}_{l}$, $\overset{\vee}{ t}_{l+1}$, $\overset{\vee}{ t}_{l+4}$ will be included in $\overset{\smile}{ {\cal A}_n}$, and the corresponding quantile fused estimators for $\eb_i$ are $\overset{\vee}{ \ef}_{l+1}$, for any $i=\overset{\vee}{ t}_{l}, \overset{\vee}{ t}_{l}+1,\cdots, \overset{\vee}{ t}_{l+1}-1$ and $\overset{\vee}{ \ef}_{l+4}$, for any $i=\overset{\vee}{ t}_{l+1}, \overset{\vee}{ t}_{l+1}+1, \cdots, \overset{\vee}{ t}_{l+2}, \cdots \overset{\vee}{ t}_{l+3} \cdots, \overset{\vee}{ t}_{l+4}-1$.

\begin{remark}
	\label{Rq_vit}
	Under Assumptions (A1)-(A8), if, moreover,  $(n \delta_n)^{-1}I^*_{min} {\underset{n \rightarrow \infty}{\longrightarrow}} \infty$,  $K^* \leq  |\overset{\smile}{ {\cal A}_n}| \leq K_{max}$, then,  by Theorems \ref{Theorem 3.1ii}, \ref{Theorem 3.2} and \ref{Theorem 3.3},  it holds that
	\[
	\sup_{\scriptsize\begin{array}{c}
		j \in \{1, \cdots,|\overset{\smile}{ {\cal A}_n} |+1   \} \\
		|\breve{ t}_j-t^*_k| \leq n \delta_n
		\end{array}} \| \overset{\smile}{ \ef}_{j} -\ef^*_k \|=O_{\PP} (b_n ),
	\]
	for any $k=1, \cdots , K^*+1$.
\end{remark}

The performance of the quantile fussed penalty can be further improved. Indeed, if we consider only the change-point estimates belonging to the set  $\overset{\smile}{ {\cal A}_n}$, then we can define weights for the adaptive penalty and to use the idea of the adaptive LASSO instead as the adaptive LASSO approach is well known for having some better selection performance in general. The adaptive fused penalty generalization is considered in more details in the next section.

 \section{Adaptive fused  quantile method} 
\label{section4}
In this section we provide an alternative method for the automatic change-point detection in the linear model and we introduce an adaptive extension for the fused penalty approach discussed in the previous section. For the purpose of this section, we suppose  that the assumptions given in Remark \ref{Rq_vit} are all satisfied. Hence, we have the following:
\begin{equation}
\label{cond_R31}
\left\{
\begin{array}{l}
\textrm{Assumptions (A1)-(A8) hold;} \\
(n \delta_n)^{-1}I^*_{min} \rightarrow \infty, \quad \textrm{as $n \to \infty$}; \\
K^* \leq  |\overset{\smile}{ {\cal A}_n}| \leq K_{max}, \quad  \textrm{ holds with probability converging to one.}
\end{array}
\right.
\end{equation}

It is clear from the assumptions above that the adaptive fused quantile  method can be only considered if there are at least as many change-points being detected by the quantile fused method as there are true change-points $K^* \in \mathbb{N}$, and, in addition, the number of estimated change-points in the model is bounded from above. 

As an extension to \eqref{eq3} let us define the adaptive version of the quantile process 
\begin{equation}
\label{aeq4}
S(\eb^n) \equiv \sum^n_{i=1} \rho_\tau(Y_i - \ex_i^\top  \eb_i) +n \lambda_n \sum^n_{i=2} \omega_i \| \eb_i-\eb_{i-1}\|,
\end{equation}
with the weights $\omega_i$ depending on the differences between two consecutive quantile fused estimators and some deterministic sequence $(d_n)$, for $i=2, \cdots , n$, such that 
\[
\omega_i \equiv \big(\max(\| \overset{\smile}{ \eth}_i  \|_\infty, d_n) \big)^{-\gamma} =\bigg( \max \big(\|\overset{\smile}{ \eb}_i - \overset{\smile}{ \eb}_{i-1} \|_\infty, d_n \big)\bigg)^{-\gamma},
\] 
where $\gamma >0$ is some positive constant. Moreover, it is assumed that the sequence $(d_n)_n$ satisfies 
\begin{equation}
\label{CD1}
d_n {\underset{n \rightarrow \infty}{\longrightarrow}} 0, \quad \quad \textrm{and} \quad \quad \frac{\lambda_n}{\delta_n\big(\max(d_n,b_n)\big)^{\gamma}}  {\underset{n \rightarrow \infty}{\longrightarrow}} 0,
\end{equation}
as $n \rightarrow \infty$. In fact, the relation  in \eqref{CD1} can be used to replace the condition in  Assumption (A7), where we need that $\lambda_n/\delta_n \rightarrow 0$. As an example of such sequences we can consider, for instance, $\lambda_n=n^{-1} (\log n)^{5/2}$, $\delta_n=n^{-1/2} (\log n)^{3}$, and $I^*_{min} = n/d_1$, with some constant $d_1 \in (0,1)$. Then, we obtain that $b_n=n^{-1/2}$ and for $d_n=n^{-1/2}$ in \eqref{CD1} we need that $\gamma \leq 1$.

The adaptive fused  quantile estimators for $\eb^n \in \mathbb{R}^{np}$ and $\eth^n \in \mathbb{R}^{np}$ respectively, are defined as  
\[
\widehat{\eb^n} \equiv  \argmin_{\eb^n \in \R^{np}} {S}(\eb^n), \qquad \textrm{and} \qquad \widehat{\eth^n} \equiv \big( \widehat{\eb}_1^\top,(\widehat{\eb}_2-\widehat{\eb}_1)^\top, \cdots , (\widehat{\eb}_n-\widehat{\eb}_{n-1})^\top  \big)^\top,
\]
where $\widehat{\eb^n} =(\widehat{\eb}_1^\top, \cdots , \widehat{\eb}_n^\top)^\top$ and $\widehat{\eth^n}=(\widehat{\eth}_1^\top, \cdots , \widehat{\eth}_n^\top)^\top$. The corresponding estimates of the change-point  locations are the observations where $\widehat{\eb}_i \neq \widehat{\eb}_{i-1} $. Let $\widehat{{\cal A}}_n$ denotes the set of indexes, such that 
\begin{equation}
\label{ahAn}
\widehat{{\cal A}}_n \equiv \{ i \in \{2, \cdots , n \}; \; \widehat{\eb}_i \neq \widehat{\eb}_{i-1}  \} = \{ \hat t_1, \cdots ,\hat t_{|\widehat{\cal A}_n|} \},
\end{equation}
where $|\widehat{\cal A}_n|$ is the cardinality of $\widehat{{\cal A}}_n$. Thus, for any $k \in \{0, \cdots , |\widehat{{\cal A}}_n|\}$, and the $(k+1)$-th model phase (e.i., observations between $\hat t_k$ and $\big(\hat t_{k+1}-1 \big)$, respectively $\hat t_{|\widehat{{\cal A}}_n|+1 }$ for the last phase) the corresponding parameter estimator within the given phase is equal to $\widehat{\ef}_{k+1}$, with $\hat t_0=1$ and $\hat t_{|\widehat{{\cal A}}_n|+1}=n$. In the following theorems we state some important properties of the estimation approach based on the adaptive fused penalty. The proofs are all postponed to the appendix.

\begin{theorem}
	\label{aTheorem 3.1}
	Under the assumptions in (\ref{cond_R31}), together with the condition in (\ref{CD1}), if $|\widehat{\cal A}_n|=K^*$, it holds that
	\[
	\lim_{n  \rightarrow \infty } \PP \left[\max_{1 \leq k \leq K^*} | \hat t_k - t^*_k | \leq n \delta_n  \right]  =  1.
	\]
\end{theorem}  

The theorem above gives the consistency property of the change-point location estimators given by the adaptive fused approach. In the following theorem  we state the convergence rate of the adaptive fused quantile  estimators of $\ef^*_k$ when the estimated number of change-points coincides with $K^*$. Compared to the rate $b_n$ of $\overset{\smile}{\ef}_k$ given by Theorem \ref{Theorem 3.1ii}, the convergence rate  of $ \widehat{\ef}_k$ depends, in addition, on $I^*_{min}$ and $\lambda_n$, and also on $b_n$ and the given sequence $d_n$.
The convergence rate for the adaptive fused quantile estimator of the unknown regression parameters is one of the main contribution of this paper.

\begin{theorem}
	\label{aTheorem 3.1ii}
	Under the same conditions as in Theorem  \ref{aTheorem 3.1}, we have,  for any $k = 1, \cdots , K^*+1$, that
	\[
	\big\| \widehat{\ef}_k-\ef^*_k \big\|=O_{\PP} \bigg(\frac{n \lambda_n \big(\max(d_n,b_n) \big)^{-\gamma}}{I^*_{min}} +\frac{1}{\sqrt{I^*_{min}}}\bigg).
	\]
\end{theorem}

\begin{remark}	
	Comparing the results of Theorems \ref{Theorem 3.1ii} and \ref{aTheorem 3.1ii}, for $n \lambda_n (I^*_{min})^{-1/2} \geq C>0$, and $n \in \mathbb{N}$  large enough, the convergence rate of the adaptive  fused quantile estimator $ \widehat{\ef}_k$ can be slower than the convergence rate of the quantile fused estimator $ \overset{\vee}{ \ef}_k$, for $k \in \{1, \cdots , K^*+1\}$. This is due to the presence of the deterministic sequence $d_n$ in the weights $\omega_i$ which makes the minimization of the quantile loss function, under constraints, possibly slower.
\end{remark}

The results above are given for a specific situations where the number of detected change-points coincides with the truth. In the following theorem we show that even if the number of change-points is overestimated then for each true change-point there is at least one  estimated change-point at a distance less than $n \delta_n$.

\begin{theorem}
	\label{aTheorem 3.2}
	Under the assumptions in (\ref{cond_R31}), together with the condition in (\ref{CD1}), if, in addition, $K^*=  |{\cal A}^*| \leq |\widehat{{\cal A}}_n| \leq K_{max} < \infty $, then 
	$$\lim_{n \rightarrow \infty}  \PP \big[{\cal E}(\widehat{{\cal A}}_n|| {\cal A}^*) \leq n \delta_n \big] = 1.$$
\end{theorem}

The upper bound, $K_{max}$, for the number of change-points in Theorem \ref{aTheorem 3.2}  may be arbitrary but bounded and it can differ from the one considered in Theorem \ref{Theorem 3.2}. In contrast to the theorem above, which deals with the situation where the number of change-points is overestimated, we can define an additional condition for the sequences $(\lambda_n)$, $(d_n)$, $(b_n)$, and $I^*_{max}$, and $\gamma > 0$, such that the adaptive fused quantile  method does not underestimate the true number of change-points. Specifically, if we require that 
\begin{equation}
\label{CD2}
\frac{n \lambda_n \big(\max(d_n,b_n) \big)^{-\gamma}}{\sqrt{I^*_{max}}} {\underset{n \rightarrow \infty}{\longrightarrow}} \infty,
\end{equation}
then it can be proved that the adaptive fused quantile estimation approach underestimates the true number of change-points with probability tending to zero, as $n \in \mathbb{N} \to \infty$  (see the next theorem).  A straightforward example of sequences $(\lambda_n)$, $(d_n)$, $(b_n)$, and $I^*_{max}$, and the value of $\gamma > 0$, which satisfy the condition in (\ref{CD2}), are, for instance,  $\gamma=1$, $b_n=d_n=n^{-1/2}$, $I^*_{max} = n/d_2$, with the constant $d_2 \in (0,1)$, $\lambda_n = n^{-1} (\log n)^{5/2}$. Taking also into account the  relation in \eqref{CD1}, we deduce by relation \eqref{CD2} that the maximum distance between two successive change-points must be much smaller than the square of the minimum distance between two successive change-points: $I^*_{max}/(I^*_{min})^2 \rightarrow 0$ as $n \rightarrow \infty$. Comparing this with the example sequences, the possibility that $I^*_{min}=n/d_1$, with, $0< d_1 \leq d_2 <1$ satisfies this condition.

\begin{theorem}
	\label{aTheorem 3.3}
	Under the assumptions in \eqref{cond_R31}, together with the conditions in (\ref{CD1}) and (\ref{CD2}),  it holds that 
	$$\lim_{n \rightarrow \infty}  \PP \big[|\widehat{{\cal A}}_n| < |{\cal A}^*| \big] = 0.$$
\end{theorem} 

\begin{corollary}
	\label{Raq_vit}
	Under the assumptions in \eqref{cond_R31}, \eqref{CD1}, and \eqref{CD2},  if the derivative of the density function $f$ is bounded in some neighborhood of zero, and for any $j \in \{1, \cdots , |\widehat{{\cal A}}_n| +1\}$,  such that $\hat t_j -\hat t_{j-1} \overset{\PP} {\underset{n \rightarrow \infty}{\longrightarrow}} \infty$, it holds that 
	$$\PP[\widehat{\cal A}_n  ={\cal A}^*] \quad {\underset{n \rightarrow \infty}{\longrightarrow}} \quad 1.$$
\end{corollary} 

The corollary above provides a very interesting result for the situation where  the distance between any two consecutive estimated change-points converges to infinity. In such case the number of change-points being estimated by the adaptive fused approach corresponds, with a probability converging to one, with the true number of change-points, $K^* \in \mathbb{N}$. 
For the density function $f$ we assume that there is some neighborhood of zero, such that $\{x \in \R; \; |x| \leq \eta  \}$, with $\eta \searrow 0$.  The corollary above also shows that a judicious choice of the sequences $(\lambda_n)$, $(\delta_n)$ leads to a consistent estimation of the number of change-points detected by the adaptive fused quantile method, which holds with probability converging to one, as $n \in \mathbb{N}$ tends to infinity.  Taking also into account Theorem \ref{aTheorem 3.1}, the (only) estimator of each change-point is at a distance less than $n \delta_n$ from the true change-point, and, moreover, the regression parameter estimates between two consecutive change-point estimates converge, by Theorem  \ref{aTheorem 3.1ii}, to the true values with the convergence rate of $ n \lambda_n \big(\max(d_n,b_n) \big)^{-\gamma} (I^*_{min})^{-1} + (I^*_{min})^{-1/2}$. If the distance between any two consecutive  change-points, estimated by the adaptive fused method, converges to infinity, then the probability    of an overestimation of the number of change-points converges to zero (see the proof of Corollary \ref{Raq_vit} and, especially, the relation in (\ref{SP_eq31})).

\section{Simulation study}
\label{section5}
In this section we empirically compare the quantile fused method and the adaptive version of the quantile fused method, which are both proposed in this paper. In addition, we also consider a competitive estimation algorithm proposed in \cite{Qian.Su.16}, which we refer to as a \textit{standard LASSO approach}. The standard LASSO approach is used to estimate the conditional expectation in the model while the proposed fused methods are both used to estimate the conditional median ($\tau = 0.5$). All three methods are  compared for a wide range of different scenarios (for instance, different error distributions, signal-to-noise ratio, change-point magnitudes, sample size, or the model selection strategy). We also considered various quantile levels $\tau \in (0,1)$ and dimensions $p \in \mathbb{N}$, however, only applied for the quantile fused approaches.  

\begin{figure}
	\centering
	\subfigure[Standard LASSSO $|$ Distribution $N(0,1)$] {\label{fig3:a}\includegraphics[width=0.48\textwidth]{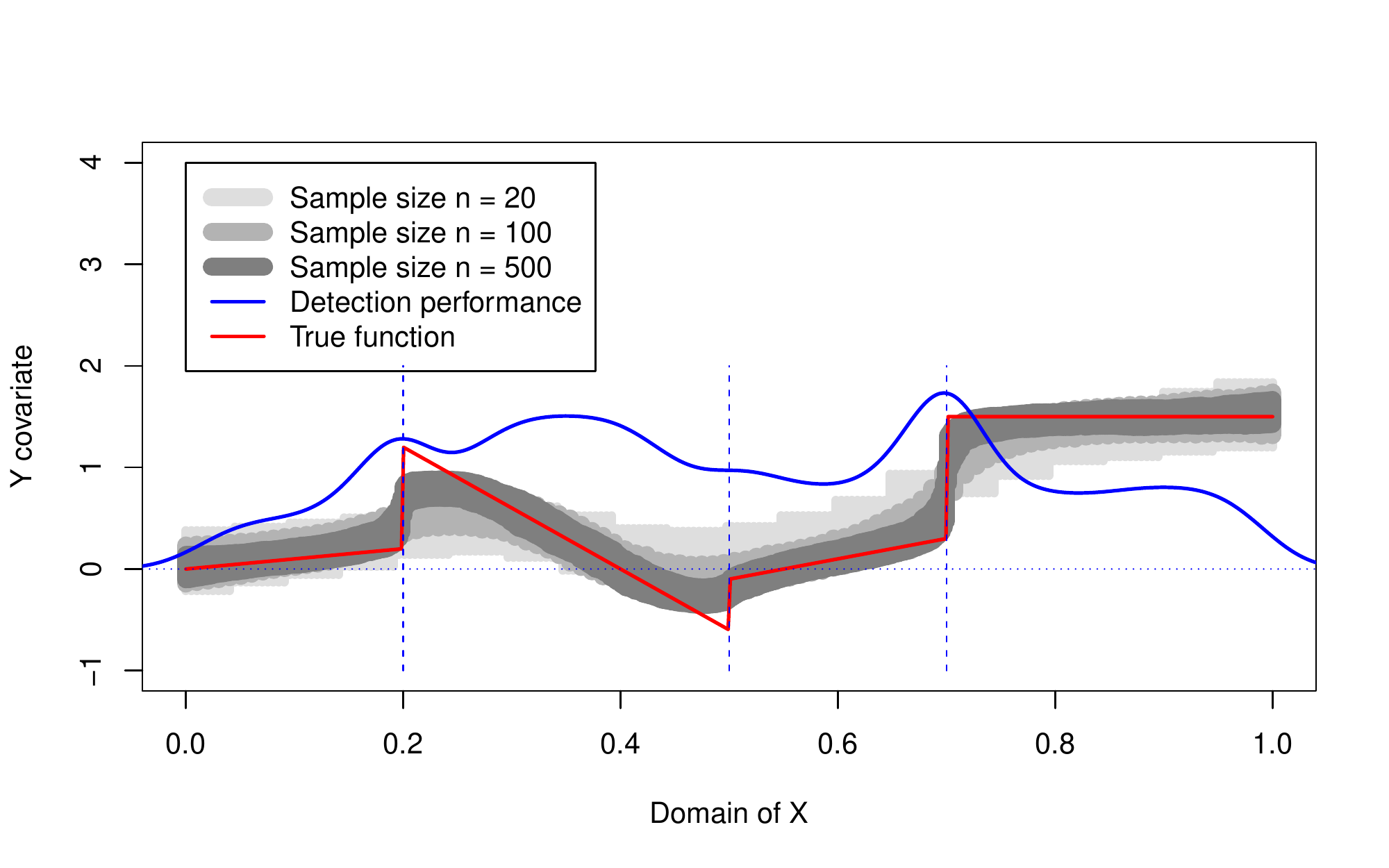}}
	\subfigure[Quantile LASSSO $|$ Distribution $N(0,1)$] {\label{fig3:b}\includegraphics[width=0.48\textwidth]{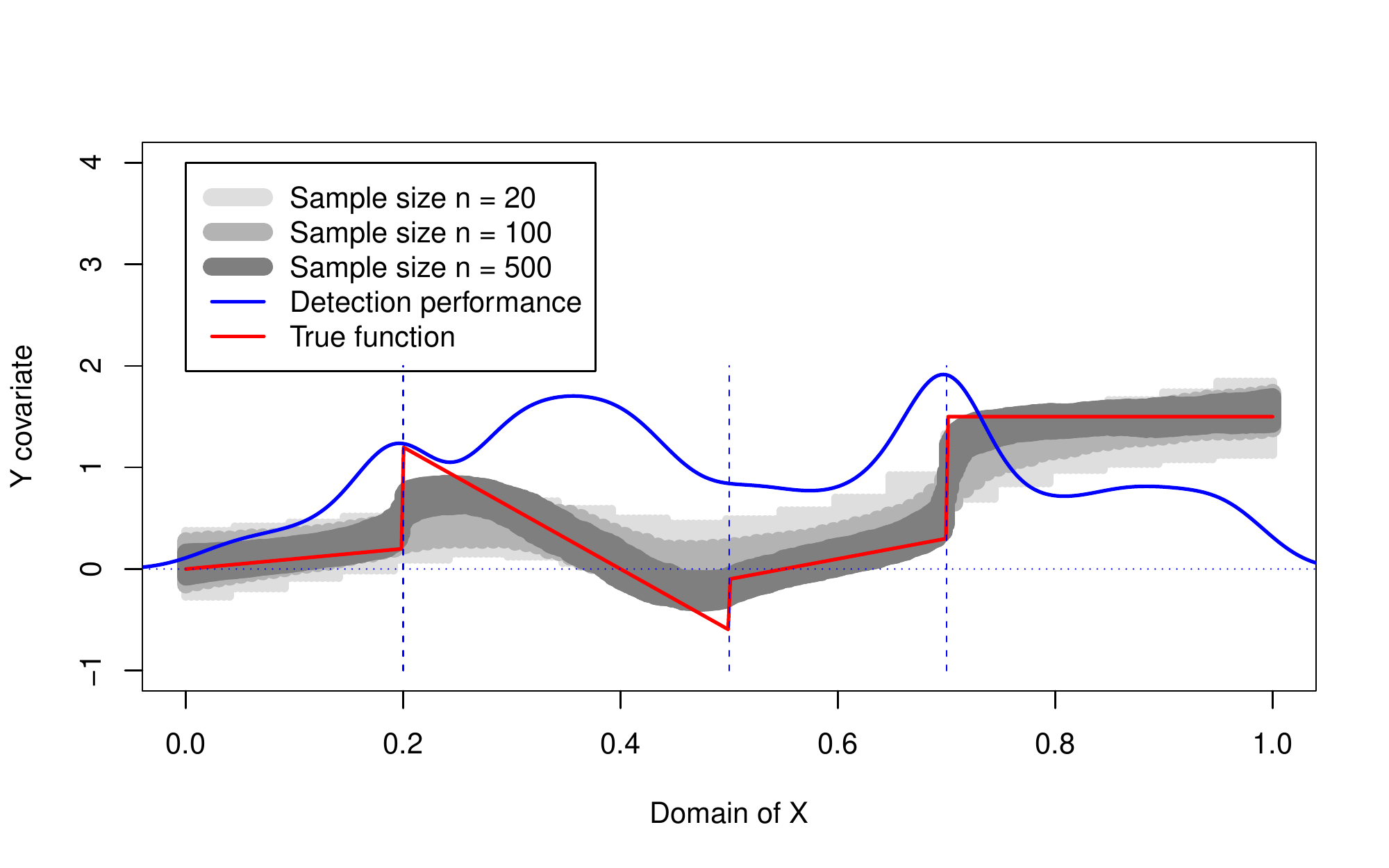}}
	
	\subfigure[Standard LASSO $|$ Distribution $t_3$] {\label{fig3:a}\includegraphics[width=0.48\textwidth]{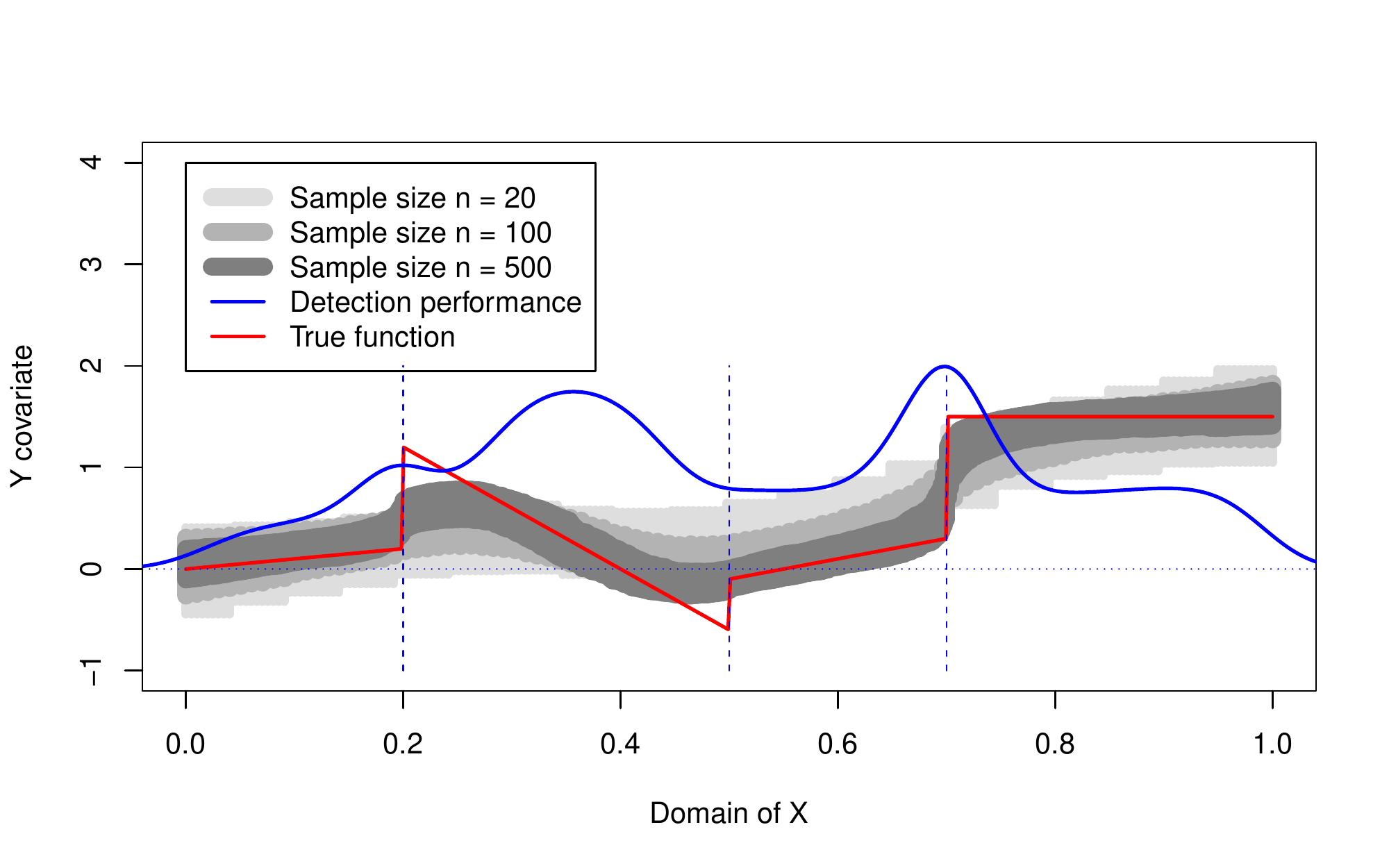}}
	\subfigure[Quantile LASSO $|$ Distribution $t_3$] {\label{fig3:b}\includegraphics[width=0.48\textwidth]{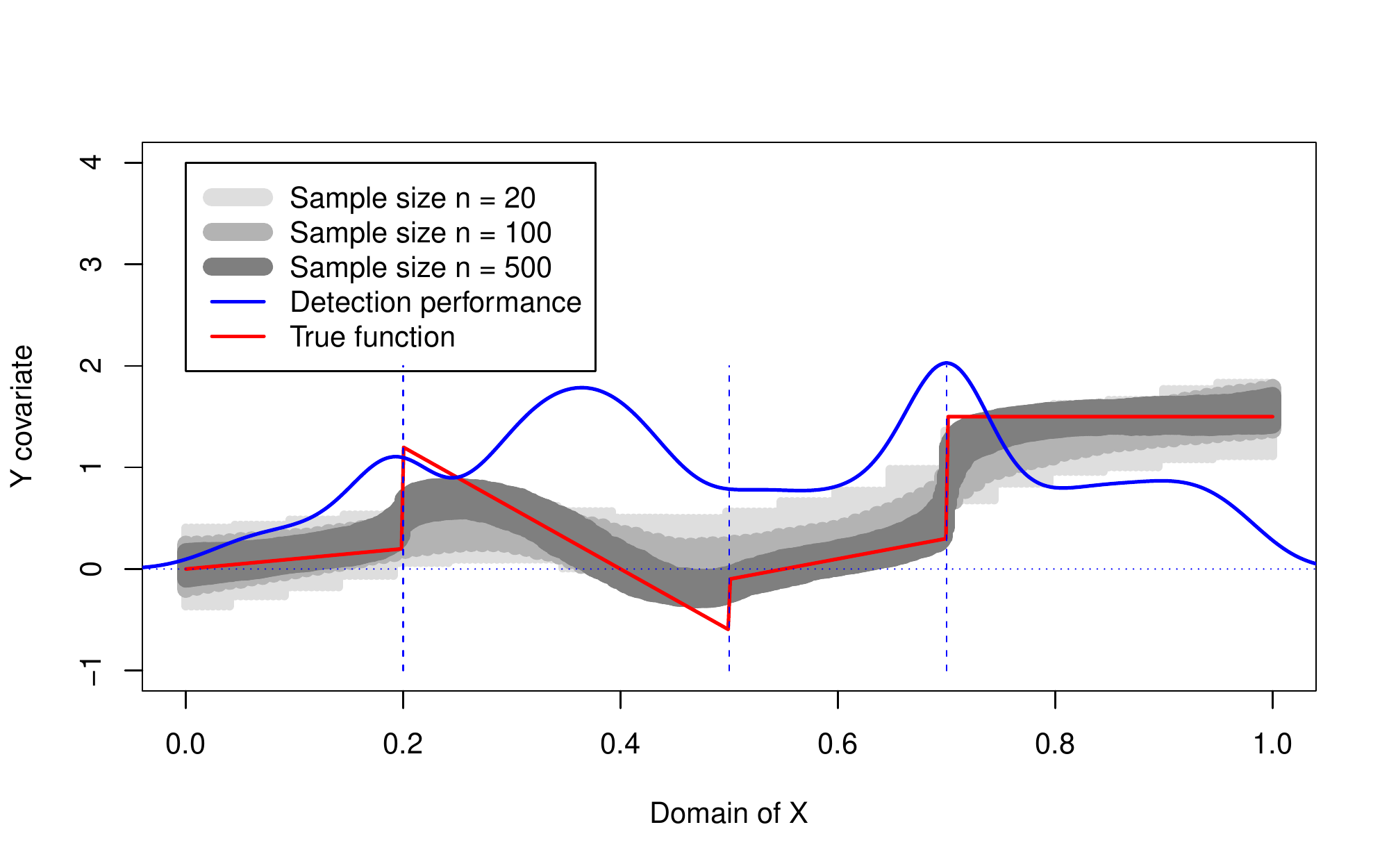}}
	
	\subfigure[Standard LASSO $|$ Distribution $C(0,1)$] {\label{fig3:a}\includegraphics[width=0.48\textwidth]{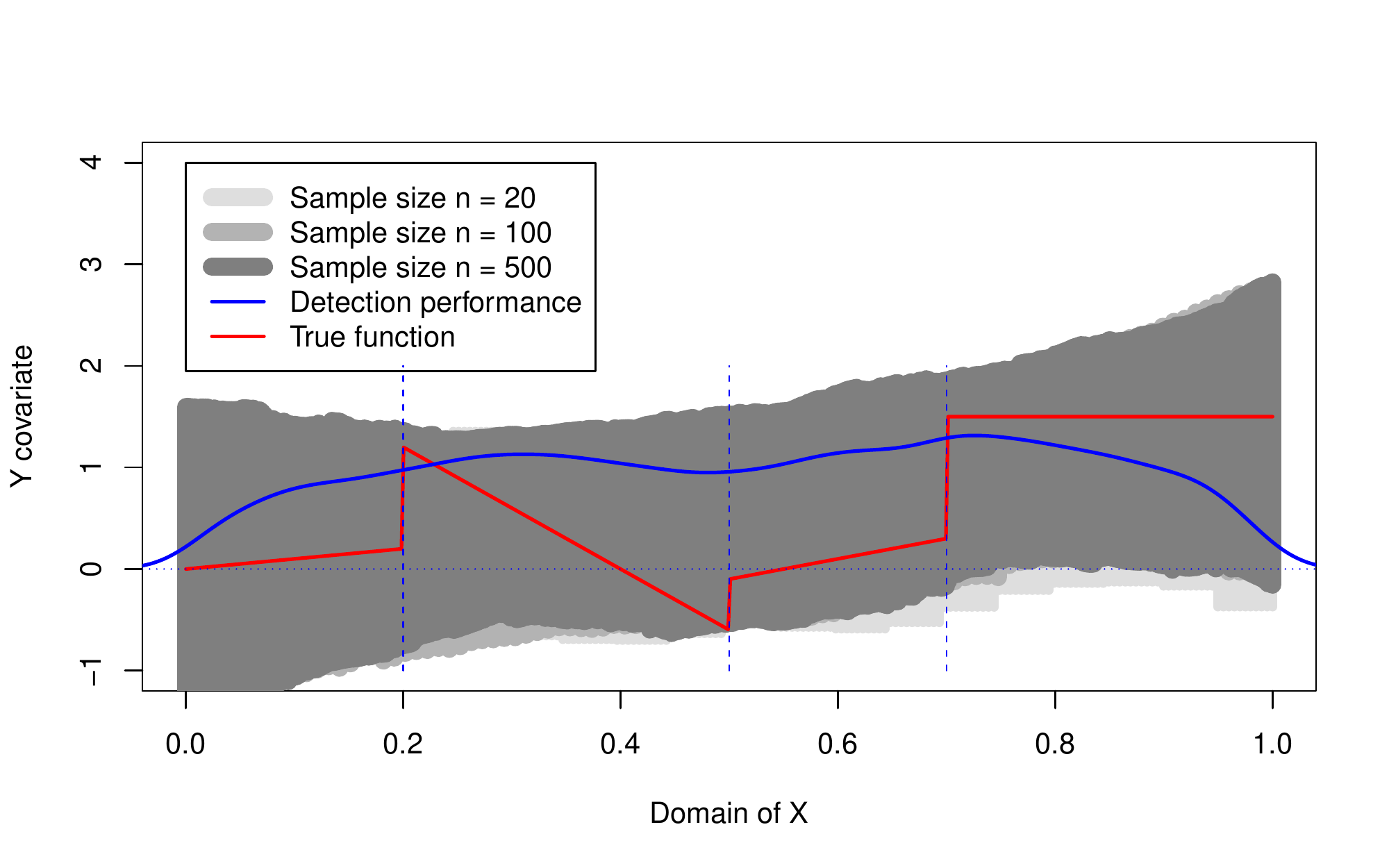}}
	\subfigure[Quantile LASSO $|$ Distribution $C(0,1)$] {\label{fig3:b}\includegraphics[width=0.48\textwidth]{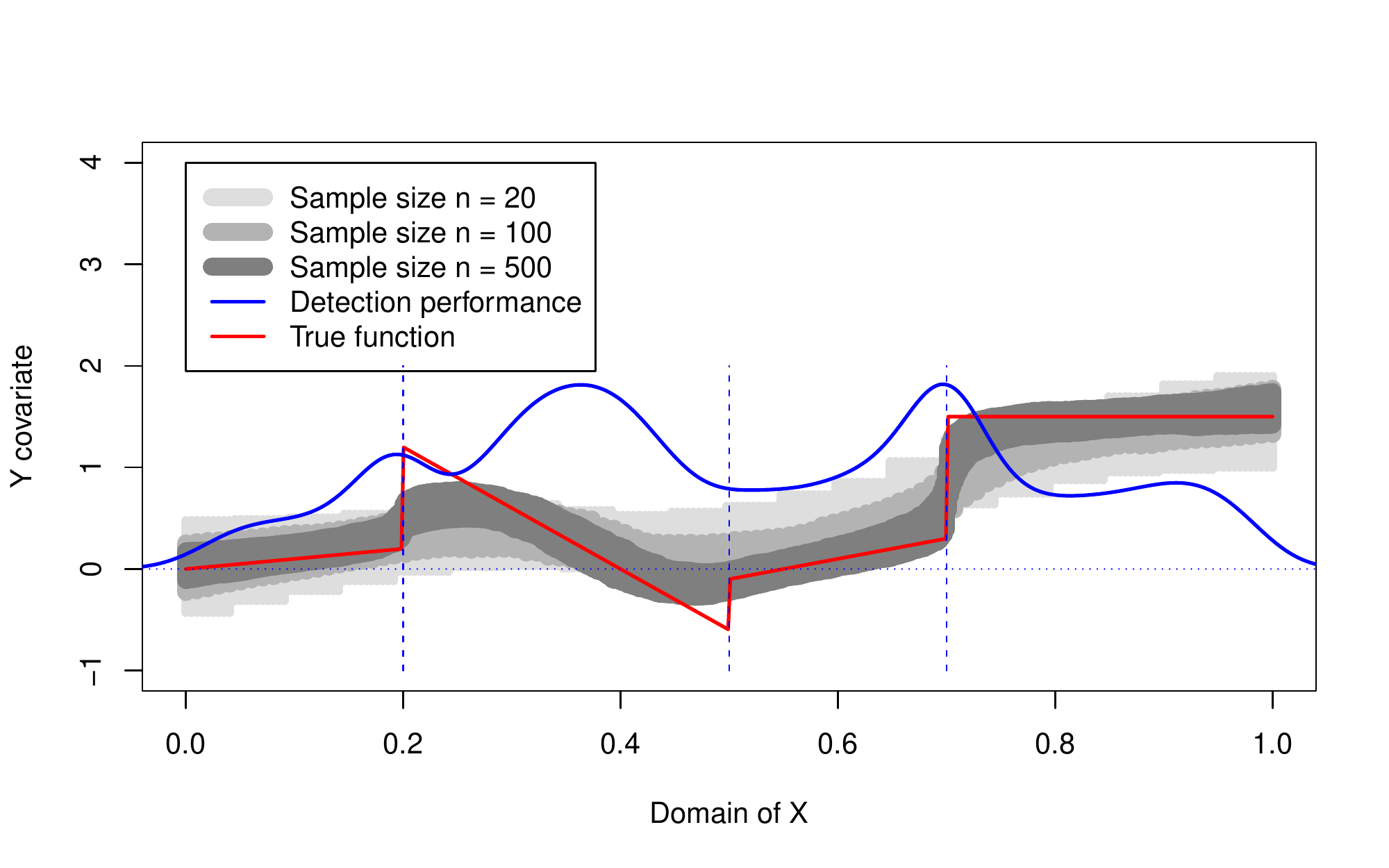}}
	
	\caption{\footnotesize Graphical comparison of the empirical performance of the standard LASSO approach and the proposed fused quantile method (denoted as Quantile LASSO) for three error distributions and three sample sizes. The point-wise interquartile bands are provided for each scenario and the overall change-point detection performance is visualized in terms of a rescaled density of estimated change-point locations for the sample size $n = 500$ calculated out of 500 Monte Carlo simulations (blue solid line). The red solid line shows the true underlying model.}
	\label{fig1}
\end{figure}

For illustration, a simple linear model for $p = 2$ and three change-points (and thus, four model stages) is considered to compare the empirical performance of three different estimation techniques. However, to be able to directly compare  models for different number of observations, $n \in \mathbb{N}$, in just one single graph, we always rescale the model defined in \eqref{eq1} and \eqref{eq2} such that each index $i \in \{1, \dots, n\}$ will be expressed as $\tilde{i} = i/n$. Hence, without any loss of generality we can use a common domain for the underlying model which will be the interval $[0,1]$. The underlying dependence takes the form

\begin{equation}
\label{sim_function}
\boldsymbol{x}_i^\top \boldsymbol{\beta}_i = \left\{ \begin{array}{ll}
\boldsymbol{x}_i^\top (0, 1)^\top & \textrm{for $\tilde{i} = i/n \in (0,0.2)$,}\\
\boldsymbol{x}_i^\top (2.4, -6)^\top & \textrm{for $\tilde{i} = i/n \in [0.2, 0.5)$,}\\
\boldsymbol{x}_i^\top (-1.1, 2)^\top & \textrm{for $\tilde{i} = i/n \in [0.5, 0.7)$,}\\
\boldsymbol{x}_i^\top (0.5, 0)^\top & \textrm{for $\tilde{i} = i/n \in [0.7, 1]$,}
\end{array}  \right.
\end{equation}

where $\boldsymbol{x}_{i} = (1, \tilde{i})^\top = (1, i/n)^\top$. The underlying function in \eqref{sim_function} is defined such that various situations are implicitly included in the model: the first change-point location $\xi_1 = t_1^*/n = 0.2$ 
introduces a relatively small jump but a huge change in the slope (respectively, $|\boldsymbol{\phi}_2 - \boldsymbol{\phi}_1| = (2.4, 7)^\top$); the second change-point location, $\xi_2 = t_2^*/n = 0.5$ 
introduces large magnitudes for the change in both, the function itself and its derivative, and, moreover, it compensates the effect of the first change in some sense (equivalently, we have $|\boldsymbol{\phi}_3 - \boldsymbol{\phi}_2| = (3.5, 8)^\top$); finally, relatively small magnitudes for the jump and the slope change are observed at the third change-point location (i.e., $|\boldsymbol{\phi}_4 - \boldsymbol{\phi}_3| = (0.6, 2)^\top$). In addition, the model phases have various lengths and different number of observations are therefore expected to occur in each phase (see Figure \ref{fig1} for more details). 

Three error distributions are considered (standard normal, $t$-distribution with three degrees of freedom, and the Cauchy distribution), three different sample sizes are used ($n \in \{20, 100, 500\}$), and the final model is obtained by one of the three selection procedures: the first procedure uses the prior knowledge of three change-points in the model and the corresponding regularization parameter is denoted as $\lambda_{(3)}$; the second model is defined by the regularization parameter $\lambda_{AS} = n^{-1} (\log n)^{5/2}$, which satisfies the theoretical assumptions needed for the proofs to hold; finally, the last model selection procedure is defined by the regularization parameter $\lambda_{MS}$ which minimizes the theoretical mean squared error quantity. The models are always compared with respect to various qualities: the estimation performance is assessed by using the empirical bias  $n^{-1} \sum_{i = 1}^n (\boldsymbol{x}_i^\top \boldsymbol{\beta}_i^* - \boldsymbol{x}_i^\top \widehat{\boldsymbol{\beta}_i})$ and the empirical mean squared error (MSE) term $n^{-1} \sum_{i = 1}^n (\boldsymbol{x}_i^\top \boldsymbol{\beta}_i^* - \boldsymbol{x}_i^\top \widehat{\boldsymbol{\beta}_i})^2$, where $\boldsymbol{\beta}_i^*$ is the true value of the parameter and $\widehat{\boldsymbol{\beta}_i^*}$ is the corresponding estimate; In addition, the change-point detection performance is assessed via the change-point detection error, defined as $1/3 \sum_{k = 1}^{3} |\widehat{t}_k - t_{k}^*|$, however, provided only in situations where at least three change-points are detected in the model. For a~more detailed comparison we also report some overall insight into the number of change-points being detected in each scenario. The results are summarized in Tables \ref{tab2} and \ref{tab3}.

\begin{table}[!ht]\footnotesize
	\begin{center}\scalebox{0.8}{	
		\begin{tabular}{ccc|r<{\hspace{-\tabcolsep}}>{\hspace{-\tabcolsep}\,}lr<{\hspace{-\tabcolsep}}>{\hspace{-\tabcolsep}\,}lr<{\hspace{-\tabcolsep}}>{\hspace{-\tabcolsep}\,}l
    r<{\hspace{-\tabcolsep}}>{\hspace{-\tabcolsep}\,}lr<{\hspace{-\tabcolsep}}>{\hspace{-\tabcolsep}\,}lr<{\hspace{-\tabcolsep}}>{\hspace{-\tabcolsep}\,}l}
\hiderowcolors
\multirow{2}{*}{$\boldsymbol{\mathcal{D}}$} & & \multirow{2}{*}{$\boldsymbol{n}$} & 
     \multicolumn{4}{c}{\textbf{Model with $\boldsymbol{\lambda_{(3)}}$}} & \multicolumn{4}{c}{\textbf{Model with $\boldsymbol{\lambda_{AS}}$}} & \multicolumn{4}{c}{\textbf{Model $\boldsymbol{\lambda_{MS}}$}} \\
 ~ & ~ & ~ & 
     \multicolumn{2}{c}{\scalebox{1}{Est. Bias}} & \multicolumn{2}{c}{MSE} & \multicolumn{2}{c}{\scalebox{1}{Est. Bias}} & \multicolumn{2}{c}{MSE} & \multicolumn{2}{c}{\scalebox{1}{Est. Bias}} & \multicolumn{2}{c}{MSE}  \\\hline\hline
\multicolumn{3}{c}{~} & \multicolumn{12}{c}{~}\\[-0.4cm]
$\boldsymbol{N}$ & & \textbf{20} & 0.00 & \textit{(0.23)} & 0.28 & \textit{(0.12)} & 0.00 & \textit{(0.23)} & 0.29 & \textit{(0.13)} & 0.00 & \textit{(0.23)} & 0.25 & \textit{(0.11)} \\
 & & \textbf{100} & 0.00 & \textit{(0.10)} & 0.17 & \textit{(0.05)} & 0.00 & \textit{(0.10)} & 0.11 & \textit{(0.04)} & 0.00 & \textit{(0.10)} & 0.10 & \textit{(0.03)} \\
 &  \multirow{-3}{*}{\rotatebox{90}{SLasso}} & \textbf{500} & 0.00 & \textit{(0.04)} & 0.15 & \textit{(0.04)} & 0.00 & \textit{(0.04)} & 0.05 & \textit{(0.01)} & 0.00 & \textit{(0.04)} & 0.04 & \textit{(0.01)} \\
\multicolumn{3}{c}{~} & \multicolumn{12}{c}{~}\\[-0.4cm]
\rowcolor{Gray} &  & \textbf{20} & 0.00 & \textit{(0.29)} & 0.37 & \textit{(0.18)} & 0.01 & \textit{(0.28)} & 0.37 & \textit{(0.19)} & 0.00 & \textit{(0.24)} & 0.29 & \textit{(0.14)} \\
\rowcolor{Gray} & & \textbf{100} & 0.00 & \textit{(0.13)} & 0.18 & \textit{(0.05)} & 0.00 & \textit{(0.13)} & 0.14 & \textit{(0.05)} & 0.01 & \textit{(0.12)} & 0.13 & \textit{(0.05)} \\
\rowcolor{Gray} & \multirow{-3}{*}{\rotatebox{90}{QLasso}} & \textbf{500} & 0.00 & \textit{(0.06)} & 0.16 & \textit{(0.04)} & 0.00 & \textit{(0.06)} & 0.06 & \textit{(0.01)} & 0.00 & \textit{(0.05)} & 0.05 & \textit{(0.01)} \\
\rowcolor{white} & \multicolumn{2}{c}{~} & \multicolumn{12}{c}{~}\\[-0.4cm]
\rowcolor{Gray2} &  & \textbf{20} & 0.01 & \textit{(0.28)} & 0.40 & \textit{(0.18)} & 0.00 & \textit{(0.28)} & 0.40 & \textit{(0.19)} & 0.00 & \textit{(0.26)} & 0.33 & \textit{(0.16)} \\
\rowcolor{Gray2} & & \textbf{100} & 0.00 & \textit{(0.13)} & 0.21 & \textit{(0.07)} & 0.00 & \textit{(0.13)} & 0.20 & \textit{(0.07)} & 0.00 & \textit{(0.13)} & 0.21 & \textit{(0.07)} \\
\rowcolor{Gray2} & \multirow{-3}{*}{\rotatebox{90}{ALasso}} & \textbf{500} & 0.00 & \textit{(0.06)} & 0.16 & \textit{(0.06)} & 0.00 & \textit{(0.06)} & 0.17 & \textit{(0.06)} & 0.00 & \textit{(0.06)} & 0.16 & \textit{(0.06)} \\
\rowcolor{white} & \multicolumn{2}{c}{~} & \multicolumn{12}{c}{~}\\[-0.4cm]
\hline\multicolumn{3}{c}{~} & \multicolumn{12}{c}{~}\\[-0.4cm]
$\boldsymbol{t_3}$ &  & \textbf{20} & 0.00 & \textit{(0.39)} & 0.67 & \textit{(0.90)} & 0.00 & \textit{(0.39)} & 1.03 & \textit{(1.55)} & 0.00 & \textit{(0.39)} & 0.46 & \textit{(0.38)} \\
 & & \textbf{100} & 0.00 & \textit{(0.17)} & 0.23 & \textit{(0.21)} & 0.00 & \textit{(0.17)} & 0.60 & \textit{(0.90)} & 0.00 & \textit{(0.17)} & 0.18 & \textit{(0.08)} \\
 & \multirow{-3}{*}{\rotatebox{90}{SLasso}} & \textbf{500} & 0.00 & \textit{(0.08)} & 0.17 & \textit{(0.05)} & 0.00 & \textit{(0.08)} & 0.52 & \textit{(0.83)} & 0.00 & \textit{(0.08)} & 0.07 & \textit{(0.03)} \\
\multicolumn{3}{c}{~} & \multicolumn{12}{c}{~}\\[-0.4cm]
\rowcolor{Gray} &  & \textbf{20} & 0.00 & \textit{(0.33)} & 0.49 & \textit{(0.42)} & 0.01 & \textit{(0.33)} & 0.48 & \textit{(0.35)} & 0.00 & \textit{(0.28)} & 0.35 & \textit{(0.20)} \\
\rowcolor{Gray} & & \textbf{100} & 0.01 & \textit{(0.14)} & 0.19 & \textit{(0.06)} & 0.01 & \textit{(0.14)} & 0.16 & \textit{(0.06)} & 0.01 & \textit{(0.12)} & 0.15 & \textit{(0.05)} \\
\rowcolor{Gray} & \multirow{-3}{*}{\rotatebox{90}{QLasso}} & \textbf{500} & 0.00 & \textit{(0.07)} & 0.16 & \textit{(0.04)} & 0.00 & \textit{(0.06)} & 0.06 & \textit{(0.02)} & 0.00 & \textit{(0.06)} & 0.05 & \textit{(0.02)} \\
\rowcolor{white} & \multicolumn{2}{c}{~} & \multicolumn{12}{c}{~}\\[-0.4cm]
\rowcolor{Gray2} &  & \textbf{20} & 0.01 & \textit{(0.33)} & 0.56 & \textit{(0.56)} & 0.01 & \textit{(0.33)} & 0.51 & \textit{(0.42)} & 0.00 & \textit{(0.29)} & 0.38 & \textit{(0.23)} \\
\rowcolor{Gray2} & & \textbf{100} & 0.00 & \textit{(0.14)} & 0.23 & \textit{(0.09)} & 0.00 & \textit{(0.14)} & 0.22 & \textit{(0.08)} & 0.00 & \textit{(0.13)} & 0.22 & \textit{(0.08)} \\
\rowcolor{Gray2} & \multirow{-3}{*}{\rotatebox{90}{ALasso}} & \textbf{500} & 0.00 & \textit{(0.07)} & 0.16 & \textit{(0.06)} & 0.00 & \textit{(0.07)} & 0.17 & \textit{(0.06)} & 0.00 & \textit{(0.07)} & 0.16 & \textit{(0.06)} \\
\rowcolor{white} & \multicolumn{2}{c}{~} & \multicolumn{12}{c}{~}\\[-0.4cm]
\hline\multicolumn{3}{c}{~} & \multicolumn{12}{c}{~}\\[-0.4cm]
$\boldsymbol{C}$ & & \textbf{20} & 1.72 & \textit{(24.45)} & 11330 & \textit{(165371)} & 1.72 & \textit{(24.45)} & 11749 & \textit{(168201)} & 1.72 & \textit{(24.45)} & 10947 & \textit{(162604)} \\
 & & \textbf{100} & -1.53 & \textit{(26.37)} & 32128 & \textit{(534460)} & -1.53 & \textit{(26.37)} & 72686 & \textit{(1025214)} & -1.53 & \textit{(26.37)} & 30918 & \textit{(528628)} \\
 & \multirow{-3}{*}{\rotatebox{90}{SLasso}} & \textbf{500} & -2.45 & \textit{(39.48)} & 642888 & \textit{(14071066)} & -2.45 & \textit{(39.48)} & 776598 & \textit{(14377993)} & -2.45 & \textit{(39.48)} & 625465 & \textit{(13726447)} \\
\multicolumn{3}{c}{~} & \multicolumn{12}{c}{~}\\[-0.4cm]
\rowcolor{Gray} &  & \textbf{20} & 0.02 & \textit{(0.49)} & 1.21 & \textit{(2.87)} & 0.03 & \textit{(0.50)} & 1.06 & \textit{(1.73)} & 0.03 & \textit{(0.38)} & 0.54 & \textit{(0.52)} \\
\rowcolor{Gray} & & \textbf{100} & 0.00 & \textit{(0.19)} & 0.22 & \textit{(0.09)} & 0.00 & \textit{(0.18)} & 0.21 & \textit{(0.10)} & 0.01 & \textit{(0.16)} & 0.18 & \textit{(0.08)} \\
\rowcolor{Gray} & \multirow{-3}{*}{\rotatebox{90}{QLasso}} & \textbf{500} & 0.00 & \textit{(0.08)} & 0.17 & \textit{(0.05)} & 0.00 & \textit{(0.08)} & 0.08 & \textit{(0.02)} & 0.00 & \textit{(0.07)} & 0.07 & \textit{(0.02)} \\
\rowcolor{white} & \multicolumn{2}{c}{~} & \multicolumn{12}{c}{~}\\[-0.4cm]
\rowcolor{Gray2} &  & \textbf{20} & 0.01 & \textit{(0.56)} & 1.79 & \textit{(4.57)} & 0.01 & \textit{(0.54)} & 1.53 & \textit{(3.91)} & 0.02 & \textit{(0.41)} & 0.62 & \textit{(0.74)} \\
\rowcolor{Gray2} & & \textbf{100} & 0.01 & \textit{(0.18)} & 0.27 & \textit{(0.13)} & 0.00 & \textit{(0.19)} & 0.26 & \textit{(0.10)} & 0.00 & \textit{(0.17)} & 0.24 & \textit{(0.17)} \\
\rowcolor{Gray2} & \multirow{-3}{*}{\rotatebox{90}{ALasso}} & \textbf{500} & 0.00 & \textit{(0.08)} & 0.17 & \textit{(0.06)} & 0.00 & \textit{(0.08)} & 0.18 & \textit{(0.07)} & 0.00 & \textit{(0.08)} & 0.18 & \textit{(0.07)} \\
\hline\hline\end{tabular}
}	
	\end{center}
	\caption{\footnotesize The empirical performance of the quantile fused method (denoted as QLASSO), adaptive fused approach (ALASSO), and the standard LASSO approach (SLASSO) given for three different error distributions (standard normal, student's distribution with three degrees of freedom and the Cauchy distribution), three different sample sizes ($n \in \{20, 100, 500\}$), and three model selection techniques: a prior knowledge of three change-points in the model with the corresponding regularization parameter $\lambda_{(3)}$; the model given by the regularization parameter $\lambda_{AS} = n^{-1} (\log n)^{5/2}$ which satisfies the theoretical assumptions considered in this paper, and, finally, the model with the regularization parameter $\lambda_{MS}$ which minimizes the theoretical mean squared error. The models are compared with respect to the empirical bias defined as $n^{-1} \sum_{i = 1}^n (\boldsymbol{x}_i^\top \boldsymbol{\beta}_i^* - \boldsymbol{x}_i^\top \widehat{\boldsymbol{\beta}_i})$ and the empirical mean squared error (MSE) term $n^{-1} \sum_{i = 1}^n (\boldsymbol{x}_i^\top \boldsymbol{\beta}_i^* - \boldsymbol{x}_i^\top \widehat{\boldsymbol{\beta}_i})^2$, where $\boldsymbol{\beta}_i^*$ is the true value of the parameter and $\widehat{\boldsymbol{\beta}_i^*}$ is the corresponding estimated.  The values in the table are reported over 500 Monte Carlo simulations with the corresponding standard error values in brackets.}
	\label{tab2}
\end{table}

\begin{table}[!ht]\footnotesize
	\begin{center}\scalebox{0.8}{
		\begin{tabular}{ccc|cccccr<{\hspace{-\tabcolsep}}>{\hspace{-\tabcolsep}\,}lr<{\hspace{-\tabcolsep}}>{\hspace{-\tabcolsep}\,}lc<{\hspace{-\tabcolsep}}>{\hspace{-\tabcolsep}\,}l}
	\hiderowcolors
\multirow{2}{*}{$\boldsymbol{\mathcal{D}}$} & & \multirow{2}{*}{$\boldsymbol{n}$} & 
    $\boldsymbol{\lambda_{AS}}$ & \multicolumn{1}{c}{$\boldsymbol{\lambda_{(3)}}$} & \multicolumn{1}{c}{$\boldsymbol{\lambda_{MS}}$} & \multicolumn{2}{c}{\textbf{Number of Jumps}} & 
      \multicolumn{6}{c}{\textbf{Change-point Detection Error}} \\
 ~ & ~ & ~ & Value & Avg. & Avg. & $\lambda_{AS}$ & $\lambda_{MS}$ & \multicolumn{2}{c}{(Model  $\lambda_{(3)}$)} & \multicolumn{2}{c}{(Model $\lambda_{AS}$)} & \multicolumn{2}{c}{(Model  $\lambda_{CV}$)} \\\hline\hline
\multicolumn{3}{c}{~} & \multicolumn{11}{c}{~}\\[-0.4cm]
$\boldsymbol{N}$ &  & \textbf{20} & 0.78 & 3.77 & 4.02 & [0$|$0$|$9] & [0$|$0$|$11] & 0.11 & \textit{(0.03)} & 0.09 & \textit{(0.03)} & 0.09 & \textit{(0.04)} \\
  &  & \textbf{100} & 2.28 & 18.02 & 7.26 & [6$|$6$|$27] & [2$|$2$|$24] & 0.11 & \textit{(0.03)} & 0.03 & \textit{(0.02)} & 0.04 & \textit{(0.03)} \\
  & \multirow{-3}{*}{\rotatebox{90}{SLasso}} & \textbf{500} & 4.81 & 88.17 & 12.02 & [24$|$24$|$318] & [11$|$11$|$470] & 0.10 & \textit{(0.01)} & 0.00 & \textit{(0.00)} & 0.01 & \textit{(0.01)} \\
\multicolumn{3}{c}{~} & \multicolumn{11}{c}{~}\\[-0.4cm]
\rowcolor{Gray} & ~ & \textbf{20} & 0.78 & 0.91 & 1.35 & [0$|$0$|$8] & [0$|$0$|$9] & 0.12 & \textit{(0.03)} & 0.09 & \textit{(0.04)} & 0.09 & \textit{(0.03)} \\
\rowcolor{Gray}  &  & \textbf{100} & 2.28 & 3.67 & 2.41 & [1$|$1$|$14] & [1$|$1$|$26] & 0.11 & \textit{(0.03)} & 0.08 & \textit{(0.03)} & 0.07 & \textit{(0.03)} \\
\rowcolor{Gray}  & \multirow{-3}{*}{\rotatebox{90}{QLasso}} & \textbf{500} & 4.81 & 17.97 & 3.46 & [3$|$3$|$37] & [4$|$4$|$56] & 0.10 & \textit{(0.02)} & 0.03 & \textit{(0.02)} & 0.02 & \textit{(0.02)} \\
\rowcolor{white} & \multicolumn{2}{c}{~} & \multicolumn{11}{c}{~}\\[-0.4cm]
\rowcolor{Gray2} & ~ & \textbf{20} & 0.78 & 0.71 & 1.06 & [0$|$1$|$4] & [0$|$1$|$3] & 0.10 & \textit{(0.03)} & 0.10 & \textit{(0.03)} & 0.08 & \textit{(0.03)} \\
\rowcolor{Gray2}  &  & \textbf{100} & 2.28 & 1.69 & 1.74 & [0$|$1$|$4] & [0$|$2$|$5] & 0.09 & \textit{(0.03)} & 0.09 & \textit{(0.03)} & 0.09 & \textit{(0.03)} \\
\rowcolor{Gray2}  & \multirow{-3}{*}{\rotatebox{90}{ALasso}} & \textbf{500} & 4.81 & 4.30 & 3.73 & [0$|$2$|$6] & [0$|$3$|$6] & 0.09 & \textit{(0.04)} & 0.09 & \textit{(0.04)} & 0.09 & \textit{(0.04)} \\
\rowcolor{white} & \multicolumn{2}{c}{~} & \multicolumn{11}{c}{~}\\[-0.4cm]
\hline\multicolumn{3}{c}{~} & \multicolumn{11}{c}{~}\\[-0.4cm]
$\boldsymbol{t_3}$ &  & \textbf{20} & 0.78 & 5.16 & 7.66 & [1$|$1$|$12] & [0$|$0$|$8] & 0.12 & \textit{(0.03)} & 0.07 & \textit{(0.04)} & 0.09 & \textit{(0.04)} \\
  &  & \textbf{100} & 2.28 & 21.49 & 15.01 & [10$|$10$|$36] & [0$|$0$|$20] & 0.12 & \textit{(0.03)} & 0.02 & \textit{(0.01)} & 0.07 & \textit{(0.04)} \\
  & \multirow{-3}{*}{\rotatebox{90}{SLasso}} & \textbf{500} & 4.81 & 96.39 & 26.96 & [48$|$48$|$468] & [3$|$3$|$463] & 0.10 & \textit{(0.02)} & 0.00 & \textit{(0.00)} & 0.03 & \textit{(0.03)} \\
\multicolumn{3}{c}{~} & \multicolumn{11}{c}{~}\\[-0.4cm]
\rowcolor{Gray} & ~ & \textbf{20} & 0.78 & 0.88 & 1.50 & [0$|$0$|$9] & [0$|$0$|$12] & 0.12 & \textit{(0.04)} & 0.09 & \textit{(0.04)} & 0.09 & \textit{(0.04)} \\
\rowcolor{Gray}  &  & \textbf{100} & 2.28 & 3.47 & 2.53 & [0$|$0$|$19] & [0$|$0$|$26] & 0.11 & \textit{(0.03)} & 0.08 & \textit{(0.03)} & 0.08 & \textit{(0.03)} \\
\rowcolor{Gray}  & \multirow{-3}{*}{\rotatebox{90}{QLasso}} & \textbf{500} & 4.81 & 16.61 & 3.66 & [4$|$4$|$38] & [4$|$4$|$46] & 0.10 & \textit{(0.02)} & 0.04 & \textit{(0.02)} & 0.02 & \textit{(0.02)} \\
\rowcolor{white} & \multicolumn{2}{c}{~} & \multicolumn{11}{c}{~}\\[-0.4cm]
\rowcolor{Gray2} & ~ & \textbf{20} & 0.78 & 0.69 & 1.20 & [0$|$0$|$4] & [0$|$0$|$3] & 0.10 & \textit{(0.04)} & 0.09 & \textit{(0.03)} & 0.08 & \textit{(0.04)} \\
\rowcolor{Gray2}  &  & \textbf{100} & 2.28 & 1.59 & 1.75 & [0$|$1$|$4] & [0$|$1$|$5] & 0.10 & \textit{(0.03)} & 0.09 & \textit{(0.03)} & 0.09 & \textit{(0.03)} \\
\rowcolor{Gray2}  & \multirow{-3}{*}{\rotatebox{90}{ALasso}} & \textbf{500} & 4.81 & 3.96 & 3.66 & [0$|$1$|$6] & [0$|$2$|$7] & 0.10 & \textit{(0.04)} & 0.09 & \textit{(0.04)} & 0.09 & \textit{(0.04)} \\
\rowcolor{white} & \multicolumn{2}{c}{~} & \multicolumn{11}{c}{~}\\[-0.4cm]
\hline\multicolumn{3}{c}{~} & \multicolumn{11}{c}{~}\\[-0.4cm]
$\boldsymbol{C}$ &  & \textbf{20} & 0.78 & 20.88 & 39.42 & [2$|$2$|$19] & [0$|$0$|$13] & 0.13 & \textit{(0.04)} & 0.05 & \textit{(0.03)} & 0.11 & \textit{(0.04)} \\
  &  & \textbf{100} & 2.28 & 733.24 & 562.65 & [25$|$25$|$95] & [0$|$0$|$99] & 0.14 & \textit{(0.04)} & 0.01 & \textit{(0.01)} & 0.13 & \textit{(0.05)} \\
  & \multirow{-3}{*}{\rotatebox{90}{SLasso}} & \textbf{500} & 4.81 & 14316.02 & 5507.13 & [141$|$141$|$499] & [0$|$0$|$499] & 0.16 & \textit{(0.04)} & 0.00 & \textit{(0.00)} & 0.10 & \textit{(0.07)} \\
\multicolumn{3}{c}{~} & \multicolumn{11}{c}{~}\\[-0.4cm]
\rowcolor{Gray} & ~ & \textbf{20} & 0.78 & 0.86 & 2.06 & [0$|$0$|$10] & [0$|$0$|$8] & 0.12 & \textit{(0.04)} & 0.09 & \textit{(0.04)} & 0.10 & \textit{(0.03)} \\
\rowcolor{Gray}  &  & \textbf{100} & 2.28 & 3.21 & 2.63 & [0$|$0$|$72] & [0$|$0$|$86] & 0.12 & \textit{(0.03)} & 0.08 & \textit{(0.04)} & 0.08 & \textit{(0.04)} \\
\rowcolor{Gray}  & \multirow{-3}{*}{\rotatebox{90}{QLasso}} & \textbf{500} & 4.81 & 14.72 & 3.95 & [4$|$4$|$464] & [1$|$4$|$475] & 0.11 & \textit{(0.03)} & 0.04 & \textit{(0.02)} & 0.03 & \textit{(0.02)} \\
\rowcolor{white} & \multicolumn{2}{c}{~} & \multicolumn{11}{c}{~}\\[-0.4cm]
\rowcolor{Gray2} & ~ & \textbf{20} & 0.78 & 0.69 & 1.45 & [0$|$0$|$4] & [0$|$0$|$3] & 0.11 & \textit{(0.03)} & 0.10 & \textit{(0.03)} & 0.09 & \textit{(0.02)} \\
\rowcolor{Gray2}  &  & \textbf{100} & 2.28 & 1.61 & 2.02 & [0$|$0$|$4] & [0$|$0$|$5] & 0.10 & \textit{(0.03)} & 0.09 & \textit{(0.03)} & 0.09 & \textit{(0.03)} \\
\rowcolor{Gray2}  & \multirow{-3}{*}{\rotatebox{90}{ALasso}} & \textbf{500} & 4.81 & 4.51 & 3.90 & [0$|$1$|$8] & [0$|$1$|$8] & 0.10 & \textit{(0.04)} & 0.09 & \textit{(0.04)} & 0.08 & \textit{(0.04)} \\
\hline\hline\end{tabular}
}
	\end{center}
	\caption{\footnotesize The empirical change-point detection performance of the quantile fused method (QLASSO), adaptive fused approach (ALASSO), and the standard LASSO approach (SLASSO) for three error distributions, three sample sizes, and three model selection techniques (see the caption of Table \ref{tab2} for more details). The number of jumps is reported in the form ''$[M|M|M]$'' which stands for the minimum, median, and maximum number of change-points detected in the model over 500 Monte Carlo simulations. The change-point detection error is given as $\frac{1}{3} \sum_{k = 1}^{k} |\widehat{t}_k - t_k^*|$ and it is considered only for situations where  at least three change-points were detected in the model. The reported values are given as an average of such cases. The corresponding standard error values are reported in brackets.}
	\label{tab3}
\end{table}

As expected, it is obvious from the results reported in Tables \ref{tab2} and \ref{tab3} that while the standard LASSO fails for heavy-tailed error distributions, the proposed fused quantile approaches are  both still able to provide reliable and (asymptotically) consistent results (mainly with respect the parameter estimation performance). In addition, the adaptive fussed approach (denoted as ALASSO in the tables) seems to perform consistently even with respect to the change-point detection and it outperforms both, the standard LASSO approach (SLASSO) and the fused quantile method (QLASSO) as they both tend to select more change-points in the model, than the truth. The detection performance of the adaptive fused method seems to work under all three distributions (the median number of change-points being detected gets close to the true number of change-points if the sample size increases) however, the best results are observed under the normally distributed errors, and the slowest detection of the true change-points is observed under the Cauchy distributed errors. 
Both estimation methods proposed in this paper are shown to outperform the standard estimation techniques especially in situations where heavy-tailed error distributions are present. The standard LASSO property of overfitting the final model is evident for the fused quantile approach however,  the adaptive fused approach is able to overcome this problem and consistent asymptotic performance is empirically observed for both, the parameter estimation and the change-point detection as well. \enlargethispage{0.5cm}

\section{Regression example}
\label{regression_example}
The proposed methodology is also applied for a real data scenario: the same semi-synthetic stock example as in  \cite{Hyun-GSell-Tibshirani-18} is considered with $n = 251$ log daily returns simulated from a linear model based on three Dow Jones Industrial Average (DJIA) stocks. The true coefficient vector $\boldsymbol{\beta}_{i} \in \mathbb{R}^3$ is piece-wise constant in each of its element  with respect to $i \in \{ 1, \dots, n\}$ with three change-points located at the point $83$ (the first element of $\boldsymbol{\beta}_{i}$ changes from $-1$ to $1$), at the point $125$ (the second element of $\boldsymbol{\beta}_{i}$ changes from $-1$ to $1$), and finally, the third location at the point $166$ (where, again, the first element of $\boldsymbol{\beta}_{i}$ changes back from $1$ to its starting level of $-1$). There is no change-point with respect to the third element of the parameter vectors $\boldsymbol{\beta}_{i}$, for $i = 1, \dots, n$.

\begin{figure}
	\centering
	\subfigure[Quantile Fused Lasso] {\label{fig4:a}\includegraphics[width=0.48\textwidth]{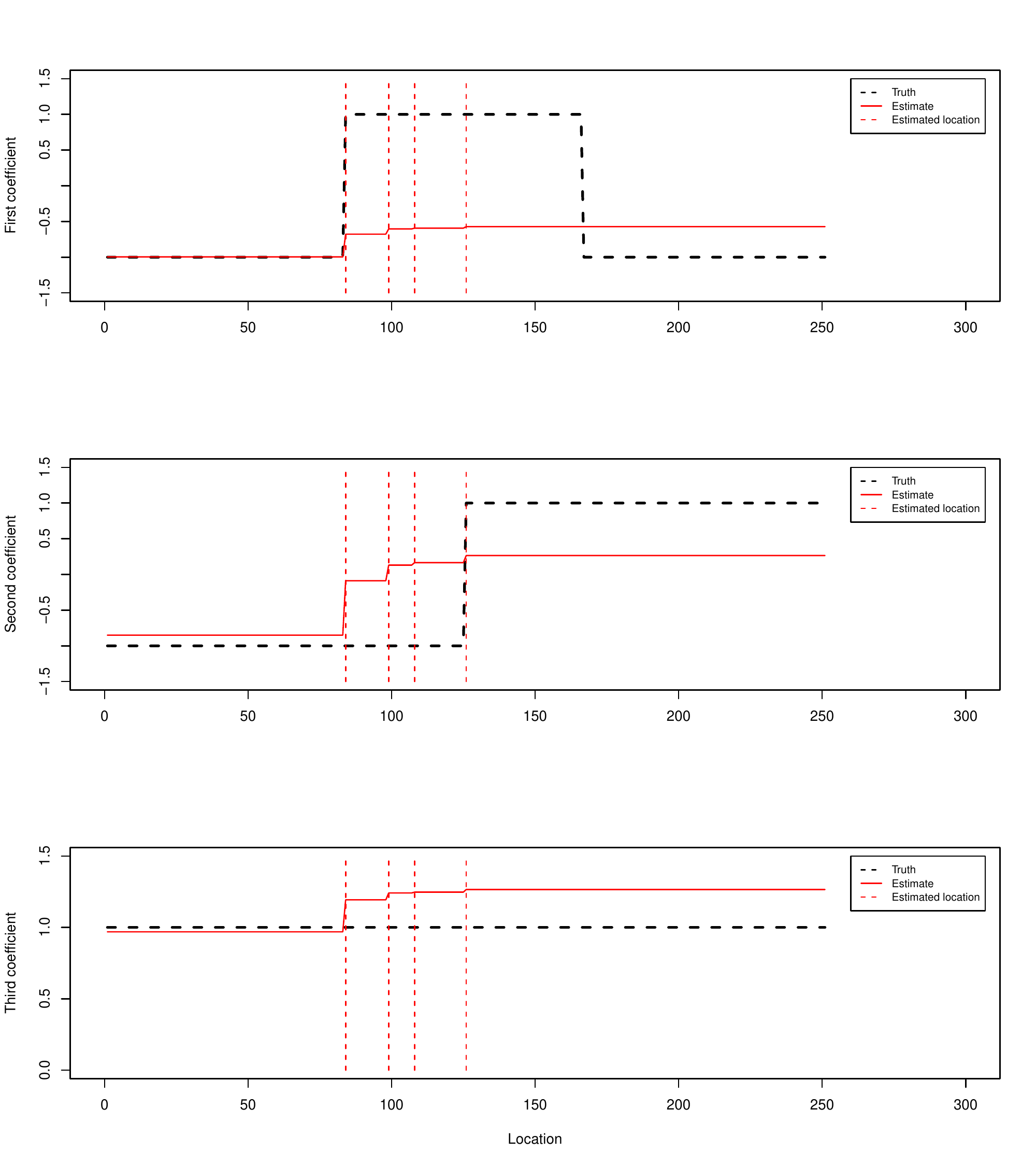}}
	\subfigure[Adaptive Quantile Fussed Lasso] {\label{fig4:b}\includegraphics[width=0.48\textwidth]{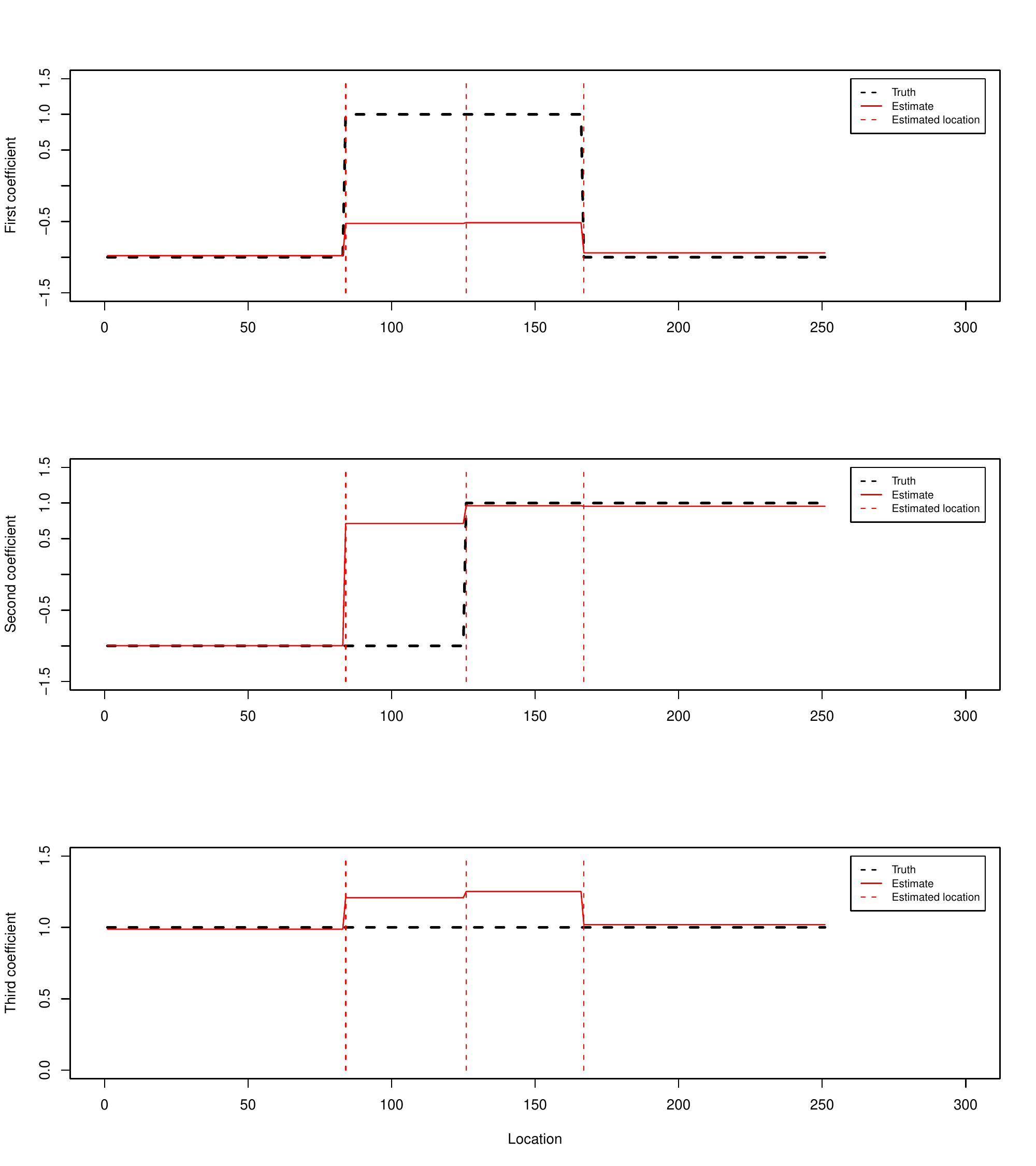}}
	
	\caption{\footnotesize An example of a semi-synthetic stock market discussed from  \cite{Hyun-GSell-Tibshirani-18}: the left panel shows the detection and estimation performance of the fused quantile method for $\tau = 0.5$. The true change-point locations are 83 and 166 for the first coefficient, there is one true change-point at the location 125 for the second coefficient and no change-point for the third coefficient. The estimated change-point locations are visualized by vertical dashed lines. The standard approach clearly overestimates the model with larger estimation bias while the adaptive quantile fused approach (right panel) detects all true and only true change-points while performing with a smaller estimation bias at the same time. }
	\label{fig4}
\end{figure}

The fussed LASSO approach used in \cite{Hyun-GSell-Tibshirani-18} with the 2-rise BIC stopping rule recovered 9 change-points in total (retaining 5 change-points after applying decluttering) with three of them being significant and roughly in a correspondence with the true change-point locations. For more details see \cite{Hyun-GSell-Tibshirani-18}. 
On the other hand, the proposed (group) quantile fused lasso  
revealed 4 change-point locations (detecting two true change-points and two false ones) while the adaptive quantile fused lasso approach correctly detected all three change-points with no false discoveries (see Figure \ref{fig4} for a comparison). However, the nature of the  group LASSO which is  used for the fused quantile LASSO and the adaptive quantile fused LASSO tends to estimate non-zero jumps for each element within a group and therefore, the detected locations in Figure \ref{fig4} always overlap across the three coefficients. In other words, the adaptive quantile fussed method correctly detects true change-point locations in the whole parameter vector $\boldsymbol{\beta}_{i} \in \mathbb{R}^3$ but it is not capable of specifying which element within the group is causing the change. This could be further improved, for instance, by adopting the idea of the sparse group LASSO approach in \cite{simon}.  From the overall point of view, however, the adaptive fused quantile LASSO clearly outperforms the fused quantile LASSO in both aspects, the change-point detection and the estimation bias, and it also seems to slightly outperform the method used in \cite{Hyun-GSell-Tibshirani-18}.

\section*{acknowledgements}
The authors are grateful to the associate editor and both reviewers for helping to improve the quality of the paper. We appreciate  their comments and suggestions. The authors also want to express their thanks to Sangwon Hyun for kindly providing us with the data for the example.  The work was partially supported by a bilateral grant between France and the Czech Republic provided by the PHC Barrande 2017 grant of Campus France (CG, grant number 38105NM) and the Ministry of Educations, Youth, and Sports in the Czech Republic (MM, Mobility grant 7AMB17FR030). The work of Mat\'u\v{s} Maciak was also supported by the Czech Science Foundation project GA\v{C}R, No. 18-00522Y.

\begin{appendices}
	\section{Proofs}
	\label{proofs}

\subsection{\textbf{Auxiliary lemmas and their proofs}}
We start with a straightforward result of the Karush-Kuhn-Tucker (KKT) optimality conditions for the quantile fused estimator defined by minimizing \eqref{eq3}, or \eqref{eq4} respectively.

\begin{lemma}
	\label{KKT}~\\
	(i) For any $l \in \{ 1, \cdots , |\overset{\vee}{\mathcal{A}}_n| \}$, all $ n \in \N$, and $\lambda_n >0$, it holds, with probability equal to one, that 
	$$\displaystyle{\tau \sum^n_{k=\overset{\vee}{t}_l} \ex_k-\sum^n_{k=\overset{\vee}{t}_l} \ex_k \e1_{\{Y_k \leq \ex^\top_k \overset{\vee}{ \eb}_k\}}=n \lambda_n \frac{\overset{\vee}{ \eth}_{\overset{\vee}{ t}_l}}{\| \overset{\vee}{\eth}_{\overset{\vee}{t}_l}\|}}.$$
	
	\noindent
	(ii) For any $j=1, \cdots n$,  all $n \in \N$, and $\lambda_n >0$, it holds, again with probability equal to one, that  $$\displaystyle{\left\|\tau \sum^n_{k=j} \ex_k-\sum^n_{k=j} \ex_k \e1_{\{Y_k \leq \ex^\top_k \overset{\vee}{ \eb}_k\}} \right\| \leq n \lambda_n}.$$
\end{lemma}  

The proof of Lemma \ref{KKT}  is similar with the proof of the following lemma and, therefore, it is omitted. The following lemma gives an analogous result, however for the adaptive fused quantile  estimator.

\begin{lemma}
	\label{aKKT}~\\
	(i) For any $l \in \{ 1, \cdots , |\widehat{\cal A}_n| \}$, all $ n \in \N$, and $\lambda_n >0$, it holds, with probability equal to one, that
	$$\displaystyle{\tau \sum^n_{k=\hat t_l} \ex_k-\sum^n_{k=\hat t_l} \ex_k \e1_{\{Y_k \leq \ex^\top_k \widehat{\eb}_k\}}=n \lambda_n \frac{\widehat{\eth}_{\hat t_l}}{\| \widehat{\eth}_{\hat t_l}\|}} \omega_{\hat t_l}.$$

	\noindent
	(ii) For any $j=1, \cdots n$,  all $n \in \N$, and $\lambda_n >0$, it holds, again with probability equal to one, that $$\displaystyle{\left\|\tau \sum^n_{k=j} \ex_k-\sum^n_{k=j} \ex_k \e1_{\{Y_k \leq \ex^\top_k \widehat{\eb}_k\}} \right\| \leq n \lambda_n} \omega_j.$$
\end{lemma}   

\noindent {\bf Proof of Lemma \ref{aKKT}}\\
We apply the  KKT optimality conditions for  $\widehat{\eth^n}$, which is the solution of  \eqref{aeq4} and by taking into account the fact that  $\sum^k_{i=1}\widehat{\eth}_i=\widehat{\eb}_k$, for any $k=1, \cdots ,n$, we obtain the assertion of the lemma. \hspace*{\fill}$\blacksquare$ \\

In the proofs of the theorems we will often use the following relation: for any vectors $\textbf{a}$, $\textbf{b}$, and $\textbf{c}$ which are of the same dimension, we have, by the triangular inequality, that 
\begin{equation}
\label{e4}
\textrm{ if }\| \textbf{a}+\textbf{b}\| \leq \| \textbf{c} \|  \textrm{ and } \| \textbf{b} \| \leq \|\textbf{c}\|, \textrm{ then } \|\textbf{a}\| \leq 2 \|\textbf{c}\|.
\end{equation}

\begin{lemma}
	\label{Lemma}~\\
	For any real  random vectors  $\textbf{A}$ and $\textbf{B}$ of the same dimension, any real value $x>0$, such that $\PP[\|\textbf{A}+\textbf{B}\| \leq x]=1$, we have that
	\begin{equation*}
	\label{AB}
	1\leq \PP\bigg[v x \geq  \|\textbf{A}\|  \bigg] +\PP\bigg[v \|\textbf{B}\| \geq (v-1)\|\textbf{A}\| \bigg],
	\end{equation*} 
	which holds for any constant  $v>1$.
\end{lemma}

\noindent {\bf Proof of Lemma \ref{Lemma}}.\\
We have: $1=\PP\big[v x \geq  \|\textbf{A}\|  \big]+\PP\big[vx <  \|\textbf{A}\|  \big]$. Taking into account:  $\PP[\|\textbf{A}+\textbf{B}\| \leq x] \leq \PP[\|\textbf{A}\|-\|\textbf{B}\| \leq x] $, we obtain that: $\PP\big[vx <  \|\textbf{A}\|  \big] = \PP\big[\big\{vx <  \|\textbf{A}\| \big\} \cap \big\{\|\textbf{A}+\textbf{B}\| \leq x \big\} \big] \leq \PP\big[\big\{vx <  \|\textbf{A}\|  \big\} \cap \big\{\|\textbf{A}\|-\|\textbf{B}\| \leq x \big\} \big]=\PP \big[v \|\textbf{B}\| \geq (v-1)\|\textbf{A}\| \big]$. 
\hspace*{\fill}$\blacksquare$ \\

The following lemma will be used to control the supremum of the averaged value of the random quantities  $\ex_i \big(\e1_{\{\varepsilon_i \leq u\}}-F(u) \big)$, for $i = 1, \dots, n$, and any $u \in \R$.

\begin{lemma}
	\label{Lemma 4}~\\
	Under Assumption (A1) imposed on the model design, Assumption (A2) for errors $(\varepsilon_i)_{1 \leqslant i \leqslant n}$, and  two positive sequences $(v_n)$, $(z_n)$ such that $v_n z^2_n (\log n)^{-1} {\underset{n \rightarrow \infty}{\longrightarrow} } \infty$, it holds that
	\[
	\lim_{n \rightarrow \infty} \PP\bigg[ \max_{
		\substack{
			1 \leq r_n< s_n \leq n\\
			s_n-r_n \geq v_n
	}} 
	\;
	\sup_{u \in \R}
	\bigg\| \frac{1}{s_n-r_n} \sum ^{s_n-1}_{i=r_n} \ex_i \big(\e1_{\{\varepsilon_i \leq u\}}-F(u) \big)\bigg\| \geq z_n
	\bigg]  = 0,
	\]
	where $F$ stands for the distribution function of the error terms $\varepsilon_i$, for $i = 1, \dots, n$.
\end{lemma}

\noindent {\bf Proof of Lemma \ref{Lemma 4}}.\\ 
We have that 
\begin{align}
\label{tx1}
\PP\bigg[ \max_{\substack{1 \leq r_n< s_n \leq n\\s_n-r_n \geq v_n}} \; \sup_{u \in \R}\bigg\| \frac{1}{s_n-r_n} & \sum ^{s_n-1}_{i=r_n} \ex_i \big(\e1_{\{\varepsilon_i \leq u\}}-F(u) \big)\bigg\| \geq z_n \bigg] \\
& \leq \sum_{\substack{1 \leq r_n< s_n \leq n\\s_n-r_n \geq v_n}} \PP\bigg[\sup_{u \in \R} \bigg\| \frac{1}{s_n-r_n} \sum ^{s_n-1}_{i=r_n} \ex_i \big(\e1_{\{\varepsilon_i \leq u\}}-F(u) \big)\bigg\| \geq z_n  \bigg].
\end{align}
Since for any $\ex \in \R^p$, it holds that $\|\ex \| \leq \sqrt{p} \| \ex\|_\infty$, we also have 
\begin{align*}
\PP & \bigg[\sup_{u \in \R} \bigg\| \frac{1}{s_n-r_n}  \sum ^{s_n-1}_{i=r_n} \ex_i \big(\e1_{\{\varepsilon_i \leq u\}}-F(u) \big)\bigg\| \geq z_n  \bigg]  \\
&\leq \PP\bigg[\sup_{u \in \R} \bigg\| \frac{1}{s_n-r_n} \sum ^{s_n-1}_{i=r_n} \ex_i \big(\e1_{\{\varepsilon_i \leq u\}}-F(u) \big)\bigg\|_\infty \geq p^{-1/2}z_n  \bigg].
\end{align*}
Let $x_{ij}$ be the $j$-th component of $\ex_i$, for $j=1, \cdots, p$. 

We use the  Hoeffding's  inequality   for independent  random variables $x_{ij} \e1_{\{\varepsilon_i \leq t\}}$, for $j=1, \cdots , p$, and we obtain, that for all $u \in \mathbb{R}$, it holds that 
\[
\PP\bigg[\bigg| \frac{1}{s_n-r_n} \sum^{s_n-1}_{i=r_n} x_{ij}(\e1_{\{\varepsilon_i \leq u\}}-F(u)) \bigg| \geq p^{-1/2}z_n  \bigg] \leq 2 \exp\big(-2   C (s_n-r_n)p^{-1}z^2_n \big),
\]
since $x_{ij}$ is bounded by  Assumption (A1). Taking into account the relation in \eqref{tx1} and the condition where $v_n z^2_n (\log n)^{-1} {\underset{n \rightarrow \infty}{\longrightarrow} } \infty$, we finally obtain
\begin{align*}
\PP&\bigg[ \max_{\substack{1 \leq r_n< s_n \leq n\\s_n-r_n \geq v_n}} \; \sup_{u \in \R}\bigg\| \frac{1}{s_n-r_n} \sum ^{s_n-1}_{i=r_n} \ex_i \big(\e1_{\{\varepsilon_i \leq u\}}-F(u) \big)\bigg\| \geq  p^{-1/2} z_n \bigg] \\
&\leq 2 n^2 \exp\big(-2   C (s_n-r_n)p^{-1}z^2_n \big) {\underset{n \rightarrow \infty}{\longrightarrow} } 0,
\end{align*}
which completes the proof of the lemma. \hspace*{\fill}$\blacksquare$ \\

Analogously to Lemma \ref{Lemma 4} we can also formulate the next lemma which controls the supremum of the average of the random quantities  $\ex_i \big(\e1_{\{\varepsilon_i \leq u_i\}}-F(u_i) \big)$, for $i = 1, \dots, n$ and any real $u_i \in{\cal B}$, where ${\cal B}$ is a bounded set in $\R$. The proof of the lemma follows similar lines as the proof above and, therefore, it is omitted.

\begin{lemma}~\\
	\label{Lemma 4bis}
	Under Assumption (A1) imposed on the model design, Assumption (A2) for errors $(\varepsilon_i)_{1 \leqslant i \leqslant n}$, two positive sequences $(v_n)$, $(z_n)$  such that $v_n z^2_n (\log n)^{-1} {\underset{n \rightarrow \infty}{\longrightarrow} } \infty$, and the distribution function $F$ of the error terms for which the limit $\lim_{n \rightarrow \infty}\big(n^{-1}\sum^n_{i=1} F(u_i))$ exists, it holds that
	\[
	\lim_{n \rightarrow \infty}\PP\bigg[ \max_{
		\substack{
			1 \leq r_n< s_n \leq n\\
			s_n-r_n \geq v_n
	}} 
	\;
	\sup_{u_i \in {\cal B}}
	\bigg| \frac{1}{s_n-r_n} \sum ^{s_n-1}_{i=r_n}  \big(\e1_{\{\varepsilon_i \leq u_i\}}-F(u_i) \big)\bigg| \geq z_n
	\bigg]  = 0,
	\]
	where $u_i \in \mathcal{B}$ and ${\cal B}$ is some bounded set in $\R$.
\end{lemma}

\subsection{\textbf{Proofs of Theorems and Corollary \ref{Raq_vit}}}

\noindent {\bf Proof of Theorem  \ref{aTheorem 3.1}}

Consider the random event $V_{n,k} \equiv \{|\hat t_k - t^*_k| \geq n \delta_n \}$, for $k \in \{1, \cdots , K^*\}$. Without any loss of generality, it is assumed that $t^*_k > \hat t_k$. The  case where $t^*_k < \hat t_k$ can be treated analogously. Since  $K^*< \infty$ by Assumption (A6), the theorem is immediately proved if we show that for any $k =1, \cdots, K^*$ we have, that $\lim_{n \rightarrow \infty} \PP [V_{n,k}]=0$.

To show that the limit equals to zero, we consider the following decomposition: $V_{n,k}=(V_{n,k} \cap W_n) \cup(V_{n,k} \cap \overline W_n)$, with   $W_n \equiv \{\max_{1 \leq k \leq K^*}|\hat t_k - t^*_k| <I^*_{min}/2  \}$ and $ \overline W_n$ being the complement of $W_n$. Thus, we obtain
\begin{equation}
\label{eq_mm1}
\PP[V_{n,k}]=\PP[V_{n,k} \cap W_n] + \PP[V_{n,k} \cap \overline W_n],
\end{equation}
and we can deal with both terms on the right side of \eqref{eq_mm1} separately.

\noindent\textit{(i)} Let us firstly study $\PP[V_{n,k} \cap W_n]$. More precisely,  we need to show that
\begin{equation}
\label{eq21}
\lim_{n \rightarrow \infty} \PP[V_{n,k} \cap W_n]=0.
\end{equation}
Suppose, again without any loss of generality, that $t^*_{k-1} \leq \hat t_k < t^*_k$, where $t^*_0=1$. Then, applying Lemma \ref{aKKT} for $j=t^*_k$ and $l=k$, we obtain that
\[
\left\|\tau \sum^n_{i=t^*_k} \ex_i-\sum^n_{i=t^*_k} \ex_i \e1_{\{Y_i \leq \ex_i^\top \widehat{\eb}_i\}} \right\| \leq n \lambda_n \omega_{t^*_k}
\]
holds with probability equal to one, and, similarly also
\[
\tau \sum^n_{i=\hat t_k} \ex_i-\sum^n_{i=\hat t_k} \ex_i \e1_{\{Y_i \leq \ex^\top_i \widehat{\eb}_i\}}=n \lambda_n \frac{\widehat{\eth}_{\hat t_k}}{\| \widehat{\eth}_{\hat t_k}\|} \omega_{\hat t_k},
\]
which again holds with probability equal to one. The last relation can be also rewritten as
\[
\tau \sum^{t^*_k-1}_{i=\hat t_k} \ex_i - \sum^{t ^*_k-1}_{i=\hat t_k} \ex_i \e1_{\{Y_i \leq \ex_i^\top \widehat{\eb}_i\}}+ \tau \sum^n_{i=t^*_k} \ex_i - \sum^n_{i=t^*_k} \ex_i \e1_{\{Y_i \leq \ex_i^\top \widehat{\eb}_i\}}= n \lambda_n \frac{\widehat{\eth}_{\hat t_k}}{\| \widehat{\eth}_{\hat t_k}\|}\omega_{\hat t_k}, 
\]
and we can directly apply the relation in \eqref{e4} for 
$
\textbf{a} =  \tau \sum^{t^*_k-1}_{i=\hat t_k} \ex_i -  \sum^{t^*_k-1}_{i=\hat t_k} \ex_i \e1_{\{Y_i \leq \ex_i^\top \widehat{\eb}_{i}\}}$,  
$ \textbf{b}  =  \tau \sum^n_{i=t^*_k} \ex_i - \sum^n_{i=t^*_k} \ex_i  \e1_{\{Y_i \leq \ex_i^\top \widehat{\eb}_{i}\}}$, and $c  =  n \lambda_n\big(\omega_{t^*_k} +\omega_{\hat t_k}\big)$,
to obtain that  
\[   
\PP \bigg[\Big\| \tau   \sum^{t^*_k-1}_{i=\hat t_k} \ex_i -  \sum^{t^*_k-1}_{i=\hat t_k} \ex_i \e1_{\{Y_i \leq \ex_i^\top \widehat{\eb}_{i}\}}  \Big\| \leq 2  n \lambda_n\big(\omega_{t^*_k} +\omega_{\hat t_k}\big) \bigg]=1.
\]  
We now use Lemma \ref{Lemma}, for $x=2n \lambda_n\big(\omega_{t^*_k} +\omega_{\hat t_k}\big)$, and some constant $v > 1$, such that \[
v>\frac{c^{(0)}\max(\tau, L_k)}{|c^{(0)}\max(\tau, L_k) - c^{(1)}\min(\tau, L_k) |} > 1,
\] where the real random vectors $\textbf{A}$ and $\textbf{B}$ are defined as follows:
\begin{itemize}
	\item  $\textbf{A} =\sum^{t^*_k-1}_{i=\hat t_k} \ex_i \e1_{\{Y_i \leq \ex_i^\top \widehat{\eb}_{i}\}} $ and $\textbf{B}= \tau \sum^{t^*_k-1}_{i=\hat t_k} \ex_i $, if $\tau < L_k$;
	\item  $\textbf{A}= \tau \sum^{t^*_k-1}_{i=\hat t_k} \ex_i $,
	and $\textbf{B}=\sum^{t^*_k-1}_{i=\hat t_k} \ex_i \e1_{\{Y_i \leq \ex_i^\top \widehat{\eb}_{i}\}} $, if $\tau > L_k$.
\end{itemize}
Now, the probability $\PP [V_{n,k} \cap W_n]$ can be expressed as 
\begin{align}
\PP [V_{n,k} \cap W_n] & =\PP\left[V_{n,k} \cap W_n \cap \left\{   
\left\| \tau   \sum^{t^*_k-1}_{i=\hat t_k} \ex_i -  \sum^{t^*_k-1}_{i=\hat t_k} \ex_i \e1_{\{Y_i \leq \ex_i^\top \widehat{\eb}_{i}\}}  \right\| \leq  2 n \lambda_n \big(\omega_{t^*_k} +\omega_{\hat t_k}\big)\right\}\right] \nonumber \\
& \textcolor{red}{\leq} {\cal P}_{1,k,n}+{\cal P}_{2,k,n}, \label{eq_mm1}
\end{align}
with ${\cal P}_{1,k,n} \equiv \PP \big[ \big\{  \|\textbf{A}\|  \leq v x \big\} \cap V_{n,k} \cap W_n \big] $, and ${\cal P}_{2,k,n} \equiv \PP \big[ \big\{v\|\textbf{B}\| \geq  (v-1) \|\textbf{A}\| \big\} \cap V_{n,k} \cap W_n \big] $.

In order to deal with these two probabilities, it is necessary to know the convergence rate of the estimator $\widehat{\ef}_{k+1}$ of $\ef_{k+1}$ obtained by minimizing (\ref{eq3}), knowing that $|\widehat{\cal A}_n|=K^*$, and that the random events $W_n$ and $V_{n,k}$ both occur. Therefore, in the following,  we study the convergence rate of $\widehat{\ef}_{k+1}$ and afterwards  we return back to study ${\cal P}_{1,k,n}$ and   ${\cal P}_{2,k,n}$.\\

\noindent \underline{\it Convergence rate of $\widehat{\ef}_{k+1}$}\\
Since the random event $W_n$ occurs, it is supposed, without any loss of generality, that we are in the following case: $t^*_{k-1}  < \hat t_k < t^*_k < (t^*_k+t^*_{k+1})/2 < \hat t_{k+1} \leq t^*_{k+1}$, with $t^*_k - \hat t_k \geq n \delta_n  $. Let us recall the convention where $t^*_0=1$ and $t^*_{K+1}=n$. By applying Lemma \ref{aKKT}  for $j=(t^*_k+t^*_{k+1})/2$ and $j=t^*_k$, and using the relation from \eqref{e4}, we get that 
\begin{align}
2 n \lambda_n \big(\omega_{t^*_k}+\omega_{(t^*_k+t^*_{k+1})/2} \big)& \geq \left\| \tau \sum^{(t^*_k+t^*_{k+1})/2}_{i=t^*_k}\ex_i - \sum^{(t^*_k+t^*_{k+1})/2}_{i=t^*_k}\ex_i \e1_{\{\varepsilon_i \leq \ex_i^\top(\widehat{\ef}_{k+1}-\ef^*_k)\}} \right\|  \nonumber \\
&\geq \left\| \sum^{(t^*_k+t^*_{k+1})/2}_{i=t^*_k}\ex_i \e1_{\{\varepsilon_i \leq \ex_i^\top(\widehat{\ef}_{k+1}-\ef^*_k)\}}\right\| - \tau \left\| \sum^{(t^*_k+t^*_{k+1})/2}_{i=t^*_k}\ex_i \right\|,\label{ac}
\end{align}
with probability equal to one.  Thus, again with probability equal to one, it holds that 
\begin{equation}
\label{A1}
\left\| \sum^{(t^*_k+t^*_{k+1})/2}_{i=t^*_k}\ex_i \e1_{\{\varepsilon_i \leq \ex_i^\top(\widehat{\ef}_{k+1}-\ef^*_k)\}}\right\|   \leq \tau \left\| \sum^{(t^*_k+t^*_{k+1})/2}_{i=t^*_k}\ex_i \right\|  + 2 n \lambda_n \big(\omega_{t^*_k}+\omega_{(t^*_k+t^*_{k+1})/2} \big) .
\end{equation}
Now, from relation \eqref{ac}, we deduce that
\[
2 n \lambda_n \big(\omega_{t^*_k}+\omega_{(t^*_k+t^*_{k+1})/2} \big) \geq \tau \left\| \sum^{(t^*_k+t^*_{k+1})/2}_{i=t^*_k}\ex_i \right\|  - \left\| \sum^{(t^*_k+t^*_{k+1})/2}_{i=t^*_k}\ex_i \e1_{\{\varepsilon_i \leq \ex_i^\top(\widehat{\ef}_{k+1}-\ef^*_k)\}}\right\|,
\]
which holds with probability one and, also, again with probability equal to one, we have
\begin{equation}
\label{A2}
\left\| \sum^{(t^*_k+t^*_{k+1})/2}_{i=t^*_k}\ex_i \e1_{\{\varepsilon_i \leq \ex_i^\top(\widehat{\ef}_{k+1}-\ef^*_k)\}}\right\| \geq \tau \left\| \sum^{(t^*_k+t^*_{k+1})/2}_{i=t^*_k}\ex_i \right\| - 2 n \lambda_n \big(\omega_{t^*_k}+\omega_{(t^*_k+t^*_{k+1})/2} \big).
\end{equation}

On the other hand, using the condition in \eqref{CD1}, we have
\begin{align} 
\left\| \sum^{(t^*_k+t^*_{k+1})/2}_{i=t^*_k}\ex_i \right\| &  \geq \left\| \sum^{(t^*_k+t^*_{k+1})/2}_{i=t^*_k}\ex_i \right\|_\infty \geq | c^{(0)}| \left(\frac{t^*_k+t^*_{k+1}}{2} - t^*_k \right) =| c^{(0)}|\frac{t^*_{k+1}-t^*_k}{2} \nonumber \\
& \geq  | c^{(0)}|\frac{n \delta_n}{2}\gg n \lambda_n \big(\max(b_n,d_n) \big)^{-\gamma} . \label{acc}
\end{align}
The relation in \eqref{acc}, together with \eqref{A1} and \eqref{A2} now imply that
\begin{equation}
\label{A3}
\left\| \sum^{(t^*_k+t^*_{k+1})/2}_{i=t^*_k}\ex_i \e1_{\{\varepsilon_i \leq \ex_i^\top(\widehat{\ef}_{k+1}-\ef^*_k)\}}\right\| = \tau \left\| \sum^{(t^*_k+t^*_{k+1})/2}_{i=t^*_k}\ex_i \right\|  +O_{\PP}\bigg( n \lambda_n\big(\max(b_n,d_n) \big)^{-\gamma} \bigg)  .
\end{equation}
Considering \eqref{acc}, we deduce that relation \eqref{A3} can be also expressed as  
\begin{equation}
\label{A4}
\sum^{(t^*_k+t^*_{k+1})/2}_{i=t^*_k}\ex_i \left(\e1_{\{\varepsilon_i \leq \ex_i^\top(\widehat{\ef}_{k+1} - \ef^*_{k+1})\}} - \tau\right)= O_{\PP}\bigg( n \lambda_n\big(\max(b_n,d_n) \big)^{-\gamma} \bigg)\eu ,
\end{equation}
for some  vector $\eu \in \R^p$, such that $\| \eu   \|=1$. The left side of \eqref{A4} can be further rewritten as 
\begin{align}
\label{A4bis}
& \sum^{(t^*_k+t^*_{k+1})/2}_{i=t^*_k}\ex_i \left(\e1_{\{\varepsilon_i \leq \ex_i^\top(\widehat{\ef}_{k+1} - \ef^*_{k+1})\}} - \tau\right) \\
& \hskip0.5cm = \sum^{(t^*_k+t^*_{k+1})/2}_{i=t^*_k} \ex_i \bigg(\e1_{\{\varepsilon_i \leq \ex_i^\top(\widehat{\ef}_{k+1} - \ef^*_{k+1})\}} -F(\ex_i^\top(\widehat{\ef}_{k+1} - \ef^*_{k+1}))+F(\ex_i^\top(\widehat{\ef}_{k+1} - \ef^*_{k+1})) - F(0) \bigg)\nonumber\\
& \hskip0.5cm = O_{\PP}\bigg( n \lambda_n\big(\max(b_n,d_n) \big)^{-\gamma} \bigg) \eu.\nonumber
\end{align}
The estimator  $\widehat{\ef}_{k+1}$ can be considered as $\|\widehat{\ef}_{k+1} - \ef^*_{k+1}\| \leq c_n  $ with probability converging to one, where $c_n$ is some deterministic sequence to be determined later. Let us define $a_{i,n} \equiv \ex_i^\top c_n \textbf{w}  $, for some vector $\textbf{w} \in \R^p$, such that  $\| \textbf{w} \| < \infty$. We prove now that the sequence $\{ a_{i,n}\}_n$ is bounded for any $i$, if $n$ is large enough. In contrary, let us assume that $\{ a_{i,n}\}_n$ is not bounded. Thus, a subsequence which converges either to  $+ \infty$ or $- \infty$ can be selected. Suppose that the subsequence $(a_{i,n_m})_{m \geq 1}$ converges to infinity. Hence, relation \eqref{A4} reduces to
\[
\sum^{(t^*_k+t^*_{k+1})/2}_{i=t^*_k}\ex_i \big(\e1_{\{\varepsilon_i \leq a_{i,n_m}\}} - \tau\big)= O_{\PP}\bigg( n_m \lambda_{n_m} \big(\max(b_{n_m},d_{n_m}) \big)^{-\gamma}\bigg)\eu,
\]
which implies, since $F(0)=\tau$ and $0<|c^{(0)}|<1$, that
$$\sum^{(t^*_k+t^*_{k+1})/2}_{i=t^*_k}c^{(0)} = O \bigg( n_m \lambda_{n_m} \big(\max(b_{n_m},d_{n_m}) \big)^{-\gamma}\bigg).$$ 
Therefore $c^{(0)}(t^*_{k+1}- t^*_k)/2=O \bigg( n_m \lambda_{n_m} \big(\max(b_{n_m},d_{n_m}) \big)^{-\gamma}\bigg)=o(n \delta_n)$ by condition \eqref{CD1} which is in a contradiction with Assumption (A4), since $c^{(0)} \neq 0$. Thus, there exists a constant $C>0$, such that $|a_{i,n}| <C$, for any $i \in \{1, \cdots , n \}$, for $n$ large enough. Next, we have that $\eE[\ex_i \e1_{\{\varepsilon_i \leq a_{i,n}\}}]= \ex_i F(a_{i,n})$ and $Var[\ex_i \e1_{\varepsilon_i \leq a_{i,n}}]= \ex_i \ex^\top_i F(a_{i,n}) (1-F(a_{i,n}))$ and we express  $\widehat{\ef}_{k+1}$ in the form $\ef^*_{k+1}+C c_n \textbf{w}$.  Using the Central Limit Theorem (CLT) together with Assumption (A1) and the fact that $a_{i,n}$ is bounded, we have that 
$$\sum^{(t^*_k+t^*_{k+1})/2}_{i=t^*_k} \ex_i  \big(\e1_{\{\varepsilon_i \leq \ex_i^\top c_n \textbf{w}\}} -F(a_{i,n}) \big)=O_{\PP}\big( ( t^*_{k+1} - t^*_k)^{1/2}  \big),$$
and also 
$$\sum^{(t^*_k+t^*_{k+1})/2}_{i=t^*_k}  \big( F(\ex_i^\top c_n \textbf{w})-F(0)  \big) =\sum^{(t^*_k+t^*_{k+1})/2}_{i=t^*_k}  c_n \ex_i \ex_i^\top f(b_{i,n}),$$ 
where $b_{i,n}$ is some value between zero and $a_{i,n}$. Then the relation in \eqref{A4bis} implies that for $n$  large enough we have
\begin{equation}
\label{mc}\left\|c_n \sum^{(t^*_k+t^*_{k+1})/2}_{i=t^*_k} \ex_i \ex_i^\top f(b_{i,n}) \textbf{w} \right\| \leq C \bigg( n \lambda_n \big(\max(b_n,d_n) \big)^{-\gamma} + (t^*_{k+1}-t^*_k)^{1/2}\bigg),
\end{equation}
due to the fact that the density function $f$ is bounded by Assumption (A2). 

Let us  define a positive definite matrix $\textbf{D} \equiv 2(t^*_{k+1}- t^*_k)^{-1} \sum^{(t^*_k+t^*_{k+1})/2}_{i=t^*_k} \ex_i \ex_i^\top f(b_{i,n})$. Using the matrix property where $\| \textbf{DC} \|= [tr (\textbf{CC}^\top \textbf{DD}^\top) ]^{1/2} \geq \mu_{min} (\textbf{D}^\top\textbf{D})^{1/2}\| \textbf{C}\|$, for $\textbf{C} = c_n \textbf{w}$, we have by Assumptions (A2) and (A3), that
$$c_n \|\textbf{w}\| \cdot \left\|\sum^{(t^*_k+t^*_{k+1})/2}_{i=t^*_k} \ex_i \ex_i^\top f(b_{i,n}) \right\| \geq m_0 c_n (t^*_{k+1} -t^*_k) \|\textbf{w}\|.$$ 
Then, taking into account  the  relation in \eqref{mc}, we obtain, for $n$  large enough, that
\begin{align*}
m_0 c_n(t^*_{k+1} -t^*_k) \|\textbf{w}\|  \leq n \lambda_n \big(\max(b_n,d_n) \big)^{-\gamma} +(t^*_{k+1}-t^*_k)^{1/2},
\end{align*} 
and also
\begin{align*}
c_n \|\textbf{w}\| \leq \frac{n \lambda_n \big(\max(b_n,d_n) \big)^{-\gamma}}{t^*_{k+1} -t^*_k} +(t^*_{k+1}-t^*_k)^{-1/2} & \leq \frac{n \lambda_n \big(\max(b_n,d_n) \big)^{-\gamma}}{I^*_{min}} +(I^*_{min})^{-1/2} \\ & \leq  \frac{n \lambda_n \big(\max(b_n,d_n) \big)^{-\gamma}}{n \delta_n} +(n \delta_n)^{-1/2}.
\end{align*}

Therefore, given the random event $W_n \cap V_{n,k}$, Assumption (A4) together with the condition in \eqref{CD1}, and the fact that $| \widehat{\cal A}_n|=K^* < \infty$, we have that $\widehat{\ef}_{k+1}$ converges to $\ef^*_{k+1}$ at the rate of order $ {n \lambda_n\big(\max(b_n,d_n) \big)^{-\gamma}}/{I^*_{min}} +(I^*_{min})^{-1/2}$.
Let us denote by $(c_n)$ the  sequence
\begin{equation}
\label{cn}
c_n \equiv   {n \lambda_n \big(\max(b_n,d_n) \big)^{-\gamma}}/{I^*_{min}} +(I^*_{min})^{-1/2}.
\end{equation} 
Now, due to Assumption (A4) and the condition in \ref{CD1}, we get that $ c_n \rightarrow 0$, as $n \rightarrow \infty$.  

Now we return back to study ${\cal P}_{1,k,n}$ and ${\cal P}_{2,k,n}$ from \eqref{eq_mm1} and we consider two separate cases here: it either holds that $\tau > L_k$, or $\tau < L_k$. 
Let us start with the situation where $\tau > L_k$ and we consider the first probability term 
$${\cal P}_{1,k,n}= \PP \left[ \bigg\{ v^{-1}\bigg\|\tau \sum^{t^*_k-1}_{i=\hat t_k} \ex_i\bigg\| \leq 2n \lambda_n\big(\omega_{t^*_k} +\omega_{\hat t_k}\big) \bigg\} \cap V_{n,k} \cap W_n \right].$$

Since $\ex_i=(c^{(0)}, x_{2i}, \cdots , x_{pi})^\top$, then
$
\big\| \sum^{t^*_k-1}_{i=\hat t_k} \ex_i  \big\|_{\infty} \geq \sum^{t^*_k-1}_{i=\hat t_k} c^{(0)}=c^{(0)}(t^*_k - \hat t_k)$ and since $\ex_i\in \R^p$, with $p$ not depending on $n$ and $\|\ex \| \geq \|\ex\|_{\infty}$, then
\[
\PP\bigg[\tau \bigg\| \sum^{t^*_k-1}_{i=\hat t_k} \ex_i\bigg\| \leq 2 v n \lambda_n\big(\omega_{t^*_k} +\omega_{\hat t_k}\big) \bigg] \leq \PP\bigg[\tau c^{(0)}(t^*_k - \hat t_k) \leq 2v n \lambda_n\big(\omega_{t^*_k} +\omega_{\hat t_k}\big)\bigg] ,
\]
with $\omega_{t^*_k}=\big( \max \big(\|\overset{\smile}{ \eb}_{t^*_k}- \overset{\smile}{ \eb}_{t^*_k-1} \|_\infty,d_n \big)\big)^{-\gamma}$ and $\omega_{\hat t_k}=\big( \max \big(\|\overset{\smile}{ \eb}_{\hat t_k}- \overset{\smile}{ \eb}_{\hat t_k-1} \|_\infty,d_n \big)\big)^{-\gamma}$. From the convergence rate of the quantile fused estimators given by Remark \ref{Rq_vit}, taking also into account Assumption (A5), we have, if $t^*_k$ is equal to $\breve{t}$, that $\|\overset{\smile}{ \eb}_{t^*_k}- \overset{\smile}{ \eb}_{t^*_k-1}\|_\infty >C>0$, and, if $t^*_k$ is not $\breve{ t}$,  then $\|\overset{\smile}{ \eb}_{t^*_k}- \overset{\smile}{ \eb}_{t^*_k-1}\|_\infty =Cb_n$.  An analogous situation also applies for $\|\overset{\smile}{ \eb}_{\hat t_k}- \overset{\smile}{ \eb}_{\hat t_k-1}\|_\infty$.  
If $\omega_{t^*_k} +\omega_{\hat t_k} >C>0$, then by taking into account Assumption (A7), the fact that $\lambda_n/\delta_n \rightarrow0$, and since $0< |c^{(0)}|<1$,  we  obtain 
\[
{\cal P}_{1,k,n} \leq \PP \bigg[ \bigg\{C \tau c^{(0)}(t^*_k - \hat t_k) \leq C n \lambda_n \bigg\} \cap \bigg\{ t^*_k -\hat t_k \geq n \delta_n \bigg\} \cap W_n \bigg] {\underset{n \rightarrow \infty}{\longrightarrow}} 0.
\] 
Alternatively, if $\omega_{t^*_k} +\omega_{\hat t_k}=\big(\max(b_n,d_n) \big)^{-\gamma}$, then by using condition \eqref{CD1}, and since $0< |c^{(0)}|<1$, we also get
\[
{\cal P}_{1,k,n} \leq \PP \bigg[ \bigg\{C \tau c^{(0)}(t^*_k - \hat t_k) \leq  n \lambda_n \big(\max(b_n,d_n) \big)^{-\gamma} \bigg\} \cap \bigg\{ t^*_k -\hat t_k \geq n \delta_n \bigg\} \cap W_n \bigg] {\underset{n \rightarrow \infty}{\longrightarrow}} 0.
\] 
Thus, the last relations imply that 
\begin{equation}
\label{PP1kn}
\lim_{n \rightarrow \infty}  {\cal P}_{1,k,n}=0.
\end{equation}

Now we deal with ${\cal P}_{2,k,n}$.  Using the fact that $\|\ex\| \geq \|\ex\|_\infty$, we immediately have
\begin{eqnarray}
{\cal P}_{2,k,n} & =& \PP \bigg[ \bigg\{\bigg\|\sum^{t^*_k-1}_{i=\hat{t}_k} \ex_i \e1_{\{Y_i \leq \ex_i^\top \widehat{\eb}_{i}\}}\bigg\| \geq  \frac{v-1}{v} \bigg\|\tau \sum^{t^*_k-1}_{i=\hat t_k} \ex_i\bigg\| \bigg\} \cap V_{n,k} \cap W_n \bigg] \nonumber \\
& \leq &  \PP \bigg[ \bigg\{\frac{v-1}{v} \tau c^{(0)} (t^*_k- \hat t_k) \leq \bigg\|\sum^{t^*_k-1}_{i=\hat t_k} \ex_i \e1_{\{Y_i \leq \ex_i^\top \widehat{\eb}_{i}\}} \bigg\|  \bigg\} \cap V_{n,k} \cap W_n \bigg] 
. \label{P1}
\end{eqnarray}
Now, by using the convergence rate $(c_n)$ given by \eqref{cn}, Assumption (A1), and the fact that $\max_{i\in \{1, \cdots ,n\}} \|\ex_i\| \leq c^{(1)}$,  we can write
\begin{align*}
\bigg\|\sum^{t^*_k-1}_{i=\hat t_k} \ex_i \e1_{\{Y_i \leq \ex_i^\top \widehat{\eb}_{i}\}} \bigg\|  & \leq \sum^{t^*_k-1}_{i=\hat t_k} \| \ex_i \|\e1_{\{Y_i \leq \ex_i^\top \widehat{\eb}_{i}\}} \leq  c^{(1)} \sum^{t^*_k-1}_{i=\hat t_k} \e1_{\{\varepsilon_i \leq \ex_i^\top (\widehat{\ef}_{k+1} -\ef^*_k)\}}\\
& = c^{(1)} \sum^{t^*_k-1}_{i=\hat t_k} \e1_{\{\varepsilon_i - \ex_i^\top(\ef^*_{k+1} -\ef^*_k)\leq \ex_i^\top (\widehat{\ef}_{k+1} -\ef^*_{k+1})\}} \\ & \leq c^{(1)} \sum^{t^*_k-1}_{i=\hat t_k} \e1_{\{\varepsilon_i - \ex_i^\top(\ef^*_{k+1} -\ef^*_k)\leq | \ex_i^\top (\widehat{\ef}_{k+1} -\ef^*_{k+1})|\}}\\
& \leq c^{(1)} \sum^{t^*_k-1}_{i=\hat t_k} \e1_{\{\varepsilon_i - \ex_i^\top(\ef^*_{k+1} -\ef^*_k)\leq \| \ex_i\| \cdot \| \widehat{\ef}_{k+1} -\ef^*_{k+1}\|\}} \\ & \leq c^{(1)} \sum^{t^*_k-1}_{i=\hat t_k} \e1_{\{\varepsilon_i - \ex_i^\top(\ef^*_{k+1} -\ef^*_k)\leq c^{(1)} c_n\}},
\end{align*}
all with probability one, except the last inequality which holds with probability converging to one as $n \rightarrow \infty$. Thus,  for \eqref{P1} we obtain that
\begin{equation}
\label{tot1}
{\cal P}_{2,k,n} \leq \PP \bigg[ \bigg\{\frac{v-1}{v} \tau c^{(0)} (t^*_k- \hat t_k) \leq  c^{(1)} \sum^{t^*_k-1}_{i=\hat t_k} \e1_{\{\varepsilon_i - \ex_i^\top(\ef^*_{k+1} -\ef^*_k)\leq c^{(1)} c_n\}}  \bigg\} \cap V_{n,k} \cap W_n \bigg] +o(1).
\end{equation}
Now, by Lemma \ref{Lemma 4bis}, since $c_n \rightarrow \infty$, we have
\[
\PP \bigg[\bigg\{ \bigg| \frac{1}{ t^*_k - \hat t_k} \sum^{t^*_k-1}_{i=\hat t_k}\e1_{\{\varepsilon_i - \ex_i^\top(\ef^*_{k+1} -\ef^*_k)\leq c^{(1)} c_n\}} -L_k \bigg| \geq \bigg| \frac{v-1}{v} \frac{\tau}{c^{(1)}} c^{(0)}-L_k \bigg| \bigg\}   \cap V_{n,k}\bigg] {\underset{n \rightarrow \infty}{\longrightarrow}} 0,
\]
and since $(v-1) \tau c^{(0)}/ \big(v c^{(1)}  \big)-L_k >0$, we also have
\[
\PP \bigg[\bigg\{  \frac{1}{ t^*_k - \hat t_k} \sum^{t^*_k-1}_{i=\hat t_k}\e1_{\{\varepsilon_i - \ex_i^\top(\ef^*_{k+1} -\ef^*_k)\leq c^{(1)} c_n\}} -L_k   \geq \frac{v-1}{v} \frac{\tau}{c^{(1)}} c^{(0)} -L_k \bigg\}  \cap V_{n,k}  \bigg] {\underset{n \rightarrow \infty}{\longrightarrow}} 0.
\]
Therefore, taking into account the relations in \eqref{P1} and \eqref{tot1} we again conclude that
\[
{\cal P}_{2,k,n}  {\underset{n \rightarrow \infty}{\longrightarrow}} 0,
\]
which proves the case for $\tau >L_k$. Now, we prove an analogous result for the situation where $\tau < L_k$. In such case the probabilities ${\cal P}_{1,k,n}$ and ${\cal P}_{2,k,n}$ can be expressed as 
\[
{\cal P}_{1,k,n} = \PP \bigg[ \bigg\{ \frac{1}{v}\bigg\|\sum^{t^*_k-1}_{i=\hat t_k} \ex_i \e1_{\{Y_i \leq \ex_i^\top \widehat{\eb}_{i}\}} \bigg\| \leq 2n \lambda_n\big(\omega_{t^*_k} +\omega_{\hat t_k}\big) \bigg\} \cap V_{n,k} \cap W_n \bigg],
\]
and
\[
{\cal P}_{2,k,n}  = \PP \bigg[ \bigg\{\tau\bigg\|\sum^{t^*_k-1}_{i=\hat t_k} \ex_i \bigg\| \geq  \frac{v-1}{v} \bigg\| \sum^{t^*_k-1}_{i=\hat t_k} \ex_i \e1_{\{Y_i \leq \ex_i^\top \widehat{\eb}_{i}\}}\bigg\| \bigg\} \cap V_{n,k} \cap W_n \bigg].
\]
Similarly as before, we study ${\cal P}_{1,k,n}$  and ${\cal P}_{2,k,n}$ separately. 
Firstly, for ${\cal P}_{1,k,n}$, we can use the fact that $\| \ex \| \geq \| \ex\|_\infty$  and, also, that the first component of $\ex_i$ is $c^{(0)}$, to obtain 
\begin{eqnarray}
{\cal P}_{1,k,n} & \leq & \PP \bigg[ \bigg\{ \bigg\|\sum^{t^*_k-1}_{i=\hat t_k} \ex_i \e1_{\{Y_i \leq \ex_i^\top \widehat{\eb}_{i}\}} \bigg\|_\infty \leq 2vn \lambda_n\big(\omega_{t^*_k} +\omega_{\hat t_k}\big) \bigg\} \cap V_{n,k} \cap W_n \bigg] \nonumber \\
& \leq & \PP \bigg[ \bigg\{ c^{(0)}\sum^{t^*_k-1}_{i=\hat t_k}  \e1_{\{\varepsilon_i \leq \ex_i^\top ( \widehat{\ef}_{k+1}- \ef^*_k)\}}  \leq 2v n \lambda_n\big(\omega_{t^*_k} +\omega_{\hat t_k}\big) \bigg\} \cap V_{n,k} \cap W_n \bigg] .
\label{eqp1t}
\end{eqnarray}
Now, since $\| \widehat{\ef}_{k+1}- \ef^*_k \| =O_{\PP}(\| \ef^*_{k+1}+c_n C- \ef^*_k \|)$, taking into account Assumptions (A1) and (A5), and since $c_n {\underset{n \rightarrow \infty}{\longrightarrow}} 0$, we have that
$$\ex_i^\top ( \widehat{\ef}_{k+1}- \ef^*_k) \leq \big| \ex_i^\top ( \widehat{\ef}_{k+1}- \ef^*_k)\big| \leq \|\ex_i\| \cdot \| \widehat{\ef}_{k+1}- \ef^*_k \| \leq c^{(1)} C,$$
which holds with probability converging to one.  Thus, with probability converging to one, we also obtain that
\[
\sum^{t^*_k-1}_{i=\hat t_k}  \e1_{\{\varepsilon_i \leq \ex_i^\top ( \widehat{\ef}_{k+1}- \ef^*_k)\}}  \geq \sum^{t^*_k-1}_{i=\hat t_k}  \e1_{\{\varepsilon_i \leq - |\ex_i^\top ( \widehat{\ef}_{k+1}- \ef^*_k)|\}}  \geq \sum^{t^*_k-1}_{i=\hat t_k}  \e1_{\{\varepsilon_i \leq - c^{(1)} C\}}. 
\]
Therefore, the relation in \eqref{eqp1t} can be further rewritten as 
\begin{eqnarray}
{\cal P}_{1,k,n} & \leq & \PP \bigg[ \bigg\{ c^{(0)}\sum^{t^*_k-1}_{i=\hat t_k}  \e1_{\{\varepsilon_i \leq - c^{(1)} C\}} \leq 2vn \lambda_n\big(\omega_{t^*_k} +\omega_{\hat t_k}\big) \bigg\} \cap V_{n,k} \cap W_n \bigg] +o(1)\nonumber \\
& \leq &  \PP \bigg[ \bigg\{ c^{(0)}\sum^{t^*_k-1}_{i= t^*_k-[n \delta_n]}  \e1_{\{\varepsilon_i \leq - c^{(1)} C\}}  \leq 2vn \lambda_n\big(\omega_{t^*_k} +\omega_{\hat t_k}\big) \bigg\} \cap V_{n,k} \cap W_n \bigg] +o(1). \nonumber
\end{eqnarray}
Now, since we have 
$$
(n \delta_n)^{-1} \sum^{t^*_k-1}_{i= t^*_k-[n \delta_n]}  \e1_{\{\varepsilon_i \leq - c^{(1)} C\}} -F\big(-c^{(1)} C\big) \overset{a.s.} {\underset{n \rightarrow \infty}{\longrightarrow}}  0,
$$
where $F(-c^{(1)} C) > 0$ and $2vn \lambda_n\big(\omega_{t^*_k} +\omega_{\hat t_k}\big)=O_{\PP}(n \lambda_n)$, we finally obtain, due to Assumption (A7), that
$${\cal P}_{1,k,n} {\underset{n \rightarrow \infty}{\longrightarrow}} 0.$$ 

On the other hand, for ${\cal P}_{2,k,n}$ we obtain by Assumption (A1) that
\begin{eqnarray}
{\cal P}_{2,k,n} & \leq & \PP \bigg[ \bigg\{\tau\bigg\|\sum^{t^*_k-1}_{i=\hat t_k} \ex_i \bigg\| \geq  \frac{v-1}{v} \bigg\| \sum^{t^*_k-1}_{i=\hat t_k} \ex_i \e1_{\{\varepsilon_i \leq \ex_i^\top ( \widehat{\ef}_{k+1}- \ef^*_k)\}}\bigg\|_\infty \bigg\} \cap V_{n,k} \cap W_n \bigg] \nonumber \\
& \leq & \PP \bigg[ \bigg\{\tau\bigg\|\sum^{t^*_k-1}_{i=\hat t_k} \ex_i \bigg\| \geq  \frac{v-1}{v} c^{(0)} \sum^{t^*_k-1}_{i=\hat t_k}   \e1_{\{\varepsilon_i \leq \ex_i^\top ( \widehat{\ef}_{k+1}- \ef^*_k)\}}  \bigg\} \cap V_{n,k} \cap W_n \bigg] \nonumber \\
& \leq & \PP \bigg[ \bigg\{\tau\sum^{t^*_k-1}_{i=\hat t_k} \|\ex_i \| \geq  \frac{v-1}{v} c^{(0)} \sum^{t^*_k-1}_{i=\hat t_k}   \e1_{\{\varepsilon_i \leq \ex_i^\top ( \widehat{\ef}_{k+1}- \ef^*_k)\}}  \bigg\} \cap V_{n,k} \cap W_n \bigg] \nonumber \\
& \leq & \PP \bigg[ \bigg\{\tau\sum^{t^*_k-1}_{i=\hat t_k} \max_i \|\ex_i \| \geq  \frac{v-1}{v} c^{(0)} \sum^{t^*_k-1}_{i=\hat t_k}   \e1_{\{\varepsilon_i \leq \ex_i^\top ( \widehat{\ef}_{k+1}- \ef^*_k)\}}  \bigg\} \cap V_{n,k} \cap W_n \bigg] \nonumber \\
& \leq & \PP \bigg[ \bigg\{  \frac{1}{t^*_k - \hat t_k}   \sum^{t^*_k-1}_{i=\hat t_k}   \e1_{\{\varepsilon_i \leq \ex_i^\top ( \widehat{\ef}_{k+1}- \ef^*_k)\}}  \leq \tau \frac{c^{(1)}}{c^{(0)}} \frac{v}{v-1} \bigg\} \cap V_{n,k} \cap W_n \bigg] .
\label{rrr}
\end{eqnarray}
By Lemma \ref{Lemma 4bis},  we have
\[
\PP \bigg[\bigg\{ \bigg| \frac{1}{ t^*_k - \hat t_k} \sum^{t^*_k-1}_{i=\hat t_k}\e1_{\{\varepsilon_i - \ex_i^\top(\ef^*_{k+1} -\ef^*_k)\leq c^{(1)} c_n\}} -L_k \bigg| \geq \bigg| \tau \frac{c^{(1)}}{c^{(0)}} \frac{v}{v-1}-L_k \bigg| \bigg\}  \cap V_{n,k} \bigg] {\underset{n \rightarrow \infty}{\longrightarrow}} 0,
\]
and since $v \geq \frac{c^{(0)} L_k}{c^{(0)} L_k - \tau c^{(1)}}$, then $\tau \frac{c^{(1)}}{c^{(0)}} \frac{v}{v-1}-L_k < 0$. Thus, the relation in \eqref{rrr} implies that
\[
\PP \bigg[ \bigg\{  \frac{1}{t^*_k - \hat t_k}   \sum^{t^*_k-1}_{i=\widehat t_k}   \e1_{\{\varepsilon_i \leq \ex_i^\top ( \widehat{\ef}_{k+1}- \ef^*_k)\}} -L_k  \leq \tau \frac{c^{(1)}}{c^{(0)}} \frac{v}{v-1} -L_k \bigg\} \cap V_{n,k} \cap W_n \bigg] {\underset{n \rightarrow \infty}{\longrightarrow}} 0,
\]
and, therefore, we can again conclude that 
$${\cal P}_{2,k,n} {\underset{n \rightarrow \infty}{\longrightarrow}} 0.$$ 

\noindent\textit{(ii)} We now study $\sum^{K^*}_{k=1}\PP[V_{n,k} \cap \overline W_n]$, with the random event $\overline W_n \equiv \big\{ \max_{1 \leqslant k \leqslant k^*} |\hat t_k -t^*_k | \geq I^*_{min/2} \big\}$. Let us define the following random events:
\begin{eqnarray}
D^{(1)}_n & \equiv & \big\{\forall k \in \{1, \cdots , K^* \}; t^*_{k-1} <\hat t_k \leq t^*_{k+1} \big\} \cap \overline W_n , \nonumber \\
D^{(2)}_n & \equiv & \big\{\exists k \in \{1, \cdots , K^* \}; \hat t_k \leq t^*_{k-1} \big\} \cap \overline W_n, \nonumber \\
D^{(3)}_n & \equiv & \big\{\exists k \in \{1, \cdots , K^* \}; \hat t_k \geq t^*_{k+1} \big\} \cap \overline W_n .\nonumber 
\end{eqnarray}
Considering the random events defined above, we can use the decomposition
\begin{equation}
\label{mm_eq2}
\sum^{K^*}_{k=1}\PP[V_{n,k} \cap \overline W_n]=\sum^{K^*}_{k=1} \PP[V_{n,k} \cap D^{(1)}_n]+\sum^{K^*}_{k=1} \PP[V_{n,k} \cap D^{(2)}_n]+\sum^{K^*}_{k=1} \PP[V_{n,k} \cap D^{(3)}_n],
\end{equation}
and we will study each term on the right side separately.  Let us start with $\sum^{K^*}_{k=1} \PP[V_{n,k} \cap D^{(1)}_n]$. We easily obtain that
\begin{align*}
\PP[V_{n,k} \cap D^{(1)}_n]&=\PP\bigg[V_{n,k} \cap D^{(1)}_n \cap \bigg\{\hat t_{k+1}-t^*_k  \geq \frac{I^*_{min}}{2} \bigg\}\bigg]+\PP\bigg[V_{n,k} \cap D^{(1)}_n \cap \bigg\{\hat t_{k+1}-t^*_k < \frac{I^*_{min}}{2} \bigg\}\bigg] \\
& \leq \PP\bigg[V_{n,k} \cap D^{(1)}_n \cap \bigg\{\hat t_{k+1}-t^*_k  \geq \frac{I^*_{min}}{2} \bigg\}\bigg]+\PP\bigg[V_{n,k} \cap D^{(1)}_n \cap \bigg\{t^*_{k+1}-\hat t_{k+1} \geq  \frac{I^*_{min}}{2} \bigg\}\bigg] ,
\end{align*}
where we used the fact that $0 \leq \hat t_{k+1}-t^*_k \leq I^*_{min}/2$ implies $$t^*_{k+1}-\hat t_{k+1}=(t^*_{k+1}- t^*_k)-( \hat t_{k+1} - t^*_k) \geq I^*_{min} -I^*_{min}/2=I^*_{min}/2.$$ Hence, we get
\begin{align}
\label{B4}
\sum^{K^*}_{k=1} \PP[V_{n,k} \cap D^{(1)}_n]& \leq \sum^{K^*}_{k=1}\PP\bigg[V_{n,k} \cap D^{(1)}_n \cap \bigg\{\hat t_{k+1}-t^*_k  \geq \frac{I^*_{min}}{2} \bigg\}\bigg] \\
&\quad +\sum^{K^*}_{k=1} \sum^{K^*-1}_{l=k+1} \PP \bigg[\bigg\{t^*_l-\hat t_l >  \frac{I^*_{min}}{2} \bigg\} \cap \bigg\{ \hat t_{l+1}- t^*_l \geq  \frac{I^*_{min}}{2} \bigg\} \cap D^{(1)}_n \bigg]. \nonumber 
\end{align}
To deal with the first additive term on the right-hand side of \eqref{B4} we firstly need to find the convergence rate of the regression parameter estimator. For this term we are in the following situation:
\[
t^*_{k-1} < \hat t_k < t^*_k - n \delta_n < t^*_k < t^*_k +I^*_{min}/2 < \hat t_{k+1} < t^*_{k+1}.
\]
Thus, we apply Lemma \ref{aKKT} for $j=t^*_k - [n \delta_n]$ and  $j=t^*_k$, and we obtain, with probability one, that
\begin{align}
2 n \lambda_n \big(\omega_{t^*_k}+\omega_{t^*_k-[n \delta_n]} \big)&  \geq  \bigg\| \tau \sum^{t^*_k -1}_{i=t^*_k - [n \delta_n]} \ex_i - \sum^{t^*_k -1}_{i=t^*_k - [n \delta_n]} \ex_i \e1_{\{Y_i \leq \ex_i^\top \widehat{\eb}_i\}}  \bigg\| \nonumber \\
&  \geq  \bigg\|  \sum^{t^*_k -1}_{i=t^*_k -[ n \delta_n]} \ex_i \e1_{\{\varepsilon_i \leq \ex_i^\top (\widehat{\ef}_{k+1}- \ef^*_k)\}}  \bigg\| - \tau   \bigg\|\sum^{t^*_k -1}_{i=t^*_k - [n \delta_n]} \ex_i  \bigg\| .\label{EVk}
\end{align}
Hence, with probability one, we also have 
\[
\bigg\|  \sum^{t^*_k -1}_{i=t^*_k - [n \delta_n]} \ex_i \e1_{\{\varepsilon_i \leq \ex_i^\top (\widehat{\ef}_{k+1}- \ef^*_k)\}}  \bigg\| \leq 2 n \lambda_n\big(\omega_{t^*_k}+\omega_{t^*_k-[n \delta_n]} \big)  +\tau  \bigg\|\sum^{t^*_k -1}_{i=t^*_k - [n \delta_n]} \ex_i  \bigg\| .
\]
As before, we again obtain in a similar way, that with probability equal to one, it holds that
\begin{align}
\tau  \bigg\|\sum^{t^*_k -1}_{i=t^*_k - [n \delta_n]} \ex_i  \bigg\| -2  n \lambda_n \big(\omega_{t^*_k}+\omega_{t^*_k-[n \delta_n]} \big) & \leq \bigg\|  \sum^{t^*_k -1}_{i=t^*_k - [n \delta_n]} \ex_i \e1_{\{\varepsilon_i \leq \ex_i^\top (\widehat{\ef}_{k+1}- \ef^*_k)\}}  \bigg\| \nonumber \\ & \leq \tau  \bigg\|\sum^{t^*_k -1}_{i=t^*_k -[ n \delta_n]} \ex_i  \bigg\| + 2 n \lambda_n \big(\omega_{t^*_k}+\omega_{t^*_k-[n \delta_n]} \big) . \label{B5}
\end{align}
Now, similarly as for \eqref{A4}, we also  obtain  that 
$$\sum^{t^*_k -1}_{i=t^*_k -[ n \delta_n]} \ex_i \big[ \e1_{\{\varepsilon_i \leq \ex^\top_i(\widehat{\ef}_{k+1} -\ef^*_k)\}}- \tau\big]=O_{\PP}\bigg( n \lambda_n \big(\max(b_n,d_n) \big)^{-\gamma}\bigg) \eu,
$$
where $\eu \in \R^p$, such that $\|\eu \|=1$. Thus,  by Assumption (A4) and the condition in  \eqref{CD1}, we conclude that 
\begin{equation}
\label{F3}
\big\| \widehat{\ef}_{k+1} -\ef^*_k\big\|=O_{\PP}\left( \frac{\lambda_n \big(\max(b_n,d_n) \big)^{-\gamma}}{\delta_n} +(n \delta_n)^{-1/2}\right)  =o_{\PP}(1).
\end{equation}
In the same way as above, applying Lemma \ref{aKKT} for $j=t^*_k +I^*_{min}/2$ and $j=t^*_k$, we have with probability 1 that
\label{B6}\begin{align}
\tau  \bigg\|\sum^{t^*_k +I^*_{min}/2 -1}_{i=t^*_k} \ex_i  \bigg\| - 2 n \lambda_n \big(\omega_{t^*_k}+\omega_{t^*_k+I^*_{min}/2} \big)& \leq \bigg\|  \sum^{t^*_k +I^*_{min}/2 -1}_{i=t^*_k} \ex_i \e1_{\{\varepsilon_i \leq \ex_i^\top (-\widehat{\ef}_{k+1}+ \ef^*_{k+1})\}}  \bigg\|  \\
& \leq \tau  \bigg\|\sum^{t^*_k +I^*_{min}/2 -1}_{i=t^*_k} \ex_i  \bigg\| +  2 n \lambda_n \big(\omega_{t^*_k}+\omega_{t^*_k+I^*_{min}/2} \big),\nonumber
\end{align}
and analogously with \eqref{F3}, we obtain that
\begin{equation}
\label{F4}
\big\| \widehat{\ef}_{k+1} - \ef^*_{k+1}\big\| =O_{\PP}\left( \frac{\lambda_n \big(\max(b_n,d_n) \big)^{-\gamma}}{\delta_n} +(n \delta_n)^{-1/2}\right) =o_{\PP}(1).
\end{equation}
However, the relations in \eqref{F3} and \eqref{F4} are in a contradiction because by Assumption  (A5)\textit{(b)} we have $\| \ef^*_{k+1} - \ef^*_k \| >c^{(b)}>0$. 
Therefore, taking into account the relation in \eqref{EVk}, we have
\[  
\lim_{n \rightarrow \infty}\PP \bigg[  V_{n,k} \cap D^{(1)}_n \cap \bigg\{\hat t_{k+1}-t^*_k  \geq \frac{I^*_{min}}{2} \bigg\}    \bigg]  =0,
\]
which also implies that the limit of the first term on the right-hand side of \eqref{B4} is  
\[\lim_{n \rightarrow \infty}\sum^{K^*}_{k=1}\PP\bigg[V_{n,k} \cap D^{(1)}_n \cap \bigg\{\hat t_{k+1}-t^*_k  \geq \frac{I^*_{min}}{2} \bigg\}\bigg] =0 .
\]

In a similar manner we can also proceed with for the second terms in \eqref{B4}.  Applying Lemma \ref{aKKT} for $j = t^*_l -I^*_{min}/2$ and $j = t^*_l$ and, afterwards, for  $j = t^*_l$ and  $j = t^*_l +I^*_{min}/2$, we obtain the contradiction. Therefore, we conclude that 
$$\lim_{n \rightarrow \infty} \sum^{K^*}_{k=1} \PP[V_{n,k} \cap D^{(1)}_n] =0.$$

For the second term in \eqref{mm_eq2} we can write
\begin{align}
\PP[V_{n,k} \cap D_n^{(2)}] & \leq \PP[  D_n^{(2)}] \leq  \sum^{K^*}_{j=1} 2^{j-1} \PP \bigg[\max \big\{l \in \{ 1, \cdots  , K^*\} ; \hat t_l \leq t^*_{l-1}\big\} =j \bigg]\nonumber\\
& \left.\begin{array}{r}
\leq K^* \sum^{K^*-1}_{j=1}2^{j-1} \PP \bigg[ \bigg\{ t^*_j - \hat t_j \geq \frac{I^*_{min}}{2} \bigg\} \cap \bigg\{\hat t_{j+1} - t^*_j \geq \frac{I^*_{min}}{2}   \bigg\}   \bigg]\\[0.4cm]
+  2^{K^* -1} \cdot K^* \cdot \PP \bigg[t^*_{K^*} - \hat t_{K^*} \geq \frac{I^*_{min}}{2} \bigg]
\end{array}\right\} \label{B9}
\end{align}
Applying again  Lemma \ref{aKKT} for $j = t^*_{K^*}$ and $j = \hat t^*_{K^*}$, we obtain that
\begin{align}
\PP & \bigg[t^*_{K^*} - \hat t_{K^*} \geq \frac{I^*_{min}}{2} \bigg] \\ & =  \PP \bigg[\bigg\{  n \lambda_n \big( \omega_{t^*_{K^*}}+ \omega_{\hat t_{K^*}} \big) \geq \bigg\| \tau \sum^{t^*_{K^*}}_{i=\hat t_{K^*}} \ex_i - \sum^{t^*_{K^*}}_{i=\hat t_{K^*}} \ex_i  \e1_{\{Y_i \leq \ex_i^\top (\widehat{\ef}_{K+1}- \ef^*_K)\}} \bigg\|  \bigg\} \cap \bigg\{ t^*_{K^*} - \hat t_{K^*} \geq \frac{I^*_{min}}{2} \bigg\}  \bigg].\nonumber
\end{align} 

Similarly as before, we again apply  Lemma \ref{aKKT} for $j = t^*_{K^*}- I^*_{min}/2$ and $j = t^*_{K^*}$, to show that \eqref{B9} converges to zero, as $n  \rightarrow \infty $.  For the first term in \eqref{B9} we have
\begin{align*}
K^* \sum^{K^*-1}_{j=1} & 2^{j-1} \PP \bigg[ \bigg\{ t^*_j - \hat t_j \geq \frac{I^*_{min}}{2} \bigg\} \cap \bigg\{\hat t_{j+1} - t^*_j \geq \frac{I^*_{min}}{2}   \bigg\}   \bigg] \\
& \leq  K^* 2^{K^* -1} \sum^{K^*-1}_{j=1}  \PP \bigg[ E_j \cap \bigg\{ t^*_j - \hat t_j \geq \frac{I^*_{min}}{2} \bigg\} \cap \bigg\{\hat t_{j+1} - t^*_j \geq \frac{I^*_{min}}{2}   \bigg\}   \bigg]\\
& \leq  K^* 2^{K^* -1} \sum^{K^*-1}_{j=1} \bigg( \PP \bigg[ \bigg\{n \lambda_n \big( \omega_{t^*_{K^*}}+\omega_{t^*_{K^*}- I^*_{min}/2}  \big) \geq \tau \bigg\| \sum^{t^*_j}_{i=\hat t_j} \ex_i \bigg\| \bigg\} \cap \bigg\{t^*_j -  \hat t_j \geq \frac{I^*_{min}}{2} \bigg\} \bigg]\\
& + \PP \left[  \bigg\{n \lambda_n \big( \omega_{t^*_{K^*}}+\omega_{t^*_{K^*}- I^*_{min}/2}  \big) \geq  \| \sum^{t^*_j}_{i=\hat t_j} \ex_i \e1_{\{\varepsilon_i \leq \ex_i^\top(\widehat{\ef}_{j+1} - \ef^*_j)\}} \| \bigg\}\cap \bigg\{t^*_j -  \hat t_j \geq \frac{I^*_{min}}{2} \bigg\} \right.\\
& \hskip10cm \left.\cap \bigg\{\hat t_{j+1} - t^*_j \geq \frac{I^*_{min}}{2}   \bigg\}  \right]  \bigg),
\end{align*}
and we will show that both probability terms in the last relation converges to zero for $n \to \infty$.  

For the first probability term we can use Assumption (A4) and the condition in \eqref{CD1}, to obtain 
\begin{align*}
\lim_{n \rightarrow \infty} \PP & \bigg[ \bigg\{n \lambda_n  \big( \omega_{t^*_{K^*}}+\omega_{t^*_{K^*}- I^*_{min}/2}  \big)\geq \tau \bigg\| \sum^{t^*_j}_{i=\hat t_j} \ex_i \bigg\| \bigg\} \cap \bigg\{t^*_j -  \hat t_j \geq \frac{I^*_{min}}{2} \bigg\} \bigg]\\
& \leq \lim_{n \rightarrow \infty} \PP \bigg[\tau c^{(0)} \frac{n \delta_n}{2}  \leq \tau c^{(0)} \frac{I^*_{min}}{2} \leq \tau c^{(0)}(t^*_j - \hat t_j) \leq  \big(\max(b_n,d_n) \big)^{-\gamma} n \lambda_n \bigg] = 0,
\end{align*}
while the second probability term can be showed to converge to 0  by applying Lemma \ref{aKKT} for $j = t^*_j - I^*_{min}/2$ and  $j = t^*_j$ and, afterwards, for $j = t^*_j$ and  $j = t^*_{j+1}$, and showing that $\|\widehat{\ef}_{j+1}-  \ef^*_j \| {\underset{n \rightarrow \infty}{\longrightarrow}} 0$, and  $\|\widehat{\ef}_{j+1}-  \ef^*_{j+1} \| {\underset{n \rightarrow \infty}{\longrightarrow}} 0$, in probability, which is a contradiction with Assumption (A5)\textit{(b)}.  Therefore, we conclude that also $\lim_{n \rightarrow \infty} \PP[V_{n,k} \cap D_n^{(2)}]=0$.

Following the same lines as above, it can be also shown that $\lim_{n \rightarrow \infty}\PP[V_{n,k} \cap D_n^{(3)}] = 0$ and therefore, we can conclude that 
\begin{equation}
\label{v2}
\lim_{n \rightarrow \infty} \sum^{K^*}_{k=1}\PP[V_{n,k} \cap \overline W_n]=0,
\end{equation}
which finally completes the proof of the theorem.  \hspace*{\fill}$\blacksquare$ \\

\noindent {\bf Proof of Theorem  \ref{aTheorem 3.1ii}}\\  
By Theorem \ref{aTheorem 3.1} we have, for any $k =1, \cdots, K^*$, that $| \hat t_k - t^*_k| =O_{\PP}(n \delta_n)$, which is also $o_{\PP}(I^*_{min})$ since $(n \delta_n)^{-1}I^*_{min} \rightarrow \infty$. Thus, for any $k \in \{1, \cdots , K^*\}$ we  either have $(t^*_{k-1}+t_k^*)/2 < \hat t_k < t^*_k$ or $t^*_k \leq \hat t_k <(t^*_k+t^*_{k+1})/2$. We suppose that $(t^*_{k-1}+t_k^*)/2 < \hat t_k < t^*_k$ and  for $\hat t_{k+1} $, we again either get  $(t^*_k+t^*_{k+1})/2 < \hat t_{k+1} < t^*_{k+1}$ or   $t^*_{k+1} \leq \hat t_{k+1} < (t^*_k+t^*_{k+1})/2$. For $(t^*_k+t^*_{k+1})/2 < \hat t_{k+1} < t^*_{k+1}$ we apply  Lemma \ref{aKKT}(ii) for $j=t^*_k$ and $j=(t^*_k+t^*_{k+1})/2$, and using the inequality in \eqref{e4}, we obtain, with probability equal to one, that 
\begin{align*}
2 n \lambda_n \big(  \omega_{t^*_k}+\omega_{(t^*_k+t^*_{k+1})/2} \big) & \geq  \bigg\| \tau \sum^{(t^*_k+t^*_{k+1})/2}_{i=t^*_k}\ex_i - \sum^{(t^*_k+t^*_{k+1})/2}_{i=t^*_k}\ex_i \e1_{\{Y_i \leq \ex_i^\top \widehat{\eb}_i\}} \bigg\|\\ &   \geq  \tau  \bigg\| \sum^{(t^*_k+t^*_{k+1})/2}_{i=t^*_k}\ex_i \bigg\|-  \bigg\|\sum^{(t^*_k+t^*_{k+1})/2}_{i=t^*_k}\ex_i \e1_{\{\varepsilon_i \leq \ex^\top_i (\widehat{\ef}_{k+1}- \ef^*_{k+1})\}} \bigg\| .
\end{align*}
Now, similarly as in the proof of Theorem \ref{aTheorem 3.1}, we show that with probability equal to one it holds that 
\[
2  n \lambda_n \big(  \omega_{t^*_k}+\omega_{(t^*_k+t^*_{k+1})/2} \big)  \geq  \bigg\|\sum^{(t^*_k+t^*_{k+1})/2}_{i=t^*_k}\ex_i \e1_{\{\varepsilon_i \leq \ex^\top_i (\widehat{\ef}_{k+1}- \ef^*_{k+1})\}} \bigg\| - \tau  \bigg\| \sum^{(t^*_k+t^*_{k+1})/2}_{i=t^*_k}\ex_i \bigg\|,
\]
which also implies
\begin{align*}
\tau  \bigg\| \sum^{(t^*_k+t^*_{k+1})/2}_{i=t^*_k}\ex_i \bigg\| - 2 n \lambda_n \big(  \omega_{t^*_k}+\omega_{(t^*_k+t^*_{k+1})/2} \big) &  \leq  \bigg\|\sum^{(t^*_k+t^*_{k+1})/2}_{i=t^*_k}\ex_i \e1_{\{\varepsilon_i \leq \ex^\top_i (\widehat{\ef}_{k+1}- \ef^*_{k+1})\}} \bigg\| \\ &  \leq  \tau  \bigg\| \sum^{(t^*_k+t^*_{k+1})/2}_{i=t^*_k}\ex_i \bigg\| + 2 n \lambda_n  \big(  \omega_{t^*_k}+\omega_{(t^*_k+t^*_{k+1})/2} \big) .
\end{align*}
By Remark \ref{Rq_vit}, we have that $\omega_{t^*_k}+\omega_{(t^*_k+t^*_{k+1})/2}=O_{\PP}\big((\max (b_n,d_n))^{-\gamma}  \big)$ and by the condition in \eqref{CD2} also
\begin{align*}
\bigg\| \sum^{(t^*_k+t^*_{k+1})/2}_{i=t^*_k}\ex_i \bigg\| \geq  \bigg\| \sum^{(t^*_k+t^*_{k+1})/2}_{i=t^*_k}\ex_i \bigg\|_\infty & \geq | c^{(0)}|\big((t^*_k+t^*_{k+1})/2 - t^*_k \big) = | c^{(0)}|(t^*_{k+1}-t^*_k)/2\\
& \geq | c^{(0)}| I^*_{min}/2 \gg n \delta_n \geq n \lambda_n  (\max (b_n,d_n))^{-\gamma}.
\end{align*}
Therefore, with probability converging to one, we obtain that  
\[
\sum^{(t^*_k+t^*_{k+1})/2}_{i=t^*_k}\ex_i \e1_{\{\varepsilon_i \leq \ex^\top_i (\widehat{\ef}_{k+1}- \ef^*_{k+1})\}}  =\tau \sum^{(t^*_k+t^*_{k+1})/2}_{i=t^*_k}\ex_i  +O_{\PP}\bigg(n \lambda_n\big(\max(b_n,d_n) \big)^{-\gamma} \bigg) \eu,
\]
for any  $\eu \in \R^p$, such that $\| \eu \| =1$. Thus, similarly as  in the proof of Theorem  \ref{aTheorem 3.1} we have
\[
\big\|  \widehat{\ef}_{k+1} - \ef^*_{k+1} \big\|=O_{\PP} \left(\frac{n \lambda_n \big(\max(b_n,d_n) \big)^{-\gamma}+\sqrt{I^*_{min}}}{I^*_{min}} \right) \ll O_{\PP} \left(\frac{n \lambda_n \big(\max(b_n,d_n) \big)^{-\gamma} +\sqrt{n \delta_n}}{n \delta_n} \right).
\]
On the other hand, if $t^*_{k+1} \leq \hat t_{k+1} < (t^*_k+t^*_{k+1})/2$, we again apply Lemma \ref{aKKT}, however, for $j=t^*_k$ and $t^*_{k+1}$, and by repeating the same arguments as above we conclude that the assertion of the theorem holds true.
\hspace*{\fill}$\blacksquare$ \\

\noindent {\bf Proof of Theorem  \ref{aTheorem 3.2}}\\  
Since  $ |\widehat{{\cal A}}_n|  \leq  K_{max} < \infty $, we have that
\begin{align}
\label{ea0}
\PP \bigg[\big({\cal E}(\widehat{{\cal A}}_n|| {\cal A}^*) \geq n \delta_n \big) \cap \{ K^* \leq |\widehat{{\cal A}}_n| \leq K_{max}  \} \bigg] & \leq \PP \bigg[ {\cal E}(\widehat{{\cal A}}_n|| {\cal A}^*) \geq n \delta_n \bigg| \; | \widehat{{\cal A}}_n|=K^* \bigg]\\
& + \sum^{K_{max}}_{K=K^*+1} \PP \bigg[{\cal E}(\widehat{{\cal A}}_n|| {\cal A}^*) \geq n \delta_n  \bigg| \; |\widehat{{\cal A}}_n |=K  \bigg].\nonumber
\end{align}
For the first term of the right-hand side of  \eqref{ea0} we have, by Theorem \ref{aTheorem 3.1}, that
\begin{equation}
\label{ea1}
\lim_{n \rightarrow \infty}\PP \bigg[ {\cal E}(\widehat{{\cal A}}_n || {\cal A}^*) \geq n \delta_n \bigg| \; |\widehat{{\cal A}}_n|=K^* \bigg] = 0,
\end{equation}
and, for the second term in \eqref{ea0}, we can write
\begin{equation}
\label{eq25}
\sum^{K_{max}}_{K=K^*+1} \PP \bigg[{\cal E}(\widehat{{\cal A}}_n|| {\cal A}^*) \geq n \delta_n\bigg| \; |\widehat{{\cal A}}_n |=K  \bigg] \leq  \sum^{K_{max}}_{K=K^*+1} \sum^{K^*}_{k=1} \bigg\{\PP[E_{K,k,1}] +\PP[E_{K,k,2}] +\PP[E_{K,k,3}] \bigg\},
\end{equation}
with the random events $E_{K,k,1}$, $E_{K,k,2}$, and $E_{K,k,3}$ being defined as 
\begin{align*}
E_{K,k,1} & \equiv \{ \forall 1 \leq l \leq K; \,\, |\hat t_l - t^*_k| \geq n \delta_n , \, \hat t_l < t^*_k  \};\\
E_{K,k,2} & \equiv \{ \forall 1 \leq l \leq K; \,\, |\hat t_l - t^*_k| \geq n \delta_n , \, \hat t_l > t^*_k  \};\\
E_{K,k,3} & \equiv \{ \exists 1 \leq l \leq K; \,\, |\hat t_l - t^*_k| \geq n \delta_n , \, |\hat t_{l+1} - t^*_k| \geq n \delta_n , \, \hat t_l < t^*_k < \hat t_{l+1} \}.
\end{align*}
We start by studying the first probability, $\PP[E_{K,k,1}]$, for some $k \in \{1, \cdots , K^* \}$ and $K \in \{K^*+1, \cdots , K_{max} \}$. It holds that
\[
\PP[E_{K,k,1}]=\PP\big[E_{K,k,1} \cap\{\hat t_K> t^*_{k-1} \} \big] +\PP\big[E_{K,k,1} \cap\{\hat t_K \leq  t^*_{k-1} \} \big].
\]  
Let us consider the random event $E_{K,k,1} \cap\{\hat t_K> t^*_{k-1} \} $. In this case we have $t^*_{k-1} < \hat t_K < t^*_k -[n \delta_n] < t^*_k$. We apply Lemma \ref{aKKT} firstly for $j=t^*_k-[n \delta_n]$ and $j=t^*_k$ and, afterwards, for $j=t^*_k$ and  $j=t^*_k+n \delta_n$. Thus, we obtain as for relation \eqref{F3} and \eqref{F4} by using Assumption (A4) and the condition in \eqref{CD1}, that 
\[
\big\|\widehat{\ef}_{k+1} -\ef^*_k\big\|=O_{\PP}\left( \frac{\lambda_n \big(\max(b_n,d_n) \big)^{-\gamma}}{\delta_n} +(n \delta_n)^{-1/2}\right)  = o_{\PP}(1),
\]
and also
\[
\big\|\widehat{\ef}_{k+1} - \ef^*_{k+1}\big\|=O_{\PP}\left( \frac{\lambda_n\big(\max(b_n,d_n) \big)^{-\gamma}}{\delta_n} +(n \delta_n)^{-1/2}\right) = o_{\PP}(1),
\]
which contradicts Assumption (A5)\textit{(b)}. Therefore, we conclude that 
$$\lim_{n \rightarrow \infty} \PP\big[E_{K,k,1} \cap\{\hat t_K> t^*_{k-1} \} \big] = 0,$$ 
and, hence
\begin{equation}
\label{pe1}
\lim_{n \rightarrow \infty}\sum^{K_{max}}_{K=K^*+1} \sum^{K^*}_{k=1}\PP\big[E_{K,k,1} \cap\{\hat t_K> t^*_{k-1} \} \big] = 0.
\end{equation}

Next, we study the probability of the random event $E_{K,k,1} \cap\{\hat t_K < t^*_{k-1} \}$. For this random event we have $\hat t_K < t^*_{k-1} < t^*_k$. Again we can consider Lemma \ref{aKKT} firstly  for $j=t^*_{k-1}$ and $j=t^*_k$ and, afterwards, for $j=t^*_k$ and $j=t^*_{k+1}$. As above, we obtain a contradiction and, therefore, we conclude that 
\begin{equation}
\label{pe2}
\lim_{n \rightarrow \infty}  \sum^{K_{max}}_{K=K^*+1} \sum^{K^*}_{k=1}\PP\big[E_{K,k,1} \cap\{\hat t_K <  t^*_{k-1} \} \big]=0.
\end{equation}
Finally, the relations in \eqref{pe1} and \eqref{pe2} imply 
\begin{equation}
\label{eq37}
\lim_{n \rightarrow \infty}  \sum^{K_{max}}_{K=K^*+1} \sum^{K^*}_{k=1}\PP\big[E_{K,k,1}  \big]=0.
\end{equation}
In a similar way it can be also proved that 
$$\lim_{n \rightarrow \infty}  \sum^{K_{max}}_{K=K^*+1} \sum^{K^*}_{k=1}\PP\big[E_{K,k,2}  \big]=\lim_{n \rightarrow \infty}  \sum^{K_{max}}_{K=K^*+1} \sum^{K^*}_{k=1}\PP\big[E_{K,k,3}  \big]=0.$$ Therefore, the proof of the theorem follows by taking into account this last relation together with \eqref{eq37}, \eqref{eq25}, \eqref{ea1}, and  \eqref{ea0}.
\hspace*{\fill}$\blacksquare$ \\

\noindent {\bf Proof of Theorem  \ref{aTheorem 3.3}}\\
If $|\widehat{\cal A}_n | < K^*$, there are at least two true consecutive change-points without any change-point estimator in between them. Without any loss of generality we assume that $K^*=2$ and $|\widehat{\cal A}_n |=1$. The theorem is proved for this case if we show that $\lim_{n \rightarrow \infty } \PP\big[(K^*=2) \cap (|\widehat{\cal A}_n |=1) \big] =0$. Without loss of generality we can assume that 
$$1 < t^*_1 < \hat t_1 \leq t^*_2-I^*_{min}/2 < t^*_2 < n.$$
Thus, we apply Lemma \ref{aKKT} for $j=t^*_2-I^*_{min}/2$ and $j=t^*_2$ and we get that
\[
n \lambda_n \big(\omega_{t^*_2}+\omega_{t^*_2-I^*_{min}/2} \big) \geq \bigg\|\tau \sum^{t^*_2}_{i=t^*_2-I^*_{min}/2} \ex_i -  \sum^{t^*_2}_{i=t^*_2-I^*_{min}/2} \ex_i  \e1_{\{Y_i \leq \ex_i^\top \widehat{\eb}_i\}} \bigg\|,
\]
which holds with probability equal to one. By Assumptions (A4) and (A7) we also have  $I^*_{min} \geq n \delta_n \gg n \lambda_n$, and therefore, we obtain with probability converging to one, as $n \rightarrow \infty$, that 
\[
\sum^{t^*_2}_{i=t^*_2-I^*_{min}/2} \ex_i  \e1_{\{Y_i \leq \ex_i^\top (\widehat{\ef}_2- \ef^*_2)\}} =\tau  \sum^{t^*_2}_{i=t^*_2-I^*_{min}/2} \ex_i \pm C n \lambda_n \big(\max(b_n,d_n) \big)^{-\gamma} ,
\]
which implies (similarly as in the proof of Theorem  \ref{aTheorem 3.1}) together with Assumption (A4) and the condition in \eqref{CD1}, that
\begin{equation}
\label{F1}
\big\| \widehat{\ef}_2- \ef^*_2 \big\|=O_{\PP} \left(\frac{n \lambda_n \big(\max(b_n,d_n) \big)^{-\gamma} }{I^*_{min}} + \frac{1}{\sqrt{I^*_{min}}} \right) = o_{\PP}(1). 
\end{equation}

Next, we take Lemma \ref{aKKT} (ii) for $j=t^*_2$ and $j=t^*_2+ I^*_{min}/2$ to get
\begin{equation}
\label{F2}
\big\| \widehat{\ef}_2- \ef^*_3\big\| =O_{\PP} \left(\frac{n \lambda_n \big(\max(b_n,d_n) \big)^{-\gamma}}{I^*_{min}} + \frac{1}{\sqrt{I^*_{min}}} \right) = o_{\PP}(1). 
\end{equation}
However, both relations in \eqref{F1} and \eqref{F2} contradicts Assumption (A5)\textit{(b)}, therefore, we conclude that $$\PP\big[(K^*=2) \cap (|\widehat{\cal A}_n |=1) \big] {\underset{n \rightarrow \infty}{\longrightarrow}} 0,$$ which  completes the proof. \hspace*{\fill}$\blacksquare$ \\

\noindent {\bf Proof of Corollary  \ref{Raq_vit}}\\
By Theorem  \ref{aTheorem 3.3} we have with probability converging to one that $|\widehat{\cal A}_n | \geq |{\cal A}^* |$. Since $\hat t_j -\hat t_{j-1} \overset{\PP} {\underset{n \rightarrow \infty}{\longrightarrow}} \infty$, for any $ j \in \{1, \cdots , |\widehat{{\cal A}}_n| +1\}$,  we have, by Theorems \ref{aTheorem 3.1ii} and \ref{aTheorem 3.2},  that
\begin{equation}
\label{ii}
\sup_{\scriptsize \begin{array}{c}
	j \in \{1, \cdots,|\widehat{\cal A}_n |+1   \} \\
	|\hat{t}_j-t^*_k| \leq n \delta_n
	\end{array}} \| \widehat{\ef}_{j} -\ef^*_k \|=O_{\PP} (c_n ),
\end{equation}
for any $k=1, \cdots , K^*+1$, and the sequence $(c_n) $ defined in \eqref{cn}. 
Taking into account the relation in \eqref{ii}, we obtain that the adaptive fused quantile   estimator $\widehat{\eb^n}$ belongs to ${\cal V }_n(\eb^*)$ with probability converging to one, where 
$${\cal V }_n(\eb^*) \equiv \{\eb^n=(\eb_1^\top, \cdots, \eb_n^\top)^\top; \| \eb_i-\eb^*_i \| \leq C c_n, \forall i=1, \cdots ,n \},$$
for some $C>0$ large enough.  

We suppose that $|\widehat{\cal A}_n | >K^*$ and we consider the following set ${\cal W}_n \equiv \{\eb^n \in {\cal V }_n(\eb^*), \| \eb_i-\eb^*_i\| >0, \forall i  \in {\overline{{\cal A}^*}} \} $. Recall that $\overline{{\cal A}^*} =\big\{ i \in \{ 2, \cdots, n \}; \eb^*_i =\eb^*_{i-1} \big\}$. We will show that
\begin{equation}
\label{SP_eq31}
\lim_{n \rightarrow \infty} \PP [ \widehat{\eb^n}  \in {\cal W}_n]=0.
\end{equation}
For this, we consider two  vector parameters $\eb=(\eb_{{\cal A}^*}^\top,\eb_{\overline{{\cal A}^*}}^\top)^\top$ and $\eb^{(1)}=(\eb^{(1)\top}_{{\cal A}^*},\eb^{(1)\top}_{\overline{{\cal A}^*}})^\top$ such that $\eb_{{\cal A}^*}=\eb^{(1)}_{{\cal A}^*}$, the sub-vector $\eb^{(1)}_{\overline{{\cal A}^*} }$ containing the elements $\{ \eb_i^{(1)}, i \in \{2, \cdots , n\}; \; \eb_i^{(1)} =\eb_{i-1}^{(1)}\}$ and the sub-vector $\eb_{\overline{{\cal A}^*}}$ such that 
\begin{equation}
\label{KMC}
K_M \equiv Card \{ i \in \overline{{\cal A}^*} ; \eb_i \neq \eb_{i-1} \} \geq 1. 
\end{equation}
Then
\begin{eqnarray}
D_n (\eb, \eb^{(1)})   \equiv  {S}(\eb) -  {S}({\eb}^{(1)}) 
=   \sum^n_{i=1} \cro{\rho_\tau(Y_i - \ex^\top_i \eb_i)-\rho_\tau(Y_i - \ex^\top_i {\eb}^{(1)}_i)} +n \lambda_n \sum_{j \in \overline{\cal A}^*}  \omega_{j} \| \eb_j- \eb_{j-1}\|. \label{SP_DD}
\end{eqnarray}
For the first term of the right-hand side of \eqref{SP_DD} we use the following identity, which holds for any $x, y \in \R$ (see \cite{Knight-98}),  
\begin{equation*}
\label{Knight}
\rho_\tau(x-y)- \rho_\tau(x)=y(\e1_{\{x \leq 0\}} - \tau)+\int^y_0 (\e1_{\{x \leq t\}} -\e1_{\{x \leq 0\}})dt.
\end{equation*}
Hence, we obtain 
\begin{align}
\sum^n_{i=1} \cro{\rho_\tau(Y_i - \ex^\top_i \eb_i)-\rho_\tau(Y_i - \ex^\top_i {\eb}^{(1)}_i)} & = \sum_{i \in  \overline{{\cal A}^*}} \ex^\top_i\big( \eb_i  -\eb_i^{(1)}\big) \bigg( \e1_{\{Y_i -  \ex^\top_i {\eb}^{(1)}_i \leq 0\}}  - \tau \bigg)\nonumber\\
& + \sum_{i \in  \overline{{\cal A}^*}} \int^{\ex^\top_i\big( \eb_i  -\eb_i^{(1)}\big)}_0 \bigg[\e1_{\{Y_i -  \ex^\top_i {\eb}^{(1)}_i \leq  v\}} - \e1_{\{Y_i -  \ex^\top_i {\eb}^{(1)}_i \leq  0\}} \bigg] dv \nonumber\\
& \equiv T_{1n} + T_{2n}.\label{mm_eq1} 
\end{align}
Note, that by the condition in \eqref{CD2}, we have $c_n \sqrt{I^*_{max}}  \rightarrow \infty$.  Now we first study $T_{1n}$: for any $ i \in \overline{{\cal A}^*}$, we have $\eb^*_i-\eb_i^{(1)}=\eb^*_{i-1}-\eb_{i-1}^{(1)}=\textbf{0}_p$, which implies that  the expectation of $T_{1n}$ is
\[
\eE[T_{1n}]=\sum_{i \in  \overline{{\cal A}^*}}\big( \eb_i  -\eb_i^{(1)}\big)^\top \ex_i  \bigg(F\big( \ex^\top_i (\eb^{(1)}_i - \eb^*_i) \big) -F(0) \bigg) =0,
\]
and for variance, it holds that 
\[
Var[T_{1n}]=\sum_{i \in  \overline{{\cal A}^*}}\big( \ex^\top_i (\eb_i - \eb^{(1)}_i) \big)^2 \tau(1-\tau).
\]
Since $\eb, \eb^{(1)} \in {\cal V }_n(\eb^*)$, and taking into account Assumptions (A1) and (A6), we have the the variance is bounded, $ Var[T_{1n}] \leq O(c^2_n I^*_{max})$. On the other hand, since $|\widehat{\cal A}_n | \leq K_{max} < \infty$, we also have that $K_M$ is bounded. Thus, by the the Law of Large Numbers for independent random variables, we obtain that $T_{1n}= o_{\PP}\big(c_n \sqrt{I^*_{max}} \big)$. The second term in \eqref{mm_eq1} can be expressed as
\[
T_{2n}=\sum_{i \in  \overline{{\cal A}^*}} \int^{\ex^\top_i\big( \eb_i  -\eb_i^{(1)}\big)}_0 \bigg[ \e1_{\{\varepsilon_i <v\}}- \e1_{\{\varepsilon_i < 0\}}\bigg] dv,
\]
and its expectation, using also the Taylor expansion, is given by
\[
\eE[T_{2n}]=\sum_{i \in  \overline{{\cal A}^*}} \int^{\ex^\top_i\big( \eb_i  -\eb_i^{(1)}\big)}_0 \bigg[ F(v)-F(0)\bigg] dv= \sum_{i \in  \overline{{\cal A}^*}} \int^{\ex^\top_i\big( \eb_i  -\eb_i^{(1)}\big)}_0 \bigg[ v f(0)+ \frac{v^2}{2} f'(\tilde v)\bigg] dv,
\]
for $\tilde v \in (0, v)$. Since the derivative $f'$ is bounded in a some neighborhood of zero, taking also into account Assumption (A1) and the fact that $\| \eb_i-\eb^{(1)}_i \| \leq C c_n$, we have
\[
\eE[T_{2n}]=\bigg( \frac{f^2(0)}{2} \sum_{i \in  \overline{{\cal A}^*}} \big( \ex^\top_i\big( \eb_i  -\eb_i^{(1)}\big)\big)^2\bigg) (1+o(1))=O\big( c^2_n I^*_{max}\big)>0.
\]
For the variance of $  T_{2n}$, since $\varepsilon_i$ are independent, we get
\begin{align*}
Var[T_{2n}] & = \sum_{i \in  \overline{{\cal A}^*}} Var \bigg[ \int^{\ex^\top_i\big( \eb_i  -\eb_i^{(1)}\big)}_0 \big[ \e1_{\{\varepsilon_i <v\}}- \e1_{\{\varepsilon_i < 0\}}\big] dv  \bigg] \\
& = \sum_{i \in  \overline{{\cal A}^*}} \eE \bigg[ \int^{\ex^\top_i\big( \eb_i  -\eb_i^{(1)}\big)}_0 \bigg( \big[ \e1_{\{\varepsilon_i <v\}}- \e1_{\{\varepsilon_i < 0\}}\big] -\big[F(v)-F(0)  \big]  \bigg)dv  \bigg]^2\\
& \leq \sum_{i \in  \overline{{\cal A}^*}} \eE \bigg[ \bigg|\int^{\ex^\top_i\big( \eb_i  -\eb_i^{(1)}\big)}_0 \bigg(  \big[ \e1_{\{\varepsilon_i <v\}}- \e1_{\{\varepsilon_i < 0\}}\big] -\big[F(v)-F(0)  \big] \bigg) dv \bigg|  \bigg] \cdot 2 \big|\ex^\top_i(\eb_i -\eb_i^{(1)})   \big|\\
& \leq 2 \sum_{i \in  \overline{{\cal A}^*}} \int^{\ex^\top_i\big( \eb_i  -\eb_i^{(1)}\big)}_0 \big(F(v)-F(0) \big) dv \cdot 2 \max_i \| \ex_i \| \cdot \| \eb_i  -\eb_i^{(1)} \|, 
\end{align*}
and using Assumption (A1) we obtain that $Var[T_{2n}] \leq 4 \eE[T_{2n}] c^{(1)} c_n$.
Since $c_n \rightarrow 0$ as $n \rightarrow \infty$, by Bienaym\'e-Tchebychev inequality we have $T_{2n}=O_{\PP}\big( c_n^2 I^*_{max} \big)$ and moreover $T_{2n} \geq C >0$ with probability converging to 1 as $n \rightarrow \infty$. Therefore, for the first term on the right-hand side of \eqref{SP_DD} we finally get 
\begin{equation}
\label{T12n}
\sum^n_{i=1} \cro{\rho_\tau(Y_i - \ex^\top_i \eb_i)-\rho_\tau(Y_i - \ex^\top_i {\eb}^{(1)}_i)} = T_{2n}+o_{\PP}(T_{2n})=O_{\PP}\big( c_n^2 I^*_{max} \big) >0.
\end{equation}
Now we study the second  term on the right-hand side of \eqref{SP_DD}. 
For any $j \in\overline{{\cal A}^*}$, we have $\eb^*_j =\eb^*_{j-1}$, $\eb_j - \eb_{j-1}=O(c_n)$. Using Remark \ref{Rq_vit}, we also have $\overset{\smile}{ \eb}_j -\overset{\smile}{ \eb}_{j-1}=O_{\PP}(b_n)$. Hence, since $K_M$ is bounded, it follows that 
\begin{equation}
\label{PenS}
0< n \lambda_n \sum_{j \in \overline{{\cal A}^*}}  \omega_{j} \| \eb_j- \eb_{j-1}\|= n \lambda_n \sum_{j \in \overline{{\cal A}^*}} \frac{c_n}{\big( \max(b_n,d_n)\big)^\gamma}\geq  n \lambda_n \frac{c_n }{\big( \max(b_n,d_n)\big)^\gamma},
\end{equation}
and by taking into account \eqref{SP_DD}, \eqref{T12n}, and \eqref{PenS} we obtain that
\begin{equation}
\label{D11}
D_n (\eb, \eb^{(1)})\geq O_{\PP}\big( c_n^2 I^*_{max} \big) + n \lambda_n \frac{c_n }{\big( \max(b_n,d_n)\big)^\gamma}=O_{\PP}\big( c_n^2 I^*_{max} \big)>0,
\end{equation}
which implies, together with the fact that $D(\eb^{(1)}, \eb^{(1)})=0$, the  relation in \eqref{SP_eq31} and also the corollary which follows. \hspace*{\fill}$\blacksquare$

\end{appendices}


\begin{thebibliography}{99}
	\bibliographystyle{abbrv}


	
	\bibitem[Ciuperca(2016)]{Ciuperca-16}
	Ciuperca, G. (2016). 
	Adaptive LASSO model selection in a multiphase quantile regression. 
	{\it Statistics}, {\bf 50}(5), 1100--1131. 
	

	
	\bibitem[Ciuperca and Maciak(2018)]{Ciuperca-Maciak-17}
	Ciuperca, G. and Maciak, M. (2018). 
	Change-point Detection by the  Quantile LASSO Method. 
	{\it Canadian Journal of Statistics}, (submitted). 
	
		\bibitem[Harchaoui and Levy(2010)]{Harchaoui.Levy.10}
	Harchaoui, Z. and L\'evy-Leduc, C. (2010).
	Multiple change-point estimation with a total variation penalty.
	{\it Journal of the American Statistical Association}, {\bf 105}(492), 1480--1493.
	
	\bibitem[Hyun et al.(2018)]{Hyun-GSell-Tibshirani-18}
	Hyun, S. and G'Sell, M. and Tibshirani, R.J. (2018)
	Exact post-selection inference for the generalized lasso path. 
	{\it Electronic Journal of Statistics}, {\bf 12}, 1053--1097. 	
	
	\bibitem[Jin et al.(2016)]{Jin-Wu-Shi-16}
   Jin, B. and Wu, Y. and Shi, X. (2016)
	Consistent two-stage multiple change-point detection in linear models. 
	{\it Canadian Journal of Statistics}, {\bf 44(2)}, 161--179. 	
		
	\bibitem[Knight(1998)]{Knight-98}
	Knight,  K. (1998). 
	Limiting distributions for L1 regression estimators under general conditions.
	{\it Annals of Statistics},  {\bf 26}(2),  755--770.
	
	\bibitem[Lee et al.(2018)]{Lee-Liao-Seo-Shin-18}
	Lee, S. and Liao, Y. and Seo M.H.  and Shin, Y. (2018)
	Oracle estimation of a change point in high-dimensional quantile regression. 
	{\it Journal of the American Statistical Association}, {\bf 113}(523), 1184--1194.	

	\bibitem[Leonardi and Buhlmann(2016)]{Leonardi-Buhlmann.16}
	Leonardi, F. and Buhlmann, P. (2016).
	Computationally efficient change point detection for high-dimensional regression.
	\textit{arxiv:1601.03704}.


	\bibitem[Lin et al.(2016)]{lin2016}
	Lin K., Sharpnack J., Rinaldo A., and  Tibshirani R. J.(2016).
	Approximate Recovery in Changepoint Problems, from $l_2$ Estimation Error Rates.
	{\it arxiv.org/abs/1606.06746}
	
	\bibitem[Qian and Su(2016)]{Qian.Su.16}
	Qian, J. and Su, L. (2016).
	Shrinkage estimation of regression models with multiple structural changes.
	{\it Econometric Theory}, {\bf 32}(6), 376--1433.

	\bibitem[Rinaldo(2009)]{Rinaldo.09}
	Rinaldo,  A. (1998).
	Properties and refinements of the fused Lasso.
	{\itshape Annals of Statistics}, {\bf 37(5B)}, 2922--2952.
	
	\bibitem[Simon et al.(2012)]{simon}
	Simon, N., Friedman, J., Hastie, T., and Tibshirani, R. (2012).
	A Sparse-Group Lasso.
	{\itshape Journal of Computational and Graphical Statistics }, {\bf 22(2)}, 231--245.
	

	\bibitem[Zhang and Xiang(2015)]{Zhang-Xiang-15}
	Zhang, C. and Xiang, Y. (2015)
	On the oracle property of adaptive group LASSO in high-dimensional linear models. 
	{\it Statistical Papers}, {\bf 57(1)}, 249--265.
	
	\bibitem[Zhang and Geng(2015)]{Zhang.Geng.15}
	Zhang, B. and Geng, J. (2015)
	Multiple change-points estimation in linear regression models via sparse group lasso. 
	{\it IEEE Transactions on Signal Processing}, {\bf 63(9)}, 2209--2224.
	
	

	\bibitem[Zheng et al.(2013)]{Zheng-Gallagher-Kulasekera-13}
	Zheng, Q. and Gallagher, C. and Kulasekera, K.B. (2013)
	Adaptive penalized quantile regression for high dimensional data. 
	{\it Journal of Statistical Planning Inference}, {\bf 143}, 1029--1038.



	\bibitem[Zheng et al.(2015)]{Zheng-Peng-He-15}
	Zheng, Q. and Peng, L. and He, X. (2015)
	Globally adaptive quantile regression with ultra-high dimensional data. 
	{\it Annals of Statistics}, {\bf 43(5)}, 2225--2258.


\end{thebibliography}
\end{document}